\documentclass[12pt, a4paper,twoside]{amsart}
\usepackage{a4wide}
\usepackage{floatflt,subfigure}
\usepackage{version}
\usepackage{ifthen}
\usepackage{graphicx}

\newcommand\proofbox{\ensuremath{\blacksquare}\relax}
\newcommand{\noproofbox}{\ensuremath{\square}\relax}
\newcommand\singlebox{\hskip10000pt minus 1fil}
\newcommand\esinglebox{\hskip10000pt minus 1fil \llap{\proofbox}}

\usepackage{amsmath}
\usepackage{amssymb}

\makeatletter
\def\swappedhead#1#2#3{%
  \thmnumber{\@upn{\the\thm@headfont #2\@ifnotempty{#1}{.~}}}%
  \thmname{#1}%
  \thmnote{ {\the\thm@notefont(#3)}}}
\makeatother

 \newtheoremstyle{changebreak}
   {9pt}
   {9pt}
   {\itshape}
   {}
   {\bfseries}
   {.}
   {\newline}
   {}

\theoremstyle{changebreak}
\swapnumbers

\newtheorem{thm}{Theorem}[section]
\newtheorem{defn}[thm]{Definition}
\newtheorem{lem}[thm]{Lemma}
\newtheorem{cor}[thm]{Corollary}

\newtheorem{prop}[thm]{Proposition}

\newtheorem{observation}[thm]{Observation}

\newtheorem{thmdef}[thm]{Theorem and Definition}
\newtheorem{deflem}[thm]{Definition and Lemma}

\swapnumbers

\newtheorem{algorithm}[thm]{Algorithm}

\def\proof{\par\smallskip
             \noindent {\sc Proof. }}
\def\proofof #1 {\par\medskip\noindent {\sc Proof of #1. }}

\def\sketch{\par\medskip\noindent{\sc Sketch of proof. }}
\def\sketchof #1 {\par\medskip\noindent {\sc Sketch of proof of #1. }}

\def\qed{\rule{0pt}{0pt}\nolinebreak\hfill\proofbox \par\medskip}
\def\qedd{\rule{0pt}{0pt}\nolinebreak\hfill\noproofbox}

\newenvironment{remark}[1][]{%
  \par\noindent \textsc{\ifthenelse{\equal{#1}{}}{Remark. }{#1. }}}{%
   \par\medskip}

\renewcommand{\mod}{\operatorname{mod}\,}

\usepackage[refpage]{nomencl}
\usepackage{version}
\excludeversion{suppress}
\includeversion{parameterarguments}

\usepackage{hyperref}

\makenomenclature

\makeindex
\renewcommand{\nomname}{List of Symbols}
\renewcommand{\pagedeclaration}[1]{\quad \dotfill\ \quad 
                      \makebox[2em][r]{#1}}

\newcommand{\Z}{\mathbb{Z}}
\newcommand{\Q}{\mathbb{Q}}
\newcommand{\R}{\mathbb{R}}

\newcommand{\C}{\ensuremath{\mathbb{C}}}
\newcommand{\Ch}{\hat{\mathbb{C}}}

\newcommand{\D}{\mathbb{D}}
\newcommand{\Ds}{\D^*}

\renewcommand{\H}{\mathbb{H}}

\newcommand{\onehalf}{\ensuremath{\frac{1}{2}}}

\newcommand{\dt}[1]{\mathcal{#1}}

\newcommand{\re}{\operatorname{Re}}
\newcommand{\im}{\operatorname{Im}}

\newcommand{\cl}[1]{\overline{#1}}
\newcommand{\interior}{\operatorname{int}}
\newcommand{\dist}{\operatorname{dist}}

\renewcommand{\theta}{\vartheta}
\renewcommand{\phi}{\varphi}
\renewcommand{\rho}{\varrho}

\newcommand{\Ek}{E_{\kappa}}

\newcommand{\ul}[1]{\underline{#1}}

\newcommand{\Orb}{\dt{O}}

\newcommand{\Sequ}{\dt{S}}

\newcommand{\Sequb}{\overline{\Sequ}}

\newcommand{\adds}{\underline{s}}

\newcommand{\s}{\adds}

\renewcommand{\u}{{\tt u\hspace{.25pt}}}
\newcommand{\m}{{\tt m\hspace{1.2pt}}}

\newcommand{\bdyit}[2]
             {{\rule{0pt}{0pt}_{\mbox{$\scriptstyle #2$}}^{\mbox{%
                   $\scriptstyle #1$}} }}
\renewcommand{\j}{{\tt j}}
\newcommand{\itj}{\bdyit{\j}{\j-1}}

\newcommand{\rt}{\tilde{r}}
\newcommand{\addrt}{\ul{\tilde{r}}}
                               
\newcommand{\addu}{\ul{\u}}
\newcommand{\ut}{\tilde{\u}}
\newcommand{\addut}{\ul{\ut}}

\newcommand{\addt}{\ul{t}}
\renewcommand{\t}{\addt}

\renewcommand{\r}{\ul{r}}
\renewcommand{\rt}{\widetilde{\r}}

\newcommand{\extaddr}{\operatorname{addr}}
\newcommand{\bifaddr}{\operatorname{addr}}
\newcommand{\ADDR}{\operatorname{Addr}}
\newcommand{\itin}{\operatorname{itin}}
\newcommand{\K}{\mathbb{K}}
\newcommand{\KP}{\K^{+}}
\newcommand{\KM}{\K^{-}}
\newcommand{\KS}{\K^{*}}
\newcommand{\Sec}{\operatorname{Sec}}
\newcommand{\Hyp}[1]{\operatorname{Hyp}(#1)}
\newcommand{\Bif}[1]{\operatorname{Bif}(#1)}

\newcommand{\gs}{g_{\adds}}

\newcommand{\periodic}[1]{\overline{#1}}
\newcommand{\per}[1]{\periodic{#1}}

\newcommand{\W}{\dt{W}}

\newcommand{\intaddnext}{\mapsto}

\newcommand{\wt}[1]{\widetilde{#1}}

\newcommand{\Hplane}{\mathcal{H}}

\newcommand{\U}{\mathcal{U}}

\newcommand{\M}{\mathcal{M}}

\newcommand{\IR}{\Gamma_{W,\hspace{.5pt}h}}
\newcommand{\IRH}[1]{\Gamma_{W,\hspace{.5pt}#1}}
\newcommand{\IRWH}[2]{\Gamma_{#1,\hspace{.5pt}#2}}

\newcommand{\Bdy}{\operatorname{Bdy}}

\newcommand{\child}{\Bif{W,h}}
\newcommand{\childH}[1]{\Bif{W,#1}}
\newcommand{\childWH}[2]{\Bif{#1,#2}}

\newcommand{\floor}[1]{\lfloor #1 \rfloor}
\newcommand{\ceil}[1]{\lceil #1 \rceil}

\newcommand{\picturedir}{.}

\title[Combinatorial Bifurcations of Exponential Maps]{%
  Combinatorics of Bifurcations \\in Exponential Parameter Space}
\author{Lasse Rempe}
\address{Mathematics Institute, University of Warwick, Coventry CV4 7AL,
United Kingdom}
\email{lasse@maths.warwick.ac.uk}
\date{today}
\thanks{The first author was supported in part
 by a postdoctoral fellowship of the 
 German Academic Exchange Service (DAAD) and
 by the German-Israeli Foundation
 for Scientific Research and Development (G.I.F.),
 grant no.\ G-643-117.6/1999}

\author{Dierk Schleicher}
\address{International University Bremen, P.O.~Box 750 561, 
28725 Bremen, Germany}
\email{dierk@iu-bremen.de}


\begin{document}

\begin{abstract}
 We give a complete combinatorial description of the
  bifurcation structure in the space of exponential maps
  $z\mapsto\exp(z)+\kappa$. This combinatorial structure is the
  basis for a number of important results about exponential parameter space.
  These include the fact that
  every hyperbolic component has connected boundary \cite{boundary,cras},
  a classification of escaping parameters \cite{markuslassedierk},
  and the fact that all dynamic and parameter
  rays at periodic addresses land \cite{landing2new,habil}.
\end{abstract}

\maketitle

\tableofcontents


\section{Introduction}

Ever since Douady and Hubbard's celebrated study of the Mandelbrot set
 \cite{orsay}, combinatorics has played a fundamental role for the
 dynamics of complex polynomials.
 In particular, the concept of \emph{external rays},
 both in the dynamical and parameter plane, and the
 landing behavior of such rays, 
 has helped in the understanding of polynomial
 Julia sets and bifurcation loci. This program has been particularly
 successful for the simplest polynomial parameter spaces: the quadratic
 family $z\mapsto z^2+c$ and its higher-dimensional cousins, 
 the unicritical families $z\mapsto z^d+c$
 \cite{orsay,fibers,eberleinschleicher}.

In this article, we consider the space of exponential maps,
  \[ \Ek:\C\to\C; z\mapsto \exp(z)+\kappa. \]
  It is well-known that a restriction on 
  the number of \emph{singular values} (i.e.,
  critical and asymptotic values) of
  an entire function 
  generally limits the amount of different dynamical features
  that can appear for the same map. Since exponential maps are
  the only transcendental entire functions which have just one
  singular value, namely the omitted value $\kappa$
  (see e.g.\ \cite[Appendix
  D]{jacktwocriticalpoints}), exponential maps
  form the  
  simplest parameter space of entire transcendental functions. 
  In addition, the exponential family can be considered as the
  limit of the polynomial unicritical families, and thus is an
  excellent candidate to apply the combinatorial methods which
  were so successful for Mandel- and Multibrot sets.


Recently, some
 progress has been made in this direction: 
 a complete classification of escaping points for exponential 
 maps in terms of
 \emph{dynamic rays} was given in
 \cite{expescaping}, and a similar construction was carried
 out to obtain \emph{parameter rays} 
 \cite{markus,markusdierk}. Also, exponential maps with attracting periodic
 cycles were classified in \cite{expattracting} using combinatorics.

 Nonetheless,
  a basic description of exponential dynamics in analogy to the
  initial study of the Mandelbrot set should involve at least
  the following results.
 \begin{enumerate}
  \item For every hyperbolic component $W$, there is a homeomorphism
   of pairs $(W,\cl{W})\to (\H,\cl{\H})$, where $\H$ is the left half
   plane. (In particular, $\partial W$ is a Jordan
     curve.) \label{item:hyperbolicbdy}
  \item Every periodic parameter ray lands at a parabolic parameter.
    \label{item:perparrays}
  \item If the singular value of $\Ek$ does not escape to $\infty$,
         then all periodic dynamic rays of $\Ek$ land.
    \label{item:landing2}
  \item If the singular value of $\Ek$ does not escape to $\infty$,
         then every repelling periodic point of $\Ek$ is the landing
         point of a periodic dynamic ray.
    \label{item:perpoints}
 \end{enumerate}

 For unicritical polynomials, the analogs of these statements all have
  relatively short analytic proofs (see e.g.\ 
  \cite{petersenryd} for (\ref{item:perparrays}) and
  \cite[Theorems 18.10 and 18.11]{jackdynamics} for (\ref{item:landing2}) and
   (\ref{item:perpoints})), but these break down
  in the exponential
  case. Nonetheless, it is possible to prove items
  (\ref{item:hyperbolicbdy}) through 
  (\ref{item:landing2}), using a novel
  approach based on a thorough study of parameter space.
  One of the goals of this
  article is to provide the first ingredients
  in this approach by obtaining a complete
  description of the combinatorial structure of
  parameter space (as given
  by bifurcations of hyperbolic components).
  In the sequel \cite{boundary}, this description is used
  to prove
  (\ref{item:hyperbolicbdy}), which, in turn, leads to proofs of
  (\ref{item:perparrays}) \cite{habil}
  and (\ref{item:landing2}) \cite{landing2new}, as well
  as some progress on
  (\ref{item:perpoints}) \cite{landing2new}.

\begin{figure}
  \begin{center}
   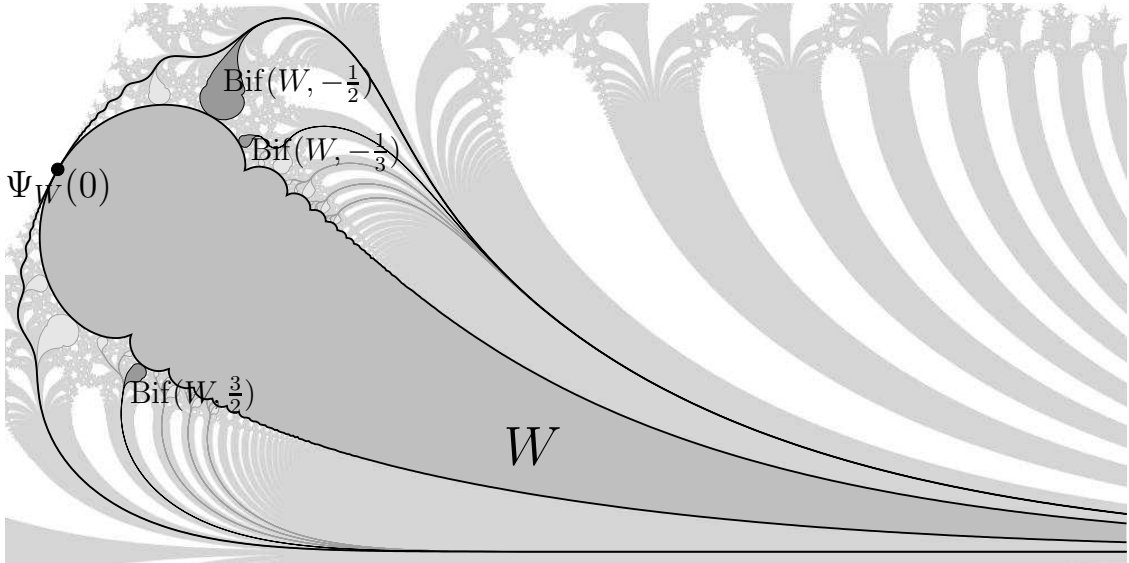
  \end{center}  
  \caption{\label{fig:childcomponents}%
    Structure of child components bifurcating from a 
    period $3$ hyperbolic component in exponential parameter space}
 \end{figure}

 To illustrate the difficulties we face, let us consider the structure
  of child components
  bifurcating from a given hyperbolic component. If we already knew
  results (\ref{item:hyperbolicbdy}) and
  (\ref{item:perparrays}) above, it would be quite easy
  to obtain the following description; compare
  Figure \ref{fig:childcomponents}.

 Let $W$ be a hyperbolic component of period $n\geq 2$, 
  and let $\mu:W\to \Ds$ be the
  multiplier map (which maps each hyperbolic parameter to the multiplier
  of its unique attracting cycle). Then there exists a conformal
  isomorphism $\Psi_W:\H\to W$ with $\mu\circ \Psi_W = \exp$ which
  extends continuously to $\partial \H$ and such that
  $\Psi_W(0)$ (the \emph{root} of $W$) is the landing point of two
  periodic parameter rays. The region containing $W$ which is enclosed by
  these two rays is called the \emph{wake} of $W$. 

 For every $h=\frac{p}{q}\in\Q\setminus \Z$, the point
  $\Psi_W(2\pi i h)$ is the root point of a (unique)
  hyperbolic component $\child$
  of period $qn$ (called a \emph{child component} of $W$).
  The component $\child$ tends to infinity above or below $W$
   depending on whether
  $h <0$ or $h>0$ (respectively). 
  If $0<h_1<h_2$ or $h_1<h_2<0$, then
  $\childH{h_1}$ tends to infinity
  below $\childH{h_2}$. Any hyperbolic component other
  than $W$
  which lies in the wake of $W$ is contained in the wake of a unique
  child component $\child$. 

 The problem we face 
  is that, without knowing (\ref{item:hyperbolicbdy}), 
  we do not know that all ``expected'' bifurcations really 
  exist. Thus we will need to obtain a purely 
  \emph{combinatorial} version 
  of this description (given by Theorem \ref{thm:bifurcationstructure}) 
  without being able to use the topological structure of parameter space. 
  This makes many arguments 
  (and statements) much more delicate. 

 Another goal of our article is to explain the relation among
  certain combinatorial objects which appear in exponential dynamics.
  In particular, there are several such objects associated to any
  hyperbolic component.
 \begin{itemize}
  \item The \emph{characteristic external addresses} of $W$
    (Definition \ref{thmdef:characteristicrays}). These are
    the addresses of the parameter rays bounding the wake of $W$.
  \item The \emph{intermediate external address} of $W$
    (Section \ref{sec:combinatorics}). This is an object which does not appear
    in Multibrot sets. It describes the combinatorial position
    of the singular value within the dynamical plane of a parameter
    in $W$. At the same time, it describes the position of the hyperbolic
    component $W$ itself in the vertical structure of parameter space.
  \item The \emph{kneading sequence} of $W$
    (Definition \ref{defn:combinatorialitinerary}). This object
    describes the itinerary of the singular orbit with respect to
    a natural \emph{dynamical} partition (as opposed to the
    \emph{static} partition used to define external 
    addresses).
  \item The \emph{internal address} of $W$ (Definition
    \ref{defn:internaladdress}). Introduced for Multibrot sets in 
    \cite{intaddr}, this address describes the
    position of $W$ within the bifurcation structure of hyperbolic
    components. Its relative, the
    \emph{angled internal address} (Definition 
    \ref{defn:angledinternaladdresses}), is decorated
    with some additional information.
 \end{itemize}
 Our study yields algorithms to convert between these different
  objects (where possible), and also to compute the address of any
  child component. These algorithms are collected in
  Appendix \ref{app:algorithms}. 

 Finally, the combinatorial objects and methods used in this article 
  have applications far beyond the scope of our 
  investigation, and are likely
  to play a significant role in further studies of the exponential family
  (as they did for the Mandelbrot set). Thus, we aim to present a 
  comprehensive exposition of these concepts which may serve
  as a reference in the future. 

\smallskip

 We should emphasize that all results of this article --- with the 
  exception of the analytical considerations of Section 
  \ref{sec:localbifurcation} --- 
  are completely combinatorial and could be formulated and 
  proved without any reference to the underlying exponential maps.  
  However, we prefer to carry out 
  an argument within an actual dynamical plane whenever possible, as we find 
  this much more intuitive. (Compare for example the 
  definition of orbit portraits in Section \ref{sec:orbitportraits},  
  as well as the 
  proofs of Lemma \ref{lem:uniqueaddress} and Proposition 
  \ref{prop:periodicintermediate}). 

 Since the combinatorial structure of exponential parameter space is
  a limit of that for unicritical polynomials, it would be possible to
  infer many of our results from corresponding facts for these families.
  However, many of these --- particularly for Multibrot sets of higher
  degrees --- are themselves still unpublished. Also, there are aspects
  of exponential dynamics, such as the \emph{intermediate external address}
  of a hyperbolic exponential map, 
  which would not feature in such an approach.
  We have thus decided to give a clean self-contained account in the 
  exponential case. 

\medskip\noindent
\textsc{Structure of this article.} 
 In the following two sections, we give a comprehensive overview of several
  combinatorial concepts for exponential maps: external addresses, 
  dynamic rays, intermediate external addresses, orbit portraits,
  characteristic rays and itineraries. In 
  Section \ref{sec:hyperboliccomponents}, we consider some basic facts
  about hyperbolic components of exponential maps, and how they are
  partitioned into \emph{sectors}.

 Section \ref{sec:localbifurcation} is the only part of
  the article in which analytical considerations are made:
  we investigate the stability of orbits at a parabolic point,
  allowing us to understand the structure of bifurcations
  occurring at such points. This provides 
  the link between our subsequent combinatorial
  considerations and the exponential parameter plane. 
  While the arguments in this section are very similar to those
  in the polynomial setting, there are
  some surprises: the combinatorics of
  a parent component can be determined with great ease from that of 
  a child component, thanks to the new feature of
  intermediate external addresses. 

 With these preliminaries, we will be in a position to 
  prove our main results.
  Section
  \ref{sec:bifurcationfromcomponent} deals with the structure
  of (combinatorial) child components of a given
  hyperbolic component, as discussed above. 
  Section 7 introduces introduces internal addresses, 
  giving a "human-readable" combinatorial structure to parameter space,
  and shows how they are related to
  the combinatorical concepts defined before.

 In Appendix \ref{app:furthertopics}, we consider some further concepts.
  These are not required for the proofs in \cite{boundary} but follow
  naturally from our discussion and will be collected for future reference. 
  Appendix \ref{app:algorithms} explicitly collects the combinatorial
  algorithms which are implied by our results.

 For the reader's convenience, a list of notation and an index of
  the relevant combinatorial concepts is provided at the end of the article.

\medskip\noindent
{\sc Some remarks on notation.} We have chosen to parametrize our
  exponential maps as $z\mapsto \Ek(z)=\exp(z)+\kappa$. Traditionally,
  they have often been parametrized as $\lambda\exp$, which is
  conjugate to $\Ek$ if $\lambda=\exp(\kappa)$. We prefer our
  parametrization mainly because the behavior of exponential maps at
  $\infty$, and in particular the asymptotics of external rays, do not
  depend on the parameter in this parametrization. Note that this is
  also the case in the usual parametrization of quadratic polynomials
  as $z\mapsto z^2+c$.
  Also, under our parametrization the picture
  in the parameter plane reflects the situation in the dynamical
  plane, which is a conceptual advantage. 
  Note that $\Ek$ and
  $E_{\kappa'}$ are conformally conjugate if and
  only if $\kappa-\kappa'\in 2\pi i \Z$. This will prove useful
  in the combinatorial description. 
  When citing known results,
  we always translate them into our parametrization.

 If $\gamma:[0,\infty)\to\C$ is a curve, we shall say that
  $\lim_{t\to\infty} \gamma(t)=+\infty$ (or, in short, call %
\index{curve to $\pm\infty$}%
  $\gamma$ a \emph{curve to $+\infty$}) if $\re
 \gamma(t)\to +\infty$ and $\im \gamma$ is bounded; analogously for
 $-\infty$. The $n$-th iterate of any function $f$ will be denoted by
 $f^n$. Whenever we write a rational number as a fraction
 $\frac{p}{q}$, we will
 assume $p$ and $q$ to be coprime. 

 We
 conclude any
 proof and any result which immediately follows from previously proved
 theorems by the symbol $\blacksquare$. A result which is cited without
 proof is concluded by $\square$.

\medskip\noindent
{\sc Acknowledgments.} We would like to thank Walter Bergweiler, Alex Eremenko,
 Markus F\"orster,
 Misha Lyubich, Jack Milnor, Rodrigo Perez, Phil Rippon and
 Juan Rivera-Letelier for many helpful
 discussions, and the Institute of Mathematical Sciences at Stony Brook
 as well as the University of Warwick 
 for continued support and hospitality.


\section{Combinatorics of Exponential Maps} \label{sec:combinatorics}

 An important combinatorial tool in the study of polynomials is the structure 
  provided by \emph{dynamic rays}, which foliate the basin of infinity. 
  Similarly, throughout this article, we will assign combinatorics to curves 
  in the dynamical plane of an exponential map, both in the set of escaping 
  points and in Fatou components. This section will review these methods, 
  which were introduced in \cite{expescaping} and \cite{expattracting}.

\subsection*{External addresses and dynamic rays} %
  \index{external address}%
  \index{external address!infinite}%
  \index{infinite external address}%
 A sequence $\s=s_1 s_2 \dots$ of integers is called an 
  (infinite) external address\footnote{%
 For brevity, we will frequently omit the adjective ``external''; 
  ``address'' will always mean ``external address'' unless explicitly
 stated otherwise}%
. If $s_1,\dots,s_n\in\Z$, then the
  address obtained by periodically repeating this sequence
  will be denoted by $\per{s_1\dots s_n}$.

 Let $\kappa\in\C$ and let 
  $\gamma:[0,\infty)\to\C$ be a curve in the dynamical plane of $\Ek$. 
  Then we say that $\gamma$ has \emph{external address} $\s$ if and only if %
   \index{external address!of a curve $\gamma$}%
   \[
      \lim_{t\to\infty} \re\Ek^{j-1}(\gamma(t))=+\infty \hspace{0.8cm} 
                                  \text{and}\hspace{0.8cm}
      \lim_{t\to\infty} \im\Ek^{j-1}(\gamma(t))=2\pi s_j \]
 for all $j\geq 1$; in this case, we also write
  $\s=\extaddr(\gamma)$. An external address 
\nomenclature[addrgamma]{$\extaddr(\gamma)$}{(external address of $\gamma$)}
  $\s$ is called \emph{exponentially bounded}
  if there exists some $x>0$ such that
  $2\pi |s_k| < F^{k-1}(x)$
  for all $k\geq 1$, where $F(t)=\exp(t)-1$ is used as a model function 
  for exponential growth.
\nomenclature[F]{$F(t)$}{(model for exponential growth)}

 The \emph{set of escaping points} of $\Ek$ is defined to be
  \[ I := I(\Ek) := \{z\in\C: |\Ek^n(z)|\to\infty\}. \]
\nomenclature[I]{$I(\Ek)$}{(set of escaping points)}
 It is known that the Julia set $J(\Ek)$ is the closure of $I(\Ek)$
\nomenclature[J]{$J(\Ek)$}{(Julia set)}
 \cite{alexescaping,alexmisha}. 
  In \cite{expescaping}, the set $I(\Ek)$ has been completely classified.
   In particular, it was shown that it consists of curves to $\infty$,
   so-called \emph{dynamic rays}. We will use this result in the
   following form. (Note that the fact that dynamic rays are 
   the path-connected components of $I(\Ek)$ was stated but
   not proved in
   \cite{expescaping}; for a proof compare \cite{markuslassedierk}.)

\begin{thmdef}[Dynamic Rays]
 Let $\kappa\in\C$. Then, for every exponentially bounded
  address $\s$, there exists a unique injective curve
  $\gs:[0,\infty)\to I(\Ek)$ or $\gs:(0,\infty)\to I(\Ek)$ which has 
  external address $\s$ and whose trace is a path-connected component
  of $I(\Ek)$. The curve $\gs$ 
  is called the \emph{dynamic ray} at address
  $\s$. %
 \index{dynamic ray}\nomenclature[gs]{$g_{\s}$}{(dynamic ray)}%

 If $\kappa\notin I(\Ek)$, then every path-connected component of $I(\Ek)$ 
  is a dynamic ray. If $\kappa\in I(\Ek)$, then every such component is either
  a dynamic ray or is mapped into a dynamic ray under finitely many
  iterations.
\end{thmdef}
\begin{remark} In order to state this theorem as given, dynamic rays
 need to be parametrized differently from \cite{expescaping}.
 In this article, we will only be using dynamic rays
 at periodic addresses $\s$, and for these our parametrization agrees with
 that of \cite{expescaping}, provided that the singular orbit does not 
 escape. 
\end{remark}

\subsection*{Intermediate external addresses}
 We shall also need to assign combinatorics to certain %
  \index{external address!intermediate}%
  \index{intermediate external address}%
  curves in Fatou components which, under finitely many iterations,
  map to a curve to $-\infty$. Let 
  $\gamma:[0,\infty)\to\C$ be a curve in the dynamical plane of 
  $\kappa$ such that, for some $n\geq 1$, 
   $\lim_{t\to\infty} \Ek^{n-1}(\gamma(t)) = -\infty$.
  Then there exist $s_1,\dots, s_{n-2}\in\Z$ and $s_{n-1}\in \Z+\onehalf$ such
   that
   \[ \lim_{t\to\infty} \re(\Ek^{j-1}(\gamma(t))) = +\infty
     \hspace{0.8cm} \text{and} \hspace{0.8cm}
      \lim_{t\to\infty} \im(\Ek^{j-1}(\gamma(t))) = 2\pi s_{j}
      \]
  for $j=1,\dots,n-1$. We call
   \begin{equation}
    \extaddr(\gamma) := s_1 s_2 \dots s_{n-1} \infty
      \label{eqn:intermediateaddr}
   \end{equation} %
   \index{external address!of a curve $\gamma$}%
   \index{intermediate external address!of a curve $\gamma$}%
  the \emph{intermediate external address} of $\gamma$. Any sequence of the
  form (\ref{eqn:intermediateaddr}) with $s_1,\dots, s_{n-2}\in\Z$ and
  $s_{n-1}\in\Z+\onehalf$ is called an intermediate external address
  (of length $n$).

 To illustrate
 the relationship between infinite and intermediate
 external addresses, 
 consider the following construction. Define
 \[ f:\R\setminus \bigl\{(2k-1)\pi : k\in\Z\bigr\}\to \R, t\mapsto \tan(t/2). \]
 Then to any (infinite) external address $\s$ we can associate a
 unique point $x$ for which $f^{k-1}(x)\in
 \bigl((2s_k-1)\pi,(2s_k+1)\pi\bigr)$ for all $k$. However, there 
 are countably many points which are not realized by any external address
 in this way, namely the
 preimages of $\infty$ under the
 iterates of $f$. Adding intermediate external addresses
 corresponds to filling in these points. The space $\Sequb$ of all
 infinite and intermediate external addresses is thus
 order-isomorphic to the circle $\overline{\R}\cong \mathbb{S}^1$. 
 We also set $\Sequ := \Sequb\setminus\{\infty\}$.
\nomenclature[S]{$\Sequ$, $\Sequb$}{(sequence spaces)}
   The \emph{shift map} is the function
   \[ \sigma: \Sequ\to\Sequb; s_1 s_2 \dots \mapsto s_2 \dots; \]
\nomenclature[sigma]{$\sigma$}{(shift map)}%
  note that $\sigma$ corresponds to the function $f$ in the above model.

\subsection*{Lexicographic and vertical order}
 The space $\Sequ$ naturally comes equipped with the 
  \emph{lexicographic order} on external addresses. (As seen above,
  this ordered space is
  isomorphic to the real line $\R$, and in particular is complete.) Similarly,
  the space $\Sequb$ carries a (complete) circular ordering.
  In our combinatorial considerations, we will routinely use the
  following fact.
  \begin{observation}[Shift Preserves Order On Small Intervals]
   \label{obs:orderpreserving}
   For every $\s=s_1 s_2 s_3\dots$ and
   $\s' := (s_1+1) s_2 s_3\dots$, the map
   $\sigma\colon [\s,\s')\to \Sequb$ preserves the circular order of
   $\Sequb$.
  \end{observation}

 Any family of pairwise disjoint curves to $+\infty$ 
  has a natural \emph{vertical order}:
  \index{vertical order}%
  among any two such curves,
  one is \emph{above} the other. More precisely, suppose that
  $\gamma:[0,\infty)\to\C$ is a curve to $+\infty$ and
  define $\Hplane_R
  := \{z\in\C: \re z > R\}$ for $R>0$. If $R$
  is large enough, then  the set $\Hplane_R \setminus \gamma$ has exactly two
  unbounded components,
  one above and one below $\gamma$. 
  Any curve $\wt{\gamma}$ to $+\infty$ which is 
  disjoint from $\gamma$
  must (eventually) tend to
  $\infty$ within one of these.

 It is an immediate consequence of the definitions that, if $\gamma$ 
   and $\wt{\gamma}$ have (infinite or intermediate) external addresses
   $\extaddr(\gamma)\neq \extaddr(\wt{\gamma})$, then $\gamma$ is above
   $\wt{\gamma}$ if and only if $\extaddr(\gamma) > \extaddr(\wt{\gamma})$.

\subsection*{Intermediate address of attracting and parabolic dynamics}
 Let us suppose that $\Ek$ 
   has an attracting or parabolic periodic point. 
  Then the singular value $\kappa$ is contained in
  some periodic Fatou component;  we call this component the 
   \emph{characteristic Fatou component}.
\index{characteristic Fatou component}%
\nomenclature{$U_1$}{(characteristic Fatou component)}%
  Let $U_0\mapsto U_1\mapsto \dots \mapsto U_n=U_0$ 
  be the cycle of periodic Fatou components, labeled such that
  $U_1$ is the characteristic component. (This will be our convention
  for the remainder of the paper.)
  Since $U_1$ contains a neighborhood of the singular value, $U_0$ contains a
  left half plane. In particular, $U_0$ contains a
  horizontal curve along which $\re(z)\to-\infty$.
  Its pullback to $U_1$ under $\Ek^{n-1}$ has an intermediate external address
  $\s$ of length $n$. (The address $\s$ does not depend on the initial
  choice of the curve to $-\infty$, since the latter is unique up to homotopy
  in $U_0$.)

  We call $\s$ the 
  \emph{intermediate external
  address of $\kappa$} and denote it by $\extaddr(\kappa)$;
  it will play a special role throughout this article. 
 \index{external address!of an attracting or parabolic map}%
 \index{intermediate external address!of an attracting or parabolic map}%
  \nomenclature[addrkappa]{$\extaddr(\kappa)$}{(intermediate 
                                       external address of $\kappa$)}%
  The following was proved independently in
   \cite{expattracting} and \cite{dfj}; the idea of the proof
   goes back to \cite[Section 7]{bakerexp}.

  \begin{prop}[Existence of Attracting Maps with Prescribed Combinatorics]
   \label{prop:hypexistence}
   Let $\s$ be an intermediate external address. Then there
    exists an attracting parameter $\kappa$ with
    $\extaddr(\kappa)=\s$. \qedd
  \end{prop}

  A converse result 
   was also proved in
   \cite{expattracting}: the external address $\extaddr(\kappa)$ 
   determines $\Ek$
   up to quasiconformal conjugacy
   (see Proposition \ref{prop:hyperboliccomponents}).

\subsection*{Attracting dynamic rays}

\begin{figure}
  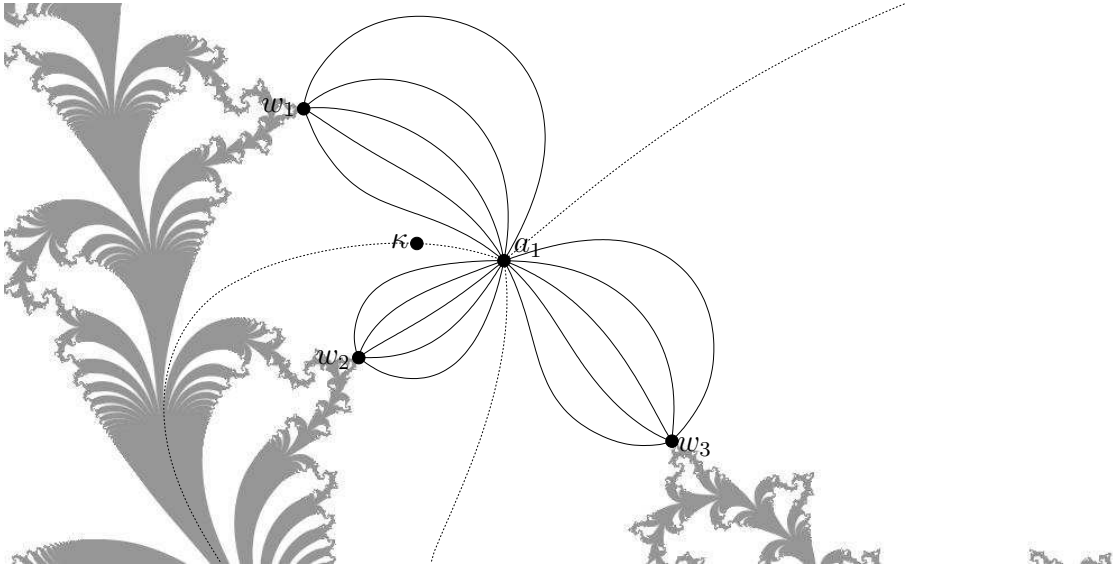
 \caption{\label{fig:attractingrays}%
   Attracting dynamic rays for a parameter where the
    attracting multiplier
    has angle $1/3$. All unbroken attracting rays land at the
    distinguished boundary orbit $\{w_i\}$; while the three
    broken attracting rays (shown as dotted lines) contain
    the singular orbit.}
\end{figure}

 We shall frequently have need for a canonical
  choice of certain curves in a Fatou component. 
  Let $\Ek$ have an attracting orbit, which we label $a_0\mapsto
  a_1\mapsto \dots \mapsto a_n=a_0$ such that $a_i\in U_i$. Note that
  we can connect $\kappa$ to $a_1$ by a straight line in linearizing
  coordinates. 
  The pullback of this curve under $\Ek^n$ along the
  orbit of $a_1$ is then a curve $\gamma\subset U_1$ which connects $a_1$ to
  $\infty$ and has $\extaddr(\gamma)=\extaddr(\kappa)$. We call $\gamma$
  the \emph{principal attracting ray} of $\Ek$. More
  \index{attracting dynamic ray!principal}%
  \index{principal attracting ray}%
  generally, any maximal
  curve in $U_1$ which starts at $a_1$ and is mapped into
  a radial line by the extended K{\oe}nigs map $\phi:U_1 \to \C$ is
  called an \emph{attracting dynamic ray} of $\Ek$. Apart from the
   \index{attracting dynamic ray}%
  principal attracting ray, the attracting dynamic ray which contains
  $\kappa$ will be particularly important. Those attracting dynamic
  rays which
  map to a curve to $-\infty$
  in $U_0$ under an iterate $\Ek^{nj-1}$, $j> 0$, will sometimes be called
  \emph{broken}; the principal attracting ray
  is one of these. Note that such a broken
  ray is mapped to a proper subpiece of the
  attracting ray containing $\kappa$ by $\Ek^{nj}$. 
   \index{attracting dynamic ray!broken}%

  It can be shown that every unbroken attracting
  dynamic ray lands at a point in $\C$ \cite[Theorem 4.2.7]{thesis}.
  We shall only need this fact in the case of
  rational multipliers, where it is simpler to prove  
  (see Figure \ref{fig:attractingrays}).
\begin{deflem}[Distinguished Boundary Orbit]
  \label{deflem:distinguishedboundaryorbit}
  \index{distinguished boundary orbit}%
 Consider an attracting exponential map $\Ek$ of period $n$ whose 
  attracting multiplier 
  has
  rational angle $\frac{p}{q}$. Then the orbit 
  of the principal attracting ray of $\Ek$
  under $\Ek^n$ consists of $q$ attracting dynamic rays and contains
  all points of the 
  singular orbit in $U_1$.

 Every other attracting dynamic ray starting at $a_1$ is periodic of
 period $q$ and lands at one of the points of a unique period $q$ orbit of
 $\Ek^n$ on $\partial U_1$.
 This orbit is called the 
  \emph{distinguished boundary orbit} on $\partial U_1$. 
\end{deflem}
\proof Analogous to the case $q=1$ \cite[Lemma 6.1]{expattracting}. \qedd


\section{Orbit Portraits and Itineraries}
 \label{sec:orbitportraits}

\subsection*{Orbit Portraits}
  Following Milnor \cite{jackrays}, we will use the notion of 
  \emph{orbit portraits} to encode the dynamics of periodic rays 
  landing 
  at common points. As in the case of quadratic polynomials, this is 
  important for understanding the structure of parameter space. 

 \begin{defn}[Orbit Portrait]
  \index{orbit portrait}%
   Let $\kappa\in\C$ and 
   let $(z_1,\dots,z_n)$ be a repelling or parabolic periodic orbit
   for $\Ek$.
  Define 
   \[ A_k:=\{\r\in\Sequ:
         \text{$\r$ is periodic and the dynamic ray $g_{\r}$ lands at $z_k$}\}. \]
   Then
  $\Orb := \{A_1,\dots,A_n\}$ is called the \emph{orbit portrait} of $(z_k)$. 
  The orbit (and the orbit portrait) is called \emph{essential} if $|A_k|>1$
   for any $k$. An essential orbit portrait
   is called of \emph{satellite type} if it
   contains only one cycle of rays; otherwise it is called \emph{primitive}.
    \index{orbit portrait!satellite}%
    \index{orbit portrait!primitive}%
 \end{defn}

 \begin{lem}[Basic Properties of Orbit Portraits]
   Let $\Orb=\{A_1,\dots,A_n\}$ be an orbit portrait. Then all $A_k$ are
   finite, and the shift map carries
    $A_k$ bijectively onto $A_{k+1}$.
   Furthermore, all addresses in the portrait share the same 
    period $qn$ (called the \emph{ray period} of the orbit) for
    some integer $q\geq 1$. \index{ray period!of an orbit portrait}%
 \end{lem}
 \proof  The proof that all rays share the
   same period and are transformed bijectively by the shift is
   completely analogous to the polynomial setting
   \cite[Lemma 18.12]{jackdynamics}.
  Let $qn$ be the common period of the rays in $\Orb$ and let 
   $\s = \per{s_1 s_2 \dots s_{qn}} \in A_1$. It is easy to see that
   \[ A_1 \subset \left\{\per{s_1' s_2'\dots s_{qn}'}:
                        |s_k'-s_k|\leq 1 \right\} \]
   (see e.g.\ \cite[Lemma 5.2]{topescapingnew}). 
    The set on the right 
    hand side is finite, as required.
   \qed

  An orbit portrait can also be defined
   as an abstract combinatorial object, without reference to any
   parameter. 
   We will not do this here, but 
   we will often suppress the
   actual choice of parameter present in its definition. 

 \subsection*{Characteristic rays}
  For quadratic polynomials, every orbit portrait has two distinguished
   rays, which are exactly the two rays which separate the critical
   value from all other rays in the portrait
   \cite[Lemma 2.6]{jackrays}. A corresponding
   statement for exponential maps is given by the following result.

\begin{deflem}[Characteristic Rays \protect{\cite[Lemma 5.2]{expattracting}}]
  \label{lem:permutation}
  \index{characteristic addresses/rays!of an orbit portrait}%
 Let $\kappa\in \C$.
  Suppose that $(z_k)$ is a repelling or parabolic periodic
  orbit with essential orbit portrait $\Orb$.
  Then there exist $\ell$ and two periodic rays $g_{\r}$ and
  $g_{\addrt}$ landing at $z_{\ell}$ (the \emph{characteristic rays of
  the orbit $(z_k)$})
  such that the curve
  $g_{\r}\cup\{z_{\ell}\}\cup g_{\addrt}$
   separates the singular value from all other rays of the orbit portrait.
  The addresses $\r$ and $\addrt$ are called the 
  \emph{characteristic addresses}
  of $\Orb$; they depend on $\Orb$ but not on $\kappa$. 
  The interval in $\Sequ$ bounded by $\r$ and $\addrt$ is called
  the \emph{characteristic sector} of $\Orb$. 
  \index{characteristic sector}%

  Furthermore, if there are at least three rays landing at each $z_k$,
  then the orbit portrait of $(z_k)$ is of satellite type.

 A pair
  $\langle\r,\addrt\rangle$ with
  $\r < \addrt$ is called a \emph{characteristic ray pair
  for $\Ek$} if $\Ek$ has an orbit portrait whose characteristic rays
  are $\r$ and $\addrt$. More generally, $\langle\r,\addrt\rangle$
  is called a \emph{characteristic ray pair} if there exists some 
  $\kappa\in\C$ with such an orbit portrait. 
\index{characteristic ray pair}%
  If $\langle\r,\addrt\rangle$ is a characteristic ray pair of period $n$,
   then $\sigma^{n-1}(\addrt) < \sigma^{n-1}(\r)$.
    \qedd
\end{deflem}
\begin{remark}
 The final claim does not appear in the statement of 
  \cite[Lemma 5.2]{expattracting}, but is immediate from its proof.
\end{remark}

\begin{thmdef}[Characteristic Rays \protect{\cite[Lemma 5.2]{expattracting}}]
 \label{thmdef:characteristicrays}
 Let $\kappa$ be an attracting or parabolic parameter and let
  $n$ be the length of 
  $\s:=\extaddr(\kappa)$.
  Then there exists
  a unique characteristic ray pair $\langle\s^-,\s^+\rangle$
  for $\Ek$ such that the common landing point of $g_{\s^-}$ and
  $g_{\s^+}$ 
  lies on the boundary of the 
  characteristic Fatou component $U_1$; both rays have period $n$. 
  This ray pair separates $\kappa$ from all other periodic points with
  essential orbit portraits. 
\nomenclature[s-]{$\s^-$, $\s^+$}{(characteristic addresses of $\s$)}

 The addresses $\s^-$ and $\s^+$ depend only on $\s$, and are called the 
   \emph{characteristic addresses} of $\s$ (or $\kappa$). The common
   landing point is called the \emph{dynamic root} of $\Ek$. \qedd
 \index{dynamic root}%
  \index{characteristic addresses/rays!of an attracting or parabolic map}%
\end{thmdef}
\begin{remark}[Remark 1]
 In particular, 
  the dynamic root is the unique boundary point of the characteristic
  Fatou component $U_1$ which is fixed under $\Ek^n$ and which is the
  landing point of at least two dynamic rays.
\end{remark}
\begin{remark}[Remark 2]
 We will later
  give an algorithm (Algorithm \ref{alg:characteristicfromintermediate})
  for determining $\langle \s^-,\s^+\rangle$, given
  $\s$.
\end{remark}
\begin{remark}[Remark 3]
  The case of parabolic parameters was not formally treated in
   \cite{expattracting}. However, the proof is the same
   as that given there for attracting parameters.
   (Alternatively, the parabolic case
    follows from the attractive case
    by using Theorem
    \ref{thm:expper} below.) 
 \begin{suppress}
   If the parabolic orbit of $\Ek$ has an essential orbit portrait,
    then in particular the characteristic rays land at a point
    of the parabolic orbit. 
 \end{suppress}
\end{remark}

\subsection*{Itineraries and kneading sequences}
  Recall that, given an attracting or parabolic exponential map
   $\Ek$, one can connect the singular value to $\infty$ in 
   the characteristic Fatou component $U_1$ by a curve
   $\gamma$ at external address
   $\extaddr(\kappa)$. The preimage $\Ek^{-1}(\gamma)$ consists of countably
   many curves in $U_0$, and these produce a partition of the
   dynamical plane. 
   (The curve $\gamma$ is unique up to
   homotopy within $U_1$, so the partition is natural except for points
   within $U_0$.)
   To any point $z\in J(\Ek)$, one can now associate an
   \emph{itinerary}, which records through which strips of this
  \index{itinerary}%
   partition the orbit of $z$ passes. For more details, see
   \cite[Section 4]{expper}.

 In this article, we will
   use the following
   combinatorial analog of this notion. 
   If $\s\in\Sequ$, then $\sigma^{-1}(\s)$ produces a partition of
   $\Sequ$, and the itinerary of any $\t\in\Sequ$ will record
   where the orbit of $\t$ under $\sigma$ maps with respect to this
   partition. 

 \begin{defn}[Itineraries and Kneading Sequences] 
   \label{defn:combinatorialitinerary}
  Let $\adds\in \Sequ$ and $\r\in\Sequb$. 
  Then
   the \emph{itinerary of $\r$ with respect to $\adds$} is
   $\itin_{\adds}(\r) :=\u_1 \u_2 \dots$, where
\nomenclature[itin]{$\itin_{\adds}(\r)$}{(itinerary of $\r$)}
   \[ 
       \begin{cases}
        \u_k = \j      & \text{if \hspace{1mm}
                                 $\j\adds < \sigma^{k-1}(\r) < (\j+1)\adds$} \\
        \u_k = \itj &\text{if \hspace{1mm}
                                 $\sigma^{k-1}(\r)= \j\adds$}  \\
        \u_k = {\tt *}       & \text{if \hspace{1mm}
                                 $\sigma^{k-1}(\r)=\infty$}.
      \end{cases}\]
  (Note that $\itin_{\adds}(\r)$ is is a finite sequence if and only if
  $\r$ is an intermediate external address.)
  We also define 
   $\itin^+_{\adds}(\r)$ and
   $\itin^-_{\adds}(\r)$ to be the sequence obtained by replacing
  each boundary symbol $\itj$ by $\j$ or ${\tt j-1}$, respectively.
  When $\kappa$ is a fixed attracting or parabolic parameter, 
   we usually abbreviate $\itin(\r):=\itin_{\extaddr(\kappa)}(\r)$.

 \index{kneading sequence}%
  We also define the \emph{kneading sequence} of $\s$ to be 
   $\K(\s) := \itin_{\s}(\s)$. Similarly, the upper and lower kneading
   sequences of $\s$ are $\K^+(\s) := \itin^+_{\s}(\s)$ and
   $\K^-(\s) := \itin^-_{\s}(\s)$. 
\nomenclature[K]{$\K(\s)$}{(kneading sequence of $\s$)}%
\nomenclature[K+]{$\K^+(\s)$, $\K^-(\s)$}{(upper and lower 
   kneading sequences of $\s$)}%
 \end{defn}
\begin{remark}[Remark 1]
  One should think of $\s$ as lying in the
   ``combinatorial parameter plane'', whereas $\r$ lies in  the
   ``combinatorial dynamical plane'' associated with $\s$. 
\begin{suppress}
    Similarly, kneading sequences
   are parameter space objects, while itineraries are dynamical objects.
\end{suppress}%
\end{remark}
\begin{remark}[Remark 2]
  In the case $\adds=\infty$, we can define itineraries
    analogously, but the addresses $\j\adds$ and $(\j+1)\adds$ 
    in the definition will
    have to be replaced by $(\j-\frac{1}{2})\infty$ and
     $(\j+\frac{1}{2})\infty$. With this
    definition, $\itin_{\infty}(\r)=\r$ for all infinite 
    external addresses $\r$.
\end{remark}
\begin{remark}[Remark 3]
  The definition of itineraries involves a noncanonical choice
   of an offset for the labelling of the partition strips. 
   Our choice was made
   so that the external
   addresses in the interval $(\per{0},\per{1})$ 
   are exactly those whose kneading
   sequences start with ${\tt 0}$. 
\end{remark}

\medskip
The significance of itineraries lies in the fact that they can be used
 to determine which periodic rays land together, as shown in
 \cite[Theorems 3.2 and 5.4, Proposition 4.5]{expper}.

\begin{thm}[Dynamic Rays and Itineraries]
  \label{thm:expper}
 Let $\kappa$ be an attracting or parabolic parameter. Then every periodic
  dynamic ray of $\Ek$ lands at a periodic point, and conversely
  every repelling or parabolic periodic point is the landing point of
  such a ray.

 Two periodic rays $g_{\r}$ and $g_{\wt{\r}}$ land at the same point
  if and only if $\itin(\r)=\itin(\wt{\r})$. \qedd
\end{thm}

 The $n$-th itinerary entry $\u_n$ is locally constant (as a function of
  $\r$) wherever it is defined
  and an integer. If the $n$-th itinerary entry
  at $\r$ is a boundary symbol
  $\u_n=\itj$, then it is $\j-1$ slightly below $\r$ and
  $\j$ slightly above $\r$. In other words, 
   \begin{align*}
     \lim_{\t\nearrow \r} \itin_{\s}(\t) &= \itin_{\s}^-(\r) 
                              \hspace{0.5cm} \text{ and } \\
     \lim_{\t\searrow \r} \itin_{\s}(\t) &= \itin_{\s}^+(\r) 
   \end{align*}
  for all infinite external addresses $\r$. 
 
 If the $n$-th itinerary entry of $\r$ is $*$, then $\r$ is 
  an intermediate external address of length $n$. In this case,
  the $n$-th itinerary entries of addresses tend to 
  $+\infty$ when approaching $\r$ from below and to $-\infty$ when
  approaching $\r$ from above.

 We will frequently be in a situation where we compare the
  itineraries of an address $\t$ with respect to two different
  addresses $\s^1,\s^2\in\Sequ$. Therefore, let us state the
  following simple fact for further reference.
 \begin{observation}[Itineraries and Change of Partition] 
    \label{obs:itinchange}
  Let $\s^1,\s^2,\t\in\Sequ$ 
   with $\s^1<\s^2$
   and let $j\geq 1$ such that $\sigma^j(\t)$ is defined. 
   Then the
   $j$-th entries of the 
   itineraries of $\t$ with respect to $\s^1$ and $\s^2$
   coincide if and only if $\sigma^{j}(\t)\notin [\s^1,\s^2]$. 

  In particular, 
   $\itin_{\s^1}^-(\t)=
     \itin_{\s^2}^+(\t)$ 
   if and only if
   $\sigma^{k}(\t)\notin (\s^1,\s^2)$
   for all $k\geq 1$.
 \end{observation}
 \proof 
  By definition, the $j$-th itinerary entry of $\t$ as a function of
   $\s\in\Sequ$ is locally constant on $\Sequ\setminus\{\sigma^j(\t)\}$,
   which proves the ``if'' part. On the other hand, this function
   jumps by $1$ as $\s$ passes $\sigma^j(\t)$, proving the
   ``only if'' part. \qed

\begin{suppress}
Let $\addu:= \itin_{\s}^-(\t)$ 
   and
    $\addut := \itin_{\wt{\s}}^+(\t)$. 
    Clearly, if
    $\u_k=*$ for some $k$, 
    then $\t$ is an intermediate external
    address of length $k$, and 
    $\wt{\u}_k=*$. Let $k\geq 1$ with
    $\u_k\neq *$. Since $\sigma$
    preserves the circular order when
    restricted to the interval
    $(\u_k\s, (\u_k+1)\s]\ni
     \sigma^{k-1}(\t)$, 
    \[ \singlebox
       \wt{\u}_k = \u_k \quad 
     \Longleftrightarrow \quad
       \sigma^{k-1}(\t)\notin
       (\u_k\s, \u_k\wt{\s})
       \quad
   \Longleftrightarrow\quad
       \sigma^{k}(\t)\notin (\s,\wt{\s}).
        \esinglebox \]
\end{suppress}

 Finally, we shall require the following fact on the existence of
  addresses with prescribed itineraries.

 \begin{lem}[Existence of Itineraries] \label{lem:itineraries}
  Let $\s\in\Sequb$. Let $\addu$ be either an infinite sequence of
   integers or a finite sequence of integers followed by $*$, and suppose that
   $\sigma^k(\addu)\notin \{\K^+(\s),\K^-(\s)\}$ 
   for all $k\geq 1$. Then there exists
   an external address $\r\in\Sequb$ with $\itin_{\s}(\r)=\addu$; if 
   $\addu$ is
   periodic then every such $\r$ is also periodic.

  Furthermore, if
  $\t\in\Sequb$ and
  $k\geq 0$, then no two
  elements of
  $\sigma^{-k}(\t)$ have
  the same itinerary
  with respect to 
  $\s$.
  (In particular, no two  
  intermediate external 
  addresses have the
   same itinerary with respect
   to $\s$.)
 \end{lem}
 \begin{remark}
  The condition $\sigma^k(\addu)\notin \{\K^+(\s),\K^-(\s)\}$ is
  necessary. There exist periodic addresses $\s$ (for example
  $\s=\per{0}$) such that
  both $\K^+(\s)$ and $\K^-(\s)$ are not realized as the itinerary of
  any external address. Similarly,  as we will see in
  Lemma \ref{lem:kneadingandwakes} (\ref{item:kneadingbif}),
  there exist nonperiodic addresses
  $\s$ with periodic kneading 
  sequences. In this case,
  $\itin_{\s}(\r)\neq \K(\s)$ for all 
  periodic addresses $\r$. 
 \end{remark}

 \proof The set $R_k$
  of all external addresses $\r\in\Sequb$ 
 for which at least one of the itineraries
  $\itin_{\s}^+(\r)$ and
  $\itin_{\s}^-(\r)$ agrees with $\addu$ 
  in the first $k$ 
  entries is easily seen
  to be compact and nonempty for all 
  $k$. 
  Thus $R:= \bigcap_{k} R_k
           \neq\emptyset$. 
  Let $\r\in R$; then 
  $\addu\in\{\itin_{\s}^+(\r),
   \itin_{\s}^-(\r)\}$.

  We claim that $\itin_{\s}(\r)$ 
    contains no boundary symbols. 
   Indeed,
   otherwise $\sigma^k(\r)=\s$ for 
   some $k\geq 1$, and thus
   $\K(\s) = \sigma^k(\itin_{\s}(\r))$. 
  Thus
   $\KP(\s) = \sigma^k(\addu)$ or 
   $\KM(\s)=\sigma^k(\addu)$, which
   contradicts the assumption.
   Consequently,  
   $\itin_{\s}(\r)=\itin_{\s}^+(\r)=
     \itin_{\s}^-(\r)=\addu$ 
   as required. 

\smallskip

 The fact that $\r$ can be 
  chosen
  to be periodic if $\addu$ 
  is periodic is 
  \cite[Lemma 5.2]{expper}. 
  The proof that no 
  aperiodic address can have 
  the same itinerary as
   a periodic address is 
  analogous to 
  \cite[Lemma 18.12]{jackdynamics}.

\smallskip

 The final statement follows
  by induction from the
  trivial
  fact that changing 
  the first entry of
  an address $\t$ by some
  integer $m$ will also
  change the first entry of
  $\itin_{\s}(\t)$ by $m$. 
  \qed

\subsection*{Properties of characteristic ray pairs}

 As a first application of the concept of itineraries, let us deduce two
  basic properties of characteristic ray pairs.

 \begin{lem}[Characteristic Ray Pairs] \label{lem:charraypairs}
  Let $\langle\r^-,\r^+\rangle$ be a characteristic ray pair of period $n$.
  Then 
   there exist
   $\u_1,\dots \u_n\in\Z$ such that
   \begin{equation}
     \K^-(\r^-)=\K^+(\r^+)=
     \itin^-_{\r^-}(\r^+)=
     \itin^+_{\r^+}(\r^-)=
     \per{\u_1 \dots \u_n}. 
   \label{eqn:itineraries} 
  \end{equation}

 Furthermore, if $\s\in\Sequ$, then the following are equivalent.
  \begin{enumerate}
   \item $\s\in (\r^-,\r^+)$; \label{item:sinwake}
   \item $\displaystyle{\itin_{\s}(\r^-) = \itin_{\s}(\r^+)=
          \per{\u_1\dots \u_n}}$; \label{item:equalitinspecific}
   \item $\displaystyle{\itin_{\s}(\r^-) = \itin_{\s}(\r^+)}$.
      \label{item:equalitineraries}
  \end{enumerate}
 \end{lem}
 \proof 
  It follows from the definition of characteristic addresses that
   $\K^-(\r^-)=\itin_{\r^-}^-(\r^-)=
    \itin_{\r^-}^-(\r^+)=\per{\u_1\dots \u_n}$ for some
   $\u_1,\dots,\u_n\in\Z$. Furthermore, 
   no forward image of
   $\r^-$ or $\r^+$ belongs to $(\r^-,\r^+)$. 
   Thus 
   Observation \ref{obs:itinchange} implies that 
    \[ \itin_{\s}(\r^-) = \itin_{\s}(\r^+)=\per{\u_1\dots\u_n} \]
   for all $\s\in (\r^-,\r^+)$. This proves 
   (\ref{eqn:itineraries}) as well as 
    ``(\ref{item:sinwake})$\Rightarrow$(\ref{item:equalitinspecific})''.

   Clearly (\ref{item:equalitinspecific}) implies
    (\ref{item:equalitineraries}). To prove 
    that (\ref{item:sinwake}) follows from
    (\ref{item:equalitineraries}), let 
    $\s\in\Sequ\setminus (\r^-,\r^+)$. Since the
    interval $[\sigma^{n-1}(\r^+),\sigma^{n-1}(\r^-)]$ is mapped
    bijectively to $\Sequb\setminus (\r^-,\r^+)$ by the shift, 
    there is an element of $\sigma^{-1}(\s)$ between
     $\sigma^{n-1}(\r^+)$ and $\sigma^{n-1}(\r^-)$. So the $n$-th itinerary
     entries of $\r^-$ and $\r^+$
     with respect to $\s$ are different, as required. 
   \qed

\begin{lem}[Unique Intermediate Addresses]
  \label{lem:uniqueaddress}
 Let $\langle\r^-,\r^+\rangle$ 
  be a characteristic ray pair of period $n$, and let
  $\u_1,\dots \u_n\in\Z$ be as in the previous lemma.
  Moreover, let $\s$ be any intermediate external address of length $n$. 
  Then the following
  are equivalent:
 \begin{enumerate}
  \item \label{item:charaddress}
    $\r^-$ and $\r^+$ are the characteristic addresses of $\s$;
  \item \label{item:itinofs}
   $\displaystyle{\itin_{\r^-}(\s)= \itin_{\r^+}(\s)=\K(\s) =
           \u_1 \dots \u_{n-1} *}$;
  \item $\s\in (\r^-,\r^+)$ and $\K(\s) = \u_1 \dots \u_{n-1} *$.
    \label{item:equalitin}
 \end{enumerate}
  Furthermore, there is
  at most one $\s$ with one (and thus all) of these properties.
\end{lem}
 \begin{remark}
  Note that we do not claim here that such an address $\s$ always
  exists. While it is not very difficult to show that this is indeed the case,
  we will not require this fact until Section \ref{sec:internaladdresses},
  and it will be proved there (Proposition \ref{prop:periodicintermediate}).
 \end{remark}
\proof 
 If $\r$ and $\r^+$ are the 
  characteristic adresses of $\s$, then 
  no
  forward image of $\s$ lies in 
  $(\r^-,\r^+)$. Indeed, let $\kappa$
  be a parameter with 
  $\extaddr(\kappa)=\s$. Then the
  characteristic rays $g_{\r^-}$ and
  $g_{\r^+}$ separate the
  characteristic Fatou component $U_1$
  from the other periodic Fatou 
  components, and thus from the remainder of the singular orbit.
  So (\ref{item:charaddress}) implies (\ref{item:itinofs}) by
  Observation \ref{obs:itinchange}. 

\smallskip

 Suppose that (\ref{item:itinofs}) holds. Then, for 
   each $k\in\{1,\dots,n-1\}$, the three addresses
    $\sigma^{k-1}(\r^-)$, 
      $\sigma^{k-1}(\r^+)$ and $\sigma^{k-1}(\s)$
    belong to the interval $(\u_k \r^- , (\u_k+1)\r^-]$.
    By
    Observation \ref{obs:orderpreserving}, the map 
    $\sigma^{n-1}$ thus
    preserves the circular order of $\r^-$, $\r^+$ and
    $\s$. Recall that 
    $\sigma^{n-1}(\r^+)<\sigma^{n-1}(\r^-)$ by Lemma
    \ref{lem:permutation}. Since
    $\sigma^{n-1}(\s)=\infty$, it follows that
    $\r^-<\s<\r^+$. By Lemma \ref{lem:charraypairs}, it follows that
    (\ref{item:equalitin}) holds.

\smallskip

 Now let us assume that $\s$ satisfies (\ref{item:equalitin}), and
  let
  $\kappa\in\C$ with 
  $\extaddr(\kappa)=\s$. Then the dynamic rays
  $g_{\r^-}$ and $g_{\r^+}$ have a common landing point $w$ by
  Lemma \ref{lem:charraypairs} and
  Theorem \ref{thm:expper}. It easily
  follows that $\r^-$ and $\r^+$ are the characteristic addresses of
  the orbit portrait of $w$.

  Now let $\s^-$ and $\s^+$ be the characteristic addresses
  of $\s$. Then
  we have
   \[ \r^- \leq \s^- < \s < \s^+ \leq \r^+, \]
  and the first $n-1$ itinerary entries of all
  these addresses with respect to
  $\s$ 
  coincide. As above, the cyclic order of this configuration is preserved
  under $\sigma^{n-1}$. Since $\sigma^{n-1}(\s)=\infty$, this means that
   \begin{equation} \label{eqn:chainofinequalities}
     \sigma^{n-1}(\s^+) \leq
     \sigma^{n-1}(\r^+) <
     \sigma^{n-1}(\r^-) \leq \sigma^{n-1}(\s^-). 
   \end{equation}
  Since the outer two addresses in (\ref{eqn:chainofinequalities}) have
   the same itinerary, this means that the $n$-th itinerary entries of
   $\s^-,\s^+,\r^-$ and $\r^+$ also agree. Thus all these addresses belong
   to the orbit portrait of $w$.
   Since both pairs $\langle\s^-,\s^+\rangle$ and 
   $\langle\r^-,\r^+\rangle$
   are characteristic ray pairs of this portrait, they must be equal.
   This proves (\ref{item:charaddress}).

\smallskip

 Finally, by Lemma \ref{lem:itineraries}, 
  there is at most one $\s\in\Sequ$ with
  $\itin_{\r^-}(\s)=\u_1 \dots \u_{n-1} *$, which completes the proof.
  \qed

\section{Hyperbolic components}
 \label{sec:hyperboliccomponents}

\begin{figure}
 \center
 \includegraphics[height=0.9\textheight]{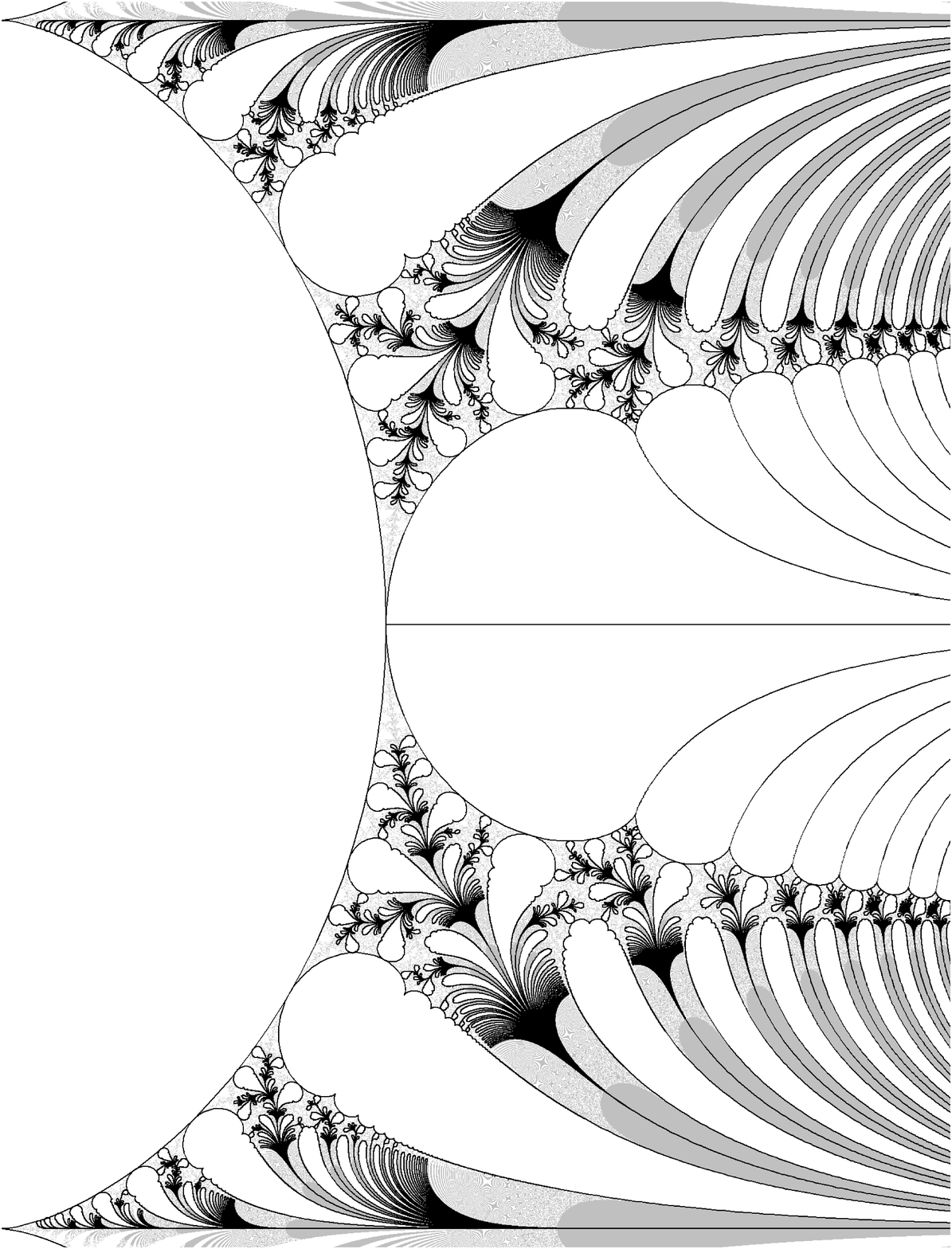}
 \caption[Hyperbolic Components]{Several hyperbolic components in the
 strip $\{\im \kappa \in [0,2\pi]\}$. Within the period two component
 in the center right of the picture, parameter rays at integer heights
 are drawn in.\label{fig:hyperboliccomponents}}
\end{figure}

A \emph{hyperbolic component} $W$ is a maximal connected
  \index{hyperbolic component}%
 open subset of parameter space
 in which each parameter has an attracting periodic orbit. It is easy
 to see that this (unique) orbit depends holomorphically on $\kappa$, and in 
 particular its period is constant throughout $W$. 

It was shown in \cite{bakerexp} that every hyperbolic
 component $W$ is simply connected and unbounded.
 Furthermore,
 the multiplier map
 $\mu:W\to\Ds$, which maps every parameter to the multiplier of its
\nomenclature[mu]{$\mu$}{(multiplier map)}
 attracting cycle, is a universal covering.
  Since $\exp:\H\to\Ds$ is also a universal covering,
  there exists a conformal isomorphism 
  $\Psi_W:\H\to W$ with $\mu\circ\Psi_W=\exp$. Note that this defines
  $\Psi_W$ uniquely only up to precomposition by
  a deck transformation of $\exp$ (i.e., an additive translation by
  $2\pi ik$, $k\in\Z$). 
  However, in \cite[Theorem 7.1]{expattracting}, it was shown
  that, when $W$ has period at least two,
  there is a unique choice of $\Psi_W$ such that, for any 
  $\kappa\in \Psi_W\bigl((-\infty,0)\bigr)$, the dynamic root lies on the
  distinguished boundary orbit. (This also follows from Theorem
  \ref{thm:labellingmap} below.) In the following, we will always fix
  $\Psi_W$ to be this \emph{preferred parametrization}. 
\index{preferred parametrization}%
\nomenclature[Psi]{$\Psi_W$}{(preferred parametrization of $W$)}

 It is well-known \cite{bakerexp}
  that there exists a unique hyperbolic component
  $W_0$ of period $1$. This hyperbolic component contains a left half plane
  and is invariant under $z\mapsto z+ 2\pi i $. 
  Since different
  choices for the parametrization of this component only correspond to
  such a translation, and thus to a relabeling of the map, there is
  no canonical choice for a preferred parametrization. By definition,
  we choose the \emph{preferred parametrization} of the period $1$
  component to be the unique map respecting the real symmetry, i.e.\ 
  \[ \Psi_W : \H \to W_0;
      z\mapsto z - \exp(z). \]

 \subsection*{Internal rays}
 We can now foliate $W$ by curves, called
  \emph{internal rays}, along which the
  \index{internal ray}%
  argument of $\mu$ is constant. More precisely, 
  the \emph{internal ray at height $h\in\R$} is the curve
  \[ \IR:(-\infty,0)\to\C, t\mapsto \Psi_W( t + 2\pi i h).
  \]
\nomenclature[Gamma]{$\IR$}{(internal ray at height $h$)}

 It is straightforward to see that, if $W$ is a hyperbolic component of
  period $\geq 2$, then 
  $\IR(t)\to+\infty$ as $t\to\infty$, and the homotopy class
  of this curve as $t\to\infty$ is independent of $h$ (see e.g.\ 
  \cite[Lemma 2.1]{expattracting}).
  We say that an internal ray $\IR$
  \emph{lands} at a point $\kappa\in\Ch$ if $\kappa=\lim_{t\nearrow 0} \IR(t)$.

 \begin{lem}[Landing of Internal Rays] \label{lem:internallanding}
  Let $W$ be a hyperbolic component of period $n$. Then 
  every internal ray $\IR$ lands at some
  point in $\Ch$, which we denote by $\Psi_W(2\pi i h)$. The set of $h$
  for which this landing point belongs to $\C$ is open and dense.
  Conversely, suppose that $\kappa_0\in\partial W\cap \C$. Then
   $\kappa_0$
   is the landing point of a unique internal ray $\IR$. Furthermore,
   $\kappa_0$ is an indifferent parameter of period 
   dividing $n$. If $a$ is a point on this indifferent periodic orbit,
   then $(E_{\kappa_0}^n)'(a)=\exp(2\pi ih)$. 
 \end{lem}
 \begin{remark}[Remark 1]
  It is not difficult to see that the extended map
   $\Psi_W$ is continuous on
   $\cl{\H}$. However, we shall not require this result here.
 \end{remark}
 \begin{remark}[Remark 2]
  This lemma leaves open the possibility that some internal rays
  land at $\infty$, disconnecting the boundary of $\C$. Proving that this
  does not happen is much more difficult. The proof of this fact,
  in 
  \cite{boundary}, uses the results of the present
  article. (An outline of the argument was given in
  \cite{cras}.) 
  This result can therefore not be used in the following sections.
 \end{remark}
 \sketch
  It is straightforward to see that every point of $\partial W\cap\C$
  has an indifferent periodic orbit of period dividing $n$. The multiplier
  map $\mu$ 
  extends to a holomorphic function on a neighborhood of $\kappa_0$ 
  (or on a finite-sheeted covering of $\kappa_0$ when $\mu(\kappa_0)=1$,
  compare \cite[Proof of Lemma 4.2]{jackrays}).
  Hence there exists some internal ray
  which lands at $\kappa_0$. If there were
  two internal rays of $W$ landing at $\kappa_0$, these could be
  connected to form a simple
  closed curve $\gamma\subset \cl{W}$ 
  which separates some part of $\partial W$ 
  from infinity. It is a standard fact that this is not possible:
  for example, it is well-known that every point of $\partial W\cap\C$ is
  structurally unstable, and thus can be approximated by 
  attracting parameters with arbitrarily high periods
  (compare \cite[Lemma 5.1.6]{thesis}). 
  The hyperbolic components containing 
  these parameters thus will be separated from $\infty$ by $\gamma$,
  which is impossible since every hyperbolic component is unbounded.
   
 The fact that the set of $h$ for which $\IR(0)$ lands in $\C$ is open
  follows easily from the above statement about the multiplier.
  That this set is also dense (and in fact has full measure)
  follows from the
  F.\ and M.\ Riesz theorem \cite[Theorem A.3]{jackdynamics}.

 Finally, it is straightforward to see that every
 finite limit point $\kappa_0$
   of an
   internal ray $\IR$ has a periodic point $a$ with
   $(E_{\kappa_0}^n)'(a)=\exp(2\pi ih)$. The set of such parameters is
   easily seen to be discrete in $\C$, proving that 
   $\IR$ lands at a point of $\Ch$. 

  For a more detailed (self-contained!) 
   proof of this lemma, compare
   \cite[Section 2]{boundary}. \qed

 When $W$ is of period at least $2$, 
 the internal ray at height $0$ is called the \emph{central
 internal ray} of $W$. If this ray lands at a point in $\C$, then its
  \index{internal ray!central}%
 landing point $\Psi_W(0)$ is called the \emph{root} of $W$.
  \index{root of a hyperbolic component}%
 The points of $\Psi_W(2\pi i \Z\setminus\{0\})\cap\C$ 
  are called \emph{co-roots} of
  \index{co-roots of a hyperbolic component}%
 $W$. For the period $1$ component, all points 
 of $\Psi_W(2\pi i \Z)$ are called co-roots.
 Recall that we will prove in \cite{boundary} that every component
 has a root (and infinitely many co-roots).
 Without this knowledge, which we may not use at this point, we cannot be
 sure that bifurcations actually exist.

\subsection*{Classification of hyperbolic components}
It is easy to see that $\extaddr(\kappa)$ depends only on the
 hyperbolic component $W$ which contains $\kappa$; this address will
 therefore also be denoted by $\extaddr(W)$. 
 \index{external address!of a hyperbolic component}%
 \index{intermediate external address!of a hyperbolic component}%
\nomenclature[addrW]{$\extaddr(W)$}{(external address of $W$)}%
 Similarly, we will talk
 about the kneading sequence, characteristic rays etc.\ of
 $W$. The following theorem, which is the main result of
 \cite{expattracting}, states that
 hyperbolic components can be completely classified terms of their
 combinatorics. (Note that the existence part of this result was
 already cited as
 Theorem \ref{prop:hypexistence}.)
\begin{prop}[Classification of Hyperbolic Components \cite{expattracting}]
 \label{prop:hyperboliccomponents}
 For every intermediate external address $\s$, there exists exactly one
 hyperbolic component $W$ with $\extaddr(W)=\s$. We denote this component
 by $\Hyp{\s}$. The vertical order of hyperbolic components coincides
 with the lexicographic order of their external addresses. \qedd
\nomenclature[Hyp]{$\Hyp{\s}$}{(hyperbolic component at address $\s$)}
\end{prop}
	
To explain the last statement, recall that, when
  $W$ is a hyperbolic component of
  period $\geq 2$, any internal ray $\IR$ satisfies
 $\IR(t)\to+\infty$ as $t\to -\infty$.
  Thus the family of central
  internal rays of hyperbolic components has a natural vertical order%
 \footnote{We should stress that we only use the ``negative'' ends of
  internal rays to define this order. A priori some internal rays might also
  tend to infinity as $t\to +\infty$. Taking these
  directions would result
  in a different order which we do not refer to.}\ 
  as described in Section \ref{sec:combinatorics}, and this is the
  order referred to in the Proposition. (Note that taking the
  central rays is not essential, as there is only one homotopy class of
  curves in $W$ along which the multiplier tends to $0$.)

\subsection*{Sectors}
 Since $\mu:W\to\Ds$ is a universal covering, parameters in $W$ are not 
  (as in the quadratic family) uniquely determined by their multiplier.
  Rather, the set $\mu^{-1}(\Ds\setminus [0,1))$ consists of countably
  many components, called \emph{sectors} of $W$. If 
   \index{sector}%
  $\kappa=\Psi_W(t+2\pi i h)$ with $h\notin\Z$, we denote the
  sector containing $\kappa$ by
  \[ \Sec(\kappa) := \Sec(W,h) := 
     \Psi_W\Bigl(\bigl\{a+2\pi i b: a < 0 \text{ and }
                   b\in ( \floor{h} , \ceil{h} ) \bigr\} \Bigr). \]  
\nomenclature[Seckappa]{$\Sec(\kappa)$}{(sector of $\kappa$)}%
\nomenclature[SecWh]{$\Sec(W,h)$}{(sector of $W$)}%
 
\begin{defn}[Sector Labels] \label{defn:sectorlabels}
  Let $W=\Hyp{\adds}$ be a hyperbolic component and
  let $\kappa\in W$ be a parameter with $\mu:=\mu(\kappa)\notin (0,1)$.
  Let $\gamma$ be the principal attracting ray of $\Ek$, and let
  $\gamma'$ be the component of $\Ek^{-1}(\gamma)$
  which starts at $a_0$. Then $\extaddr(\gamma')$ is of the form
  $s_* \adds$ with $s_*\in\Z$ (resp.\ $s_*\in\Z+\frac{1}{2}$ if
  $\adds=\infty$).
  The entry $s_* = s_*(\kappa)$ is called
  the \emph{sector label} of $\kappa$.
   \index{sector label}%
   \nomenclature[s*]{$s_*(\kappa)$}{(sector label)}%
\end{defn}
\begin{remark}
 There are two more ways to label sectors which will appear in this
  article: \emph{sector numbers} and
  \emph{kneading entries}; both will be introduced in Section
  \ref{sec:internaladdresses}. 
  We should warn the reader that our terminology is somewhat different
  from that of \cite{intaddr}, where the term \emph{sector label} is
  used to refer to what we call \emph{kneading entries}.
 \end{remark}

The following results justify the term ``sector label''; compare
 Figure \ref{fig:seclabels}.

\begin{thm}[Behavior of Sector Labels]
 \label{thm:labellingmap}
  The map $\kappa\mapsto s_*(\kappa)$ is constant on sectors of $W$.
   When $\kappa$ crosses a sector boundary so that $\mu$ passes 
   through $(0,1)$ in
   positive orientation, then $s_*(\kappa)$ increases exactly by $1$.
   In particular the induced map from sectors to indices is
   bijective. The unique sector with a given sector label $s_*$
   will be denoted by $\Sec(W,s_*)$. 
   \nomenclature[SecWs]{$\Sec(W,s_*)$}{(sector with sector label $s_*$)}%

 If the period of $W$ is at least $2$ and
  $\kappa$ is a parameter on the internal ray between
  $\Sec(W,j)$ and $\Sec(W,j+1)$, then the distinguished
  boundary fixed point of $\Ek^n$ on $\partial U_1$ has itinerary
  $\per{\u_1 \dots \u_{n-1} j}$ (where $\K(W) = \u_1 \dots \u_{n-1} *$). 
  In particular, the central internal ray of $W$ is 
  the boundary between $\Sec(W,j)$ and $\Sec(W,j+1)$, where $j$ and
  $j+1$ are
  the $n$-th entries of the characteristic addresses
  $\s^+$ and $\s^-$, respectively. 

 If $W$ is the unique period one component, a similar statement holds:
  the boundary between $\Sec(W,j-\frac{1}{2})$ and
  $\Sec(W,j+\frac{1}{2})$ is given by the internal ray
  $\{t+2\pi i j: t\leq -1\}$. For parameters on this ray, the
  distinguished boundary fixed point of $\Ek$ has itinerary
  $\per{j}$. 
\end{thm}
\sketch
 The linearizing coordinate used to define attracting
 dynamic rays depends holomorphically on $\kappa$. It easily follows
 that, as long as the principal attracting
 ray $\gamma$ does not pass through $\kappa$, its preimage 
 $\gamma'=\gamma'(\kappa)$
 from Definition \ref{defn:sectorlabels} varies continuously,
 which shows that $s_*$ is constant on sectors.

 In the following, let us restrict to the case where the period of
  $W$ is at least two; the case of the period one component is
  handled analogously. 
  Let $\kappa_0$ be a parameter with positive
  real multiplier, and let us set
  $\s := \extaddr(W)$. Then the principal attracting ray $\gamma$
  contains the singular value. Denote
  the piece of $\gamma$ which connects $a_1$ to $\kappa$ by 
  $\gamma_0$ and the piece which connects $\kappa$ to $+\infty$ by
  $\gamma_1$. We can define a branch $\phi$ of $\Ek^{-1}$ on 
  $U_1\setminus \gamma_1$ which takes $a_1$ to $a_0$. The range of $\phi$
  is then a strip $S$ of $U_0$ bounded by two consecutive preimages
  ${\gamma'}_1^-$ and ${\gamma'}_1^+$, at external addresses
  $j\s$ and $(j+1)\s$ for some $j\in\Z$.

 Let $\alpha\subset U_1\setminus \gamma_1$ be an unbroken attracting dynamic 
  ray, connecting $a_1$ to the distinguished boundary fixed point
  $w\in\partial U_1$. Then the image of $\alpha$ under
  $\phi$ is a curve connecting $a_0$ to $\Ek^{n-1}(w)$. Thus 
  $\Ek^{n-1}$ lies between ${\gamma'}_1^-$ and ${\gamma'}_1^+$, so 
  $j$ is the $n$-th itinerary entry of $w$ as claimed.

 Denote the preimage of $\gamma_0$ in $S$ by $\gamma'_0$. We then
  define two curves (in $\Ch$),
  \[ {\gamma'}^{\pm} := \gamma_0'\cup \{-\infty\} \cup  
                      {\gamma'}_1^{\pm}. \]
  By continuity of the linearizing coordinate, it then follows that,
  as $\kappa\to\kappa_0$ through parameters at positive
  (resp.\ negative) multiplier angles, the curves 
  $\gamma'(\kappa)$ converge uniformly 
  to ${\gamma'}^+$ (resp.\ ${\gamma'}^-$), which completes the proof.
  (Compare \cite[Theorem 5.5.3]{thesis} for more details.) \qedd

\section{Local Bifurcation Results}
 \label{sec:localbifurcation}

 Throughout this section, let $\kappa_0$ be a parabolic parameter of period 
  $n$ and intermediate external address $\s$. If the
  parabolic orbit portrait of $E_{\kappa_0}$ is essential, then we call 
  $\kappa_0$ a \emph{satellite} or a \emph{primitive} parameter, depending
  on the type of this orbit portrait. Similarly, we will refer to the ray
  period of this orbit portrait as the ray period of $\kappa_0$. Note that this
  ray period is also the period of the repelling (or attracting) petals
  of the parabolic orbit. 
    \index{ray period!of a parabolic parameter}%
    \index{satellite parameter}%
    \index{primitive parameter}%

  In this section, we will study what happens when $\kappa_0$ is perturbed
   into an adjacent hyperbolic component. 
   For this purpose, we will use 
   the following well-known
   statement about the analytic structure near $\kappa_0$.
   This result is beautifully exposed, and proved using elementary
   complex analysis, in \cite[Section 4]{jackrays}. All that this
   local analysis requires is that there is only one singular value,
   and thus only one single cycle of
   petals at a parabolic periodic point. 

 \begin{prop}[Perturbation of Parabolic Orbits]
   \label{prop:parabolicperturbation}
  Let $\kappa_0$ be a parabolic parameter of period $n$, with ray period
   $qn$.
  \begin{itemize}
   \item \textbf{(Primitive and Co-root case)}
     If $q=1$ (so the multiplier of the parabolic orbit is $1$), then, under
     perturbation, the parabolic orbit splits up into two orbits of period
     $n$ that can be defined as holomorphic functions of a two-sheeted cover
     around $\kappa_0$. 

     Any hyperbolic component whose boundary contains $\kappa_0$
     corresponds to one of these orbits becoming attracting (and therefore
     has period $n$).
   \item \textbf{(Satellite Case)}
     If $q\geq 2$, then, under perturbation, the parabolic orbit splits
      into one orbit of period $n$ and one of period $qn$. The period
      $n$ orbit can be defined as a holomorphic function in a neighborhood
      of $\kappa_0$, as can the multiplier of the period $qn$-orbit. The
      $qn$-orbit itself can be defined on a $q$-sheeted covering around
      $\kappa_0$. 

      Any hyperbolic component whose boundary contains
       $\kappa_0$ corresponds to
       one of these orbits becoming attracting (and therefore has
       period $n$ or $qn$). \qedd
  \end{itemize}
 \end{prop}
  Any hyperbolic component of period $qn$ that touches $\kappa_0$ 
   is called a \emph{child
   component}; note that at least one such component always exists.
   In the satellite case, any period $n$ component touching $\kappa_0$
   is called a \emph{parent component}. 
   \index{child component!of a parabolic parameter}%
   \index{parent component}%
  Note that any satellite parameter has at least one child and
   at least one parent component; a primitive parameter has at least
   one child component but no parent components. 
   (We will show in Theorems 
   \ref{thm:bifurcationchild} and
   \ref{thm:bifurcationparent}
   below that ``at least one'' can be 
   replaced by ``exactly one''.)

 We will also require the following statement on the landing behavior
  of periodic rays as $\kappa_0$ is perturbed. The proof is analogous
  to that in the case of quadratic polynomials, which can also be found
  in \cite[Section 4]{jackrays}.
  (Recall that, by the previous proposition, the parabolic orbit of
   $\kappa_0$ breaks up into two orbits under perturbation. If we
   perturb $\kappa_0$ into an adjacent hyperbolic component, one of these
   orbits becomes attracting, so
   there is a unique repelling orbit
   created in the bifurcation.)

 \begin{prop}[Orbit Stability under Perturbation]
   \label{prop:orbitstability}
  Under perturbation of a parabolic parameter $\kappa_0$
   into a child or parent
   component,
   all repelling periodic points retain the same orbit portraits. 

  Furthermore, under perturbation into a child component, the 
   repelling periodic orbit created in the bifurcation has the same
   orbit portrait as the parabolic orbit of $E_{\kappa_0}$. Under
   perturbation into a parent component, the rays landing at the
   parabolic orbit are split up, landing at distinct points of the
   newly created repelling
   orbit. \qedd
 \end{prop}

 We are now ready to describe the combinatorics of child and parent 
   components of $\kappa_0$ (and, in particular, show that 
   there is at most one of each, see Corollary \ref{cor:parabolicbifurcation}.

\begin{figure}
 \subfigure[parent component]{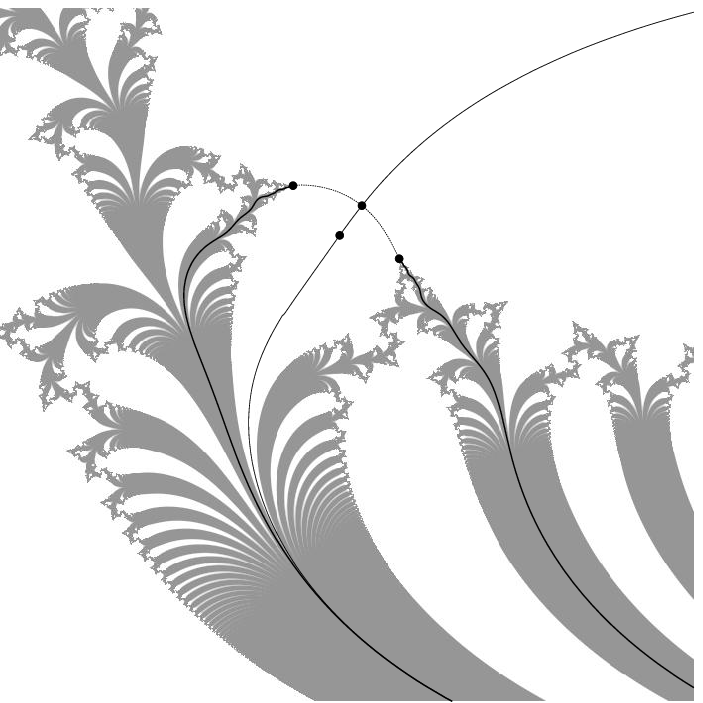%
     \label{fig:parentcomponent}}\hfill
 \subfigure[child component]{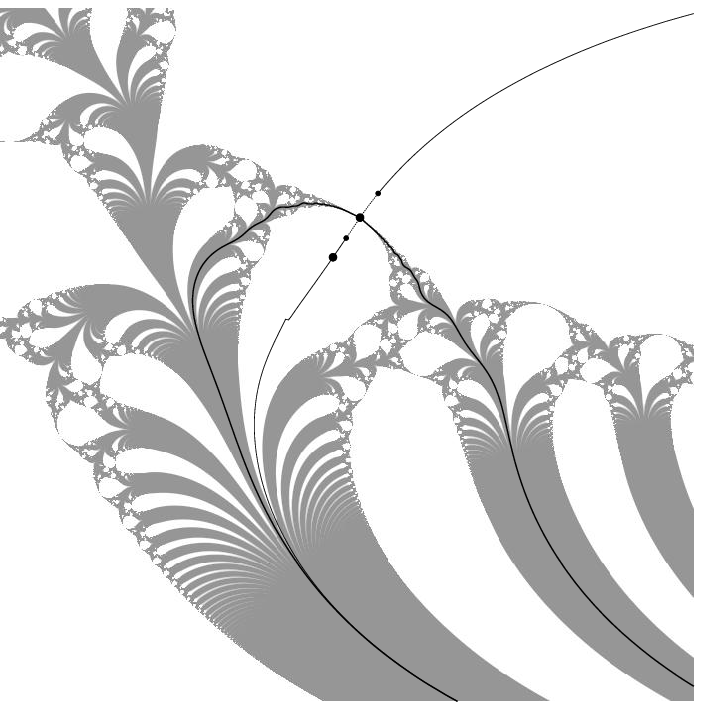%
     \label{fig:childcomponent}}
 \caption{The dynamical plane just before and after a 
  bifurcation, illustrating Theorems
  \ref{thm:bifurcationchild} and \ref{thm:bifurcationparent}. 
  In both pictures, unbroken attracting rays are dotted lines,
  broken attracting rays are solid lines, and
  dynamic rays are strong solid lines.}
\end{figure}

 \begin{thm}[Combinatorics in a Child Component]
   \label{thm:bifurcationchild}
  Let $\kappa_0$ be a parabolic parameter of period $n$ and ray period $qn$,
   and let $W$ be a child component of $\kappa_0$.  Then 
   $\extaddr(W)=\extaddr(\kappa_0)$; i.e., $W$ is the unique
   component at address $\extaddr(\kappa_0)$.

   Furthermore, for
   points on the internal ray of $W$ landing at $\kappa_0$, the
   repelling point created in the bifurcation is the distinguished boundary
   fixed point. Therefore
   $\kappa_0$ is the root point of $W$ if and only if its parabolic
   orbit portrait is essential; otherwise, $\kappa_0$ is a co-root of
   $W$.
 \end{thm}
 \proof
   If $qn=1$, then $W=\Hyp{\infty}$ is
   the unique component of period $1$. Now suppose that $qn>1$. Let
   $\IR$ be the unique internal ray
   landing at $\kappa_0$ (compare Lemma
   \ref{lem:internallanding}); this ray
   has integer height $h\in\Z$. Let
   $\kappa:=\IR(t)$ be a parameter on this ray.
   By 
   Proposition
   \ref{prop:orbitstability}, $\kappa_0$ and $\kappa$ have the same
   orbit portraits and thus they have the same characteristic addresses.
   Lemma \ref{lem:uniqueaddress} then yields
   $\extaddr(\kappa_0)=\extaddr(\kappa)=\extaddr(W)$.

  Now let $w$ be the newly created repelling point and let
   $w'$ be the distinguished boundary fixed point of $\kappa$.
   Let $\alpha$ be the 
   piece of the principal attracting ray of $\Ek$ which connects
   $\kappa$ to $\infty$. Recall from Lemma 
   \ref{deflem:distinguishedboundaryorbit} that $a_0$ and
   $\Ek^{n-1}(w')$ can be connected by an unbroken dynamic ray of $\Ek$ and
   thus belong to the same component of $\C\setminus \Ek^{-1}(\alpha)$. 
   We will show that $a_0$ and $\Ek^{n-1}(w)$ are also not separated by
   $\Ek^{-1}(\alpha)$. This implies that $w$ and $w'$ have the same
   itinerary and are therefore equal by Theorem
   \ref{thm:expper}, as required.

 Let $\Phi:U_1\to \C$ be the linearizing coordinate for $\Ek$,
  normalized so that $\Phi(\kappa)=1$, and let
  $V\subset U_1$ be the component of the 
  preimage of $\Phi^{-1}(\D(0,\frac{1}{\mu}))$
  which contains $a_1$. (Note that 
  $\D(0,\frac{1}{\mu})$ is the largest disk on which $\Phi^{-1}$ exists, and
  that $\mu=\exp(t)$.)
  Then $V$ contains the
  curve $\gamma:=\Ek^n(\alpha)$, 
  and  by definition $\Phi(\gamma) = (\mu,1)$. Since
  $\Phi|_V:V\to \D(0,\frac{1}{\mu})$ is a conformal isomorphism, the
  hyperbolic length of $\gamma$ in $V$ is
   \[ \ell_V(\gamma)=
      \int_{\mu}^1 \frac{\mu|dz|}{1-(\mu|z|)^2} \leq
      \frac{\mu(1-\mu)}{1-\mu^2}=\frac{\mu}{1+\mu}\leq 2. \]
  In other words, the hyperbolic length of $\gamma$ within $U_1$
   stays bounded
   as $t\to 0$. Since the euclidean length of
   $\gamma$ (which connects $\kappa$ and
   $\Ek^n(\kappa)$) is bounded below as $t\to 0$, 
   it follows from standard estimates on the hyperbolic metric,
   using the fact that $w\in \partial U_1$,
   that the euclidean distance
   $\dist(w,\gamma)$ is also bounded below as $t\to 0$. 

  Since $\Ek^{n+1}$ is continuous in $z$ and $\kappa$,
    \[ \liminf_{t\to 0}
        \dist(\Ek^{n-1}(w),\Ek^{-1}(\alpha)) > 0 \]
   (recall that $\Ek^n(w)=w$). 
  On the other hand, the attracting point $a_0$ and 
   the repelling point $\Ek^{n-1}(w)$ are created from the same parabolic
   point, so the distance between them tends
   to $0$ as $t\to 0$. Thus, for small enough $t$, these two points
   are not separated by $\Ek^{-1}(\alpha)$. Therefore we have shown
   $w=w'$. 

  To prove the final statement, observe that by Proposition
   \ref{prop:orbitstability},
   the parabolic orbit portrait of
   $E_{\kappa_0}$ is essential if and only if the orbit portrait of
   $w=w'$ is essential for $E_{\kappa}$. By definition, the latter is the
   case if and only if $\kappa$ belongs to the central internal ray;
   i.e., if $\kappa_0$ is the root point of $W$. \qed

Let us now turn to determining the parent component of a satellite
 parameter.
 \begin{thm}[Combinatorics in a Parent Component]
   \label{thm:bifurcationparent}
  Let $\kappa_0$ be a satellite parabolic parameter of period $n$ and
   multiplier $e^{2\pi i \frac{p}{q}}$; let $W$ be a parent component of
   $\kappa_0$. 
   If $\kappa$ is a parameter
   on the internal ray in $W$ which lands at
   $\kappa_0$, then the distinguished
   boundary orbit for $\kappa$ is the repelling
   period $qn$ orbit created from the parabolic point. 

  Furthermore, for such a parameter $\kappa$,
   let $\adds'$ be the address of the attracting dynamic ray
   which contains the singular value.
   Then $\extaddr(\kappa_0)=\adds'$. 
 \end{thm}
 \proof
  To prove the first statement, 
   choose some small wedge-shaped repelling petals for the
   parabolic orbit of $E_{\kappa_0}$. If
   $\kappa$ is close enough to $\kappa_0$ on the internal ray
   landing at $\kappa_0$, then, for $\kappa$, these petals
   (after a translation moving them to the attracting period $n$ orbit
    $(a_k)$
    created in the bifurcation)
   are still
   backward invariant and contain the newly created repelling
   period $qn$ orbit $(w_k)$
   (see \cite[Section 4]{jackrays}).
   Choose any attracting dynamic ray
   which approaches $a_1$ through one of these petals. Recall that this ray is
   periodic under $\Ek^n$, of period $q$. Pulling back, it
   must land at the unique fixed point of the return map of this petal.
   That is, the points of the repelling orbit $(w_k)$ are 
   landing points of a cycle of attracting dynamic rays; thus
   $(w_k)$ is the distinguished boundary orbit. 

  To prove the second part of the theorem, 
   let $\r,\rt$ be the characteristic
   addresses of $\kappa_0$.
   We will now find the combinatorial features of
   $\kappa_0$ within the attracting dynamics of $\kappa$, using 
   attracting dynamic rays 
   (compare Figure 
   \ref{fig:childcomponent}.)   
   Let $\gamma_0$ be the attracting dynamic ray containing
   the singular value (recall that
   $\adds'=\extaddr(\gamma_0)$). Then the attracting rays
   $\gamma_j:=\Ek^{nj}(\gamma_0)$, $j=0,\dots,q-1$ completely contain the
   singular orbit in $U_1$.
   By Proposition \ref{prop:orbitstability} and the first part
   of the theorem, 
   the rays $g_{\r}$ and $g_{\addrt}$ 
   do not land together for $\Ek$, but rather
   land separately on two points of the
   distinguished boundary cycle on $\partial U_1$. These landing points can
   be
   connected to $a_1$ by
   two
   attracting dynamic rays.
   Let $h_{\r}$ and $h_{\addrt}$ denote the curves
   obtained by extending $g_{\r}$ and $g_{\addrt}$ to $a_1$ by these
   attracting rays.

   Now consider the preimages of that part of $\gamma_0$ 
   which connects $\kappa$ to
   $\infty$; these curves are straight lines in the linearizing coordinate of 
   $U_0$, and connect
   $-\infty$ to $+\infty$ with external addresses of the form
   $m\s'$. The
   images of $h_{\r}\cup\gamma_0\cup h_{\addrt}$ 
   do not intersect these curves. Consequently for every $j\geq 0$,
   the $j$-th iterated images of $g_{\r}$,
   $g_{\addrt}$ and $\gamma_0$ belong to a common strip of this
   partition.

  Thus, $\itin_{\adds'}(\r)=\itin_{\adds'}(\addrt)$; so in particular
   $\s\in (\r,\addrt)$ by Lemma \ref{lem:charraypairs}. Furthermore,
   this itinerary agrees
   with $\K(\adds')$ on its first
   $qn-1$ entries, which implies
   by Lemma \ref{lem:uniqueaddress} that 
   $\r$ and $\addrt$ are the characteristic addresses
   of $\adds'$. Since $\r$ and $\addrt$ are also the characteristic
   addresses of $\kappa_0$, we conclude that
   $\s'=\extaddr(\kappa_0)$. \qed

\begin{cor}[Bifurcation Structure at a Parabolic Point]
  \label{cor:parabolicbifurcation}
 Suppose that $\kappa_0$ is a parabolic parameter of period $n$ and
  ray period $qn$. Let $\s := \extaddr(\kappa)$ and
  $\K(\s)=\u_1 \dots \u_{qn -1 }*$. Then
  $\partial\Hyp{\s}$ is the unique child component of $\kappa_0$; 
  if $\kappa_0$ is of satellite type, 
  then $\Hyp{\sigma^{(q-1)n}(\s)}$ is 
  the unique parent component of $\kappa_0$.
  No other hyperbolic component contains $\kappa_0$ on
   its boundary. 

  Furthermore, $\kappa_0$ is the root point of
   $\Hyp{\s}$ if and only if
   the parabolic orbit of $E_{\kappa_0}$ has essential orbit portrait
   (which is always the case if $q>1$). Otherwise,
   $\kappa_0$ is a co-root of  $\Hyp{\s}$.
\end{cor}
\proof Let $W$ be a child component of $\kappa_0$ (at least one such
  component always exists). By Theorem \ref{thm:bifurcationchild}, 
  $\extaddr(W)=\s$. By the classification of hyperbolic components
  (Proposition \ref{prop:hyperboliccomponents}), this component is unique.
  Similarly, suppose that 
  $q>1$; i.e., that $\kappa_0$ is a satellite parameter.
  By Theorem \ref{thm:bifurcationparent}, 
  $\Hyp{\sigma^{(q-1)n}(\s)}$ is the unique parent component of $\kappa_0$.
  By Proposition \ref{prop:parabolicperturbation}, no other hyperbolic
  component contains $\kappa$ on its boundary. The final statement
  was already proved in Theorem \ref{thm:bifurcationchild}.  \qed

\section{Bifurcation from 
   a Hyperbolic Component} 
  \label{sec:bifurcationfromcomponent}

 Suppose that $\kappa_0$
  is a satellite
  parabolic parameter
  with parent component
  $W$ and child component
  $V$. It follows from Theorem
  \ref{thm:bifurcationparent}
  that
  the address of the
  child component $V$
  is determined by $W$
  and the height 
  $h\in\Q\setminus\Z$
  of the internal ray
  $\IR$ landing at
  $\kappa_0$: 
  it is the
  intermediate external
  address of a curve
  situated in the 
  dynamical plane of
  $\kappa\in \IR$.
  However, this address
  is defined for
  \emph{every}
  $h\in\Q\setminus\Z$,
  regardless of whether
  we know that the
  corresponding internal
  ray has a landing point!
  
 We now
  use this idea to prove
  a combinatorial analog
  of the description of
  the structure
  of components bifurcating
  from a hyperbolic
  component $W$, as outlined
  in the introduction. 
  More precisely,
  for every $h\in\Q\setminus \Z$,
  we define the address
  $\bifaddr(W,h)$ of the
  hyperbolic component
  which ``would'' 
  bifurcate from $W$ at
  height $h$ if the
  corresponding
  bifurcation parameter 
  existed.
  Since these addresses 
  are naturally defined in 
  terms of curves in the
  dynamical plane of $\Ek$
  for $\kappa\in W$, 
  it 
  is straightforward to 
  compute their 
  itineraries with respect
  to
  $\s :=\extaddr(W)$. 
  Using this
  information, we can show
  that $\bifaddr(W,h)$
  behaves as we expect,
  working exclusively 
  within the
  combinatorial dynamical 
  plane associated to 
 $\s$.

 \subsection*{Combinatorial
 Bifurcation} 
 As discussed above,
  Theorem \ref{thm:bifurcationparent}
 suggests the following
 definition.

\begin{defn}[Combinatorial 
   Bifurcation] 
  \label{defn:combinatorialbifurcation}
 Let $W=\Hyp{\s}$ be a 
  hyperbolic component of 
   period $n$, and let
 $\kappa$ be a parameter on the internal ray $\IR$ for some
 $h\in\Q\setminus\Z$. 
 Let $\s'$ be the external address of the attracting
 dynamic ray of $\Ek$ which contains the singular value. Then we say that
 $\Hyp{\s'}$ \emph{bifurcates combinatorially} from $W$ (at height
   \index{combinatorial bifurcation}%
 $h$). We denote the address of this component by
   \[ \bifaddr(W,h) := \bifaddr(A,p/q) := 
      \bifaddr(\s,s_*,p/q) := \s' \]
\nomenclature[addrWh1]{$\bifaddr(W,h)$, $h\in\Q\setminus\Z$}{(combinatorial
      bifurcation at height $h$)}%
\nomenclature[addrApq]{$\bifaddr(A,p/q)$}{(combinatorial
      bifurcation in Sector $A$ at angle $p/q$)}%
\nomenclature[addrss*pq]{$\bifaddr(\s,s_*,p/q)$}{(combinatorial
      bifurcation in Sector $\Sec(\Hyp{\s},s_*)$ at angle $p/q$)}%
  where $A := \Sec(W,s_*)$ 
   is the sector containing $\kappa$
   and $p/q=h-\lfloor h \rfloor$ is the fractional part of $h$. 
   The component
   $\Hyp{\s'}$ is called a \emph{(combinatorial) child component} of
   $W$ and denoted by $\child$.
   \index{child component!of a hyperbolic component}%
\nomenclature[Bif]{$\child$}{(child component of $W$ at height $h$)}%
\end{defn}

 \begin{prop}[Child 
            Components] 
   \label{prop:childcomponents}
 Let $W$ be a hyperbolic 
  component of period $n$,
    and let 
   $h\in\Q\setminus\Z$. 
  Then the following hold. 
 \begin{enumerate}
  \item 
   \label{item:rootformula}
 $\displaystyle{\kappa_0
   := \Psi_W(2\pi ih) 
     = \Psi_{\child}(0)}
        \quad(\in \Ch)$.
  \item 
    \label{item:onlytwocomponents}
   If $\kappa_0\in\C$, 
    then $W$ and $\child$ 
    are the only
    hyperbolic components 
    containing $\kappa_0$ 
    on their boundaries. 
  \item If two hyperbolic 
   components have a common   
   parabolic boundary 
   point,
   then one of these 
   components is a child 
   component of the other.
  \label{item:touchingimplieschild}
  \item If $\kappa\in \IR$,
  then the dynamic rays of 
  $\Ek$ at 
  the characteristic 
  addresses of
  $\child$ land on the 
  distinguished boundary 
  cycle of 
  $\Ek$. For
  parameters in $\child$, 
  the
  dynamic root has (exact)
  period 
  $n$. 
  \label{item:childcharacteristic}
  \item No two hyperbolic 
   components have a common
   child component.
   \label{item:nocommonchildcomponents}
 \end{enumerate}
 \end{prop}
\proof
 If $\kappa_0:=\Psi_W(2\pi ih)\in\C$, then 
 $\kappa_0$ is a parabolic 
 parameter of
 period $n$ and rotation 
 number 
 $p/q:=h-
   \lfloor h \rfloor$; in 
  particular, the ray 
  period of $\kappa_0$ is 
  $qn$. Thus $W$ is
  a parent component of 
  $\kappa_0$.
  By
  Corollary
  \ref{cor:parabolicbifurcation}, $\kappa_0$ is the root of 
  $\child$, and no other hyperbolic components
  contain $\kappa_0$ on their boundaries. This proves 
  (\ref{item:onlytwocomponents});
 it also proves
 (\ref{item:rootformula})
 provided that
 $\kappa_0\in\C$. 
 Part
  (\ref{item:touchingimplieschild}) 
  was proved in
  Corollary 
  \ref{cor:parabolicbifurcation}. 

  To prove 
  (\ref{item:childcharacteristic}), 
 let $\kappa\in \IR$; also
  choose
   some parameter 
  $\kappa_1\in\child$. Let 
  $\r$ be the address of a
  periodic
   ray landing at a point
   $z_0\in \partial U_1$
   of the distinguished 
   boundary orbit of $\Ek$.
   Recall from Theorem
   \ref{thmdef:characteristicrays} 
   that the
   dynamic root of $\Ek$ is
   the only periodic point 
   on $\partial U_1$ which
   has an essential orbit 
   portrait. Thus the
   period of $\r$ is
   the same as the period
   of its landing point,
   which is $qn$.
   For 
   $1\leq j\leq n$,
   define
    \[ A_j := 
     \{\sigma^{mn + j - 1}(\r):0\leq m < q \}. \] 
   Then the landing points 
  of the dynamic rays at
    the addresses in $A_j$
    are those points of
    the orbit of $z_0$ 
    which belong to 
    $\partial U_j$. 
    Similarly as in the
    proof of Theorem 
   \ref{thm:bifurcationparent},
    we can connect these
    landing points to the periodic point $a_j$ using attracting dynamic rays
    such that the resulting ``extended rays'' $\wt{g}_{\sigma^j(\r)}$ do
    not intersect except in their endpoints.
    Again, it follows that
    all addresses in $A_1$ have the same itinerary under 
    $\wt{\s} := \bifaddr(W,h)$, and this itinerary agrees with
    that of $\wt{\s}$ in the first $qn-1$ entries. It is easy to see that
    in fact $\Orb := \{A_1,\dots,A_n\}$ is an orbit portrait for
    $E_{\kappa_1}$, and that the characteristic rays of this
    orbit portrait both belong to $A_1$. (All that needs to be checked is that
    rays at addresses in different $A_j$ cannot have the same landing point,
    which follows readily from the way that these rays are permuted by
    $\sigma^n$.) It now follows from Lemma \ref{lem:uniqueaddress} that
    these characteristic rays are in fact the characteristic rays of
    $\wt{\s}$, and (\ref{item:childcharacteristic}) is proved.

  To complete the proof of (\ref{item:rootformula}), 
  suppose now that
   $\kappa_0' := 
    \Psi_{\child}(0)\in \C$. 
  Then by 
   (\ref{item:childcharacteristic}) and Theorem
   \ref{thm:bifurcationchild}, this parabolic
   parameter has period $n$ and ray period $qn$. Thus Corollary 
   \ref{cor:parabolicbifurcation} implies that $\kappa_0'=\kappa_0$.   

  Finally, let us prove (\ref{item:nocommonchildcomponents}). 
  Suppose that $W'=\Hyp{\s'}$, of period $m$, 
  is a child component of 
  $W$. 
  By 
  (\ref{item:childcharacteristic}), 
  the period $n$ of $W$
  is uniquely determined by
  $W'$.
  By the definition of 
  child components,
 \[ \extaddr(W)=\sigma^{m-n}(\s'); \]
  so $\extaddr(W)$, and
  thus $W$, is uniquely
  determined by $W'$.
  \qed

\subsection*{Analysis of 
      Itineraries} 
 We can now describe the 
  addresses
  of child components of a given hyperbolic component $W=\Hyp{\s}$ 
  in terms of their
  itineraries under $\s$. This will be the key to understanding their behavior.
  
 \begin{lem}[Itineraries of Bifurcation Addresses] 
   \label{lem:itinerarybifurcation}
  Let $\s$ be an 
  intermediate external
  address of length $n$
  with 
  kneading sequence
   $\u_1\dots \u_{n-1}*$. 
  Furthermore, let $s_*\in\Z$ 
  (resp.\ 
   $s_*\in\Z+\onehalf$ if 
   $\s=\infty$) and 
   $\alpha=p/q
     \in\Q\cap(0,1)$. 
  If $\s\neq\infty$,
  then 
  $\s':=\bifaddr(\s,s_*,p/q)$ 
  is the unique 
  intermediate address 
  which satisfies
   \[ \itin_{\s}(\s') =
\u_1 \dots \u_{n-1} \m_{1} 
\u_1 \dots \u_{n-1} \m_{2} 
     \dots
\u_1 \dots \u_{n-1} 
                 \m_{q-1} 
\u_1 \dots \u_{n-1} \m_q, \]
  where 
  $\m_q = *$,
  $\m_{q-1}=
    \bdyit{s_*}{s_*-1}$ and
 \[ \m_{j} = 
    \begin{cases}
     s_*
      & \text{if } 
          j\alpha\in 
          [1-\alpha,0]
          \quad (\mod 1) \\
     s_*-1
      & \text{otherwise}
    \end{cases} \]
  for $j=1,\dots, q-2$.

 If $\s=\infty$, then
  the above is still 
  true, except that
  $\s_*$ is replaced
  by $\s_*+1/2$ in the
  definition of the
  $m_j$.
 \end{lem}
 \proof  As in Definition
 \ref{defn:combinatorialbifurcation}, 
 let $\kappa$ be
 a parameter on the
 internal ray at angle
 $\alpha$ in the sector
 $\Sec(\Hyp{\s},s_*)$
and let
    $\gamma$ be the 
attracting dynamic ray of 
$\Ek$ which contains the
singular value. Then
$\s' = \extaddr(\gamma)$.
Set 
$\addut := \itin_{\s}(\s')$. 
    
Since $\Ek^{j-1}(\kappa)$ belongs to the Fatou component
    $U_j$, it is clear that
$\ut_{kn+j}=\u_j$ for $k\in\{0\dots q-1\}$ and
    $j\in\{1\dots n-1\}$. So we only need to show that
    the values $\m_j:= \ut_{jn}$ have the stated form.

 For $j=1,\dots,q$, let us set
  $\gamma_j := \Ek^{jn-1}(\gamma)$. 
  Note that $\gamma_q$ is 
  the
  attracting dynamic ray to
 $-\infty$ in $U_0$.
 Thus $\m_q=*$ and 
 $\Ek(\gamma_{q-1})$ is
 the principal
 attracting ray. Thus 
 $\extaddr(\gamma_{q-1})= 
  s_*\s$ by definition of 
 $s_*$; so
 $\m_{q-1} = 
  \bdyit{s_*}{s_*-1}$. 
 Attracting dynamic rays 
 do not intersect each other or their
 $2\pi i$-translates; thus
 any other entry
 $\m_{j}$ is either $s_*$ 
 or $s_*-1$, depending on 
 whether 
 $\gamma_j$ is above or 
 below $\gamma_{q-1}$. 
 Since $\Ek^n$ permutes
 the curves $\gamma_j$ 
 cyclically with rotation
 number $\alpha$, 
 $\gamma_j$ is above
 $\gamma_{q-1}$ if and
 only if $j\alpha\in
 [1-\alpha,1]\ (\mod 1)$.
 \qed

 In this and the
  following section,
  we will frequently
  be concerned with the
  question when two addresses
  whose itineraries
  coincide must in fact 
  be the same. Let us
  therefore state the
  following simple
  fact
  for further reference.

 \begin{observation}[%
  Agreeing Itineraries]
  \label{obs:agreeingaddresses}
  Let $\s\in\Sequb$ and
  $k > 0$.
  Let $\r^1,\r^2$ be
  addresses with
  $\sigma^k(\r^1)
   \leq \sigma^k(\r^2)$
   whose itineraries
  (under $\s$)
  agree in the first $k$
  entries.
  Suppose furthermore
  that, for
  $j=0,\dots,k-1$, 
  $\sigma^j(\s)\notin
    [\sigma^k(\r^1),
     \sigma^k(\r^2)]$.

 Then $\r^1\leq\r^2$ and
  $\sigma^k$ maps
  the interval
 $[\r^1,\r^2]$ bijectively
  onto
   $[\sigma^k(\r^1),
    \sigma^k(\r^2)]$;
  in other words, 
  the addresses
  $\r^1$ and $\r^2$ agree
  in the first $k$ entries.
 \end{observation}
\proof
 Note that it is sufficient
  to deal with the
  case $k=1$; the general case
  follows by induction.  

 So suppose that
  $\r^1$ and $\r^2$ are as in
  the statement, with $k=1$,
  and let $I$ denote the
  interval in $\Sequ$ bounded
  by $\r^1$ and $\r^2$. Since
  both addresses have the
  same first itinerary entry, 
  $\sigma$ maps $I$ bijectively
  either to $[\sigma(\r^1),\sigma(\r^2)]$ or to
  $\Sequb\setminus \bigl(\sigma(\r^1),\sigma(\r^2)\bigr)$. 
  Since $\s\notin \sigma(I)$, it follows from the
  hypotheses that the former must be the case, as required. \qed

\begin{suppress}
 Consider an address
  $\addrt^2$ which agrees
  with $\r^1$ in the first
  $k$ entries and with
  $\r^2$ in the remaining
  entries.
  Then clearly
  $\sigma^{k-j}$ maps
  $[\sigma^{j}(\r^1),
    \sigma^{j}(\addrt^2)]$
  bijectively to
  $[\sigma^{k}(\r^1),
    \sigma^{k}(\r^2)]$
  for
  $j=0,\dots,k$. 

 Thus by hypothesis
  $\s\notin
  [\sigma^j(\r^1),
    \sigma^j(\addrt^2)]$ for
  $j=1,\dots,k$. 
  So the first 
  $k$ itinerary entries 
  of $\addrt^2$ and 
  $\r^1$, and thus of 
  $\addrt^2$ and $\r^2$, 
  agree. Since 
  $\sigma^k(\r^2)= 
   \sigma^k(\addrt^2)$, 
  this implies 
  $\addrt^2=\r^2$ 
  by Lemma 
  \ref{lem:itineraries}. \qed 
\end{suppress}

 We can now use this 
  observation to relate 
  the behavior of 
  the itineraries of 
  the addresses 
  $\bifaddr(W,h)$ 
  to that of 
  the addresses 
  themselves. 

 \begin{prop}[Monotonicity
           of Itineraries] 
    \label{prop:itinerarymonotonicity}
  Let $W=\Hyp{\s}$ be a 
 hyperbolic component of 
 period $n\geq 2$,
   with kneading sequence 
 $\K(W)=\u_1\dots 
           \u_{n-1} *$.

  \begin{enumerate}
\item 
    Let $I\subset\Sequ$ denote either 
    the interval 
  $[\s^-,\s]$ or the 
  interval
   $[\s,\s^+]$
  (where 
 $\s^- < \s^+$
 are, as usual, the
 characteristic addresses
 of $\s$). 
 Suppose that $\r^1,\r^2\in I\setminus\{\s\}$ 
    have the following property:
    if $\ell\in\{1,2\}$ and $j\geq 0$ are such that 
    $\sigma^{jn}(\r^{\ell})$ is defined, then 
    $\sigma^{jn}(\r^{\ell})\in I$ and
    $\itin_{\s}(\sigma^{jn}(\r^{\ell}))$ starts with
    $\u_1\dots \u_{n-1}$.
    Then 
     \[ \r^1 \leq \r^2\ \Longleftrightarrow\ \itin_{\s}(\r^1) 
                                          \leq \itin_{\s}(\r^2).\]
   \label{item:itinerarymonotonicity}

   \item Let 
     $h=\frac{p}{q}\in \Q
         \setminus\Z$ 
     and $\r := 
     \bifaddr(W,h)$. If 
     $h>0$ (resp.\ if $h<0)$, then
     $\sigma^{jn}(\r)\in (\s^-,\s]$ (resp.\
     $\sigma^{jn}(\r)\in [\s,\s^+)$) for all
  $0\leq j < q$.
     \label{item:alliteratesabove}

  \end{enumerate}
 \end{prop}
 \begin{remark} We should remark on the lexicographic 
  order of itineraries referred to in
  (\ref{item:itinerarymonotonicity}). There is a natural order between
  integer itinerary entries and boundary symbols: 
  $\m < \itj$ if and only if $\m\leq \j-1$. However, it is not clear
  how the symbol $*$ should fit into this order. 
  We will
  fix the convention
  that the symbol $*$ is 
  incomparable to any other
  itinerary entry, which gives our claim
  the strongest possible meaning.

 (In fact, this is not
  relevant for our considerations:
  the itineraries of any two
  addresses $\r^1,\r^2\in\Sequ$
  will be comparable unless
  at least one of them is
  an intermediate external
  address which is not a
  preimage of $\s$. Clearly
  this cannot happen in our case.)
\end{remark}
 \proof 
  To prove item 
   (\ref{item:itinerarymonotonicity}); 
  let us fix our
   ideas by supposing
   that $I=[\s^-,\s]$. 
   As already remarked above, 
   $\itin_{\s}(\r^1)$ and 
   $\itin_{\s}(\r^2)$ are 
   comparable.
\begin{suppress}
   indeed, this is trivial 
   if they
   are both infinite 
   external addresses. 
   If at least one of the 
   addresses is 
   intermediate, then it 
   is an iterated
   preimage of $\s$ and 
   thus its itinerary 
   contains a boundary
   symbol. It follows 
   easily that the two 
   itineraries are 
   comparable. 
\end{suppress}
   Note also that the orbit of
   $\s$ does not enter
   the interval $I$ by the definition
   of characteristic addresses
   (this will
   enable us to apply the
   previous observation.)    

  Suppose first that
   $\itin_{\s}(\r^1)=
      \itin_{\s}(\r^2)$.
  If $\r^1$ and $\r^2$ 
    are intermediate,
    then
  $\r^1=\r^2$ by Lemma 
   \ref{lem:itineraries}.
  On the other hand,
   if $\r^1$ and $\r^2$ are
   infinite,
   then for any
   $\ell\geq 0$, the hypotheses
   of
   Observation
   \ref{obs:agreeingaddresses} 
   are satisfied 
   with
   $k=n\ell$. Thus
   the first $n\ell$
   entries of $\r^1$ and $\r^2$ 
   agree; since $\ell$ is arbitrary,
   this means that
   $\r^1=\r^2$.
  
  So now suppose that 
  $\itin_{\s}(\r^1)\neq 
   \itin_{\s}(\r^2)$, say
   $\itin_{\s}(\r^1)<
    \itin_{\s}(\r^2)$. 
  Let $\ell\geq 1$ be
   such that these
   itineraries first
   differ in the
   $n\ell$-th entry. Then
   $\r^1$ and $\r^2$
   agree in the first
   $(\ell-1)n$ entries by
   Observation 
   \ref{obs:agreeingaddresses}.
  Thus we may suppose, by
   passing to the
   $(\ell-1)n$-th iterates,
   that $\ell=1$.

 Then
  $\sigma^{n-1}(\r^1) < 
   \sigma^{n-1}(\r^2)$.
 Since $\sigma^{n-1}$ 
  preserves the circular 
  order of
  $\s$, $\r^1$ and $\r^2$,
  and since 
  $\r^1,\r^2 \leq \s$, it 
  follows that
  $\r^1 < \r^2$, 
 as required.

\smallskip

 Now let us prove 
  (\ref{item:alliteratesabove}).
  We again
  fix our ideas
  by supposing that 
  $h=p/q>0$. Let $\u_n$ 
  denote the common $n$-th 
   itinerary
   entry of $\s^-$
   and $\s^+$ with respect to
   $\s$. Then
  the $n$-th entries of
  the addresses
  $\s^-$ and $\s^+$ are
  $\u_n+1$ and $\u_n$,
  respectively.

  Since $h>0$, we thus
  have 
  $s_*\geq \u_n+1$ by Theorem
  \ref{thm:labellingmap}. So for
  $j=1,\dots,q-1$,
  the $nj$-th itinerary
  entries $\m_j$ of
  $\r$ satisfy 
   $\m_j\geq \u_n$ 
   by
   Lemma 
   \ref{lem:itinerarybifurcation}.
    We claim that, 
 for every $j=0,\dots,q-2$,
  \[ \sigma^{(j+1)n}(\r)\in (\s^-,\s] \ \Longrightarrow
      \ \sigma^{jn}(\r)\in (\s^-,\s]. \]
   Since 
  $\sigma^{(q-1)n}(\r)=\s$,
 part (\ref{item:alliteratesabove})
  then follows by 
 induction. 

 To prove the claim,
  suppose that 
 $\sigma^{(j+1)n}(\r)\in 
    (\s^-,\s]$. The first
  $n-1$ itinerary
  entries of
 $\s^-, \sigma^{jn}(\r)$ 
 and $\s$ are the same,
  so
 $\sigma^{n-1}$ preserves 
 the circular order of 
 these addresses
 by Observation
\ref{obs:orderpreserving},
   Thus it is sufficient to
 show that
  $\sigma^{(j+1)n-1}(\r)>
   \sigma^{n-1}(\s^-)$. 
 This is trivial if
   $\m_j>\u_n$, and 
 follows from the fact that
 $\sigma^{(j+1)n}(\r)\in 
 (\s^-,\s]$ if $\m_j=\u_n$.
 \qed

 \begin{cor}[Monotonicity of $\bifaddr(W,h)$] \label{cor:monotonicity}
   Let $W$ be a hyperbolic component. Then
    the function
     $h\mapsto \bifaddr(W,h)$ is strictly increasing
     on each of the intervals $\{h>0\}$ and $\{h<0\}$.
 \end{cor}
 \proof Consider the function
   $p/q \mapsto \m_{1} \dots \m_{q-1}$, where 
   $\m_{j}$ are the numbers from Lemma
   \ref{lem:itinerarybifurcation}. 
  It is an easy exercise
  to check that
   this function is strictly increasing (with respect to 
   lexicographic order). 

 If $W\neq\Hyp{\infty}$,
  the claim follows
  from
   Proposition \ref{prop:itinerarymonotonicity} 
   (\ref{item:itinerarymonotonicity}).
 If $W=\Hyp{\infty}$,
  the claim follows
  directly since
  itineraries and
  external addresses
  coincide in this case (compare
  Remark 2 after Definition
   \ref{defn:combinatorialitinerary}.
    \qed
 
 \begin{prop}[Continuity 
   Properties] 
     \label{prop:limits}
  Let $W$ be a hyperbolic 
  component and 
   $\s:=\extaddr(W)$.
  Then
 $\displaystyle{%
 \lim_{h\to\pm\infty}
     \extaddr(W,h)}=\s$.
  Furthermore, if 
   $h_0\in\R$, then the 
   behavior of 
   $\extaddr(W,h)$ for 
  $h\to h_0$ is as follows.
  \begin{enumerate}
   \item \label{item:rational}
      If $h_0\in \Q\setminus\Z$, then 
      \begin{align*}
        \lim_{h\nearrow h_0} \bifaddr(W,h) &= \r^- \text{ and }  \\
        \lim_{h\searrow h_0} \bifaddr(W,h) &= \r^+,
      \end{align*}
    where $\r := \bifaddr(W,h_0)$.
  \item \label{item:characteristic}
     If $h_0 = 0$ and $\s\neq\infty$, then
     \begin{align*}
        \lim_{h\nearrow 0} \bifaddr(W,h) &= \s^+ \text{ and } \\\
        \lim_{h\searrow 0} \bifaddr(W,h) &= \s^-.
     \end{align*}
  \item \label{item:irrational}
   Otherwise, 
   the limit
  $\bifaddr(W,h_0):=
   \lim_{h\to h_0}
   \bifaddr(W,h)$ exists.
\nomenclature[addrWh2]{$\bifaddr(W,h)$, $h\neq 0$}{(combinatorial
               bifurcation at height $h$)}
 \end{enumerate}
 \end{prop}
 \proof Let 
  $\u_1\dots \u_{n-1}*$ be 
  the kneading sequence 
  of $W$. 
  Let $h_0=p/q\in 
      \Q\setminus\Z$, and 
 $\r\:=\bifaddr(W,h_0)$. 
 Recall that,
  for parameters 
  $\kappa\in \IRH{h_0}$, 
  the landing point 
  of the dynamic ray 
  $g_{\r^+}$ lies
  on the distinguished 
  boundary orbit of $\Ek$. 
  Since the rays landing
  on this cycle are 
  permuted cyclically with
  rotation number
  $p/q$, it follows easily
  (as in Lemma
   \ref{lem:itinerarybifurcation})
  that 
   \begin{align*}
    \itin_{\s}(\r^+)&= 
        \per{\u_1 \dots \u_{n-1} \m_1 \dots
             \u_1 \dots \u_{n-1} \m_{q-2}
             \u_1 \dots \u_{n-1} s_* 
             \u_1 \dots \u_{n-1} (s_*-1)},
    \end{align*}
  where the $\m_j$ are as 
  in Lemma 
  \ref{lem:itinerarybifurcation}. 
  It follows from Lemma 
  \ref{lem:itinerarybifurcation} 
  that the 
  limit address 
 $\lim_{h\nearrow h_0}\bifaddr(W,h)$ 
   has the same itinerary
   under $\s$ as $\r^+$. 
 If $n=1$, then the two
  addresses must trivially 
  be equal.
  Otherwise,
   this 
  follows easily from
  Observation
  \ref{obs:agreeingaddresses}.
 
 All other parts of the 
 proposition are proved 
 analogously. For each 
 part, we need to prove 
 the equality
 of two external addresses,
 at least one of which is 
 given as a
 monotone limit of 
 addresses $\bifaddr(W,h)$.
 In each case, it is
 easy to verify that the 
 corresponding itineraries 
 are equal,
 which, as above, implies 
 that the same is true for 
 the
 external addresses. \qed

 We are now ready to state 
 and prove our main theorem
 on the structure of the 
 child components of a 
 given
 hyperbolic component. If 
 $\s:=\extaddr(W)
      \neq\infty$, then
 the \emph{wake} of $W$ is
  the set
    \index{wake!of a hyperbolic component}
 $\W(W) := (\s^-,\s^+)$; 
 if $\s=\infty$, 
  $\W(W)$ is defined to 
  be all of $\Sequb$.
\nomenclature[WW]{$\W(W)$}{(wake of $W$)}
  In the following theorem, we use this only for simpler notation,
  but wakes will play an important role in the next section.

 \begin{thm}[Bifurcation 
               Structure]
    \label{thm:bifurcationstructure}
 Let $W$ be a hyperbolic 
  component and 
  $\s:=\extaddr(W)$. 
 If $n\geq 2$, then the map
 $\bifaddr(W,\cdot):
    \Q\setminus\Z\to\Sequ$
 has the following 
    properties.
  \begin{enumerate} 
   \item 
  $\bifaddr(W,\cdot)$ is 
  strictly increasing on
  $\{h>0\}$ and 
  (separately)
  on $\{h<0\}$;
  \label{item:monotonicity}
 \item
  $\bifaddr(W,\frac{p}{q})$
  is an intermediate 
  external address of
  length $qn$ (for $\frac{p}{q}$ in lowest terms); 
 \label{item:intermediateaddress}
 \item if 
  $h\in\Q\setminus\Z$ such 
  that
  $\Psi_W(2\pi i h)\in\C$, 
  then the parameter 
  $\Psi_W(2\pi i h)$ lies 
  on
  the boundary of 
  $\child= 
    \Hyp{\bifaddr(W,h)}$;
  \label{item:bifurcation}
 \item 
  $\displaystyle{%
     \cl{\W(W)}=
     \cl{\bigcup_{%
        h\in\Q\setminus\Z}
       \W(\child)}}$.
 \label{item:subwakesfillwake}
 \item 
  $\displaystyle{%
    \lim_{h\to+\infty}
       \bifaddr(W,h) = 
    \lim_{h\to-\infty} 
       \bifaddr(W,h)= \s}$;
  \label{item:addressofW}
 \item 
  $\displaystyle{%
     \lim_{h\nearrow 0} 
  \bifaddr(W,h)=\s^+}$ and
  $\displaystyle{%
     \lim_{h\searrow 0} 
  \bifaddr(W,h)=\s^-}$.
 \label{item:charaddresses}
  \end{enumerate}
 These properties uniquely determine the map $\bifaddr(W,\cdot)$, and no
  such map exists if the preferred parametrization
  $\Psi_W$ is replaced by some other conformal
  parametrization $\Psi:\H\to W$ with $\mu\circ\Psi = \exp$.

 If $n=1$, then the map
 $\extaddr(W,\cdot)$ is strictly increasing on all of
 $\Q\setminus\Z$ and satisfies properties
 (\ref{item:intermediateaddress}) to (\ref{item:subwakesfillwake}) above,
 and no other map has these properties. 
 \end{thm}
\proof Let us assume for 
 simplicity that $n>1$; 
 the proofs for the
 case $n=1$ are completely 
 analogous. 
 Property 
 (\ref{item:monotonicity}) is just the statement of
  Corollary \ref{cor:monotonicity},
 and 
 (\ref{item:intermediateaddress})
 holds by definition.
 Property 
  (\ref{item:bifurcation}) 
  is Proposition
   \ref{prop:childcomponents} (\ref{item:rootformula}). Properties
  (\ref{item:addressofW}) and (\ref{item:charaddresses}) were proved in
  Proposition \ref{prop:limits}. 

 To establish 
  (\ref{item:subwakesfillwake}), 
  note first that the inclusion ``$\supset$'' is clear. To prove
  ``$\subset$'', 
  let
  $\r\in \cl{\W(W)}$. If $\r\in\{\s^-,\s^+,\s\}$, then we are done by
  (\ref{item:addressofW}) and (\ref{item:charaddresses}). 
  Otherwise, there exists
  $h_0\in\R\setminus\{0\}$ such that 
  \[ \t^-:=\lim_{h\nearrow h_0} \bifaddr(W,h)\leq \r \leq
     \lim_{h\searrow h_0} \bifaddr(W,h) =: \t^+. \]
  If $h\notin \Q\setminus\Z$, then Proposition
  \ref{prop:limits} (\ref{item:irrational}) shows that
  $\t^-=\r=\t^+$, and we are done. Otherwise, Proposition \ref{prop:limits}
  (\ref{item:rational}) implies that 
  $\r\in \cl{\W(\bifaddr(W,h))}$, which completes the proof of
  (\ref{item:subwakesfillwake}).

 Let us now prove the uniqueness statements. Suppose that
  $a:\Q\setminus\Z\to\Sequ$ 
  also satisfies properties
  (\ref{item:monotonicity})
   to
  (\ref{item:subwakesfillwake}).
  Then, by item (\ref{item:bifurcation}) and Corollary
  \ref{cor:parabolicbifurcation}, we know that
  $a(h)=\bifaddr(W,h)$ whenever
  $\Psi_W(ih)\in\C$. By Lemma \ref{lem:internallanding},
  the set of such $h\in\R$ is open and dense,
  so there is a dense set of rationals on which
  $a$ and $\bifaddr(W,\cdot)$ agree. 
  Properties (\ref{item:monotonicity}) and (\ref{item:subwakesfillwake})
   then
   easily imply that Proposition \ref{prop:limits} is also
   true for the map $a$, which shows that both maps
   must be equal everywhere.

 Finally, if the preferred
  parametrization
  $\Psi_W$ is replaced by 
  some other 
  parametrization of $W$,
  then clearly no map 
  $a$
  which satisfies
  (\ref{item:bifurcation})
  can also satisfy 
 (\ref{item:monotonicity}).
 Indeed, such a map
  must
  agree with
  $\bifaddr(W,m+\cdot)$
  on a dense set
  for some 
  $m\in\Z\setminus\{0\}$.
  Therefore
  $a$
  is not
  monotone near 
  $m$ by 
  (\ref{item:charaddresses}).
 \qed

 We note the following 
 consequence of the 
 previous 
 theorem for reference.

 \begin{cor}[Subwakes Fill Wake] \label{cor:subwakesfillwake}
  Suppose that $\s\in\W(W)\setminus\{\extaddr(W)\}$. 
  Then there exists 
   a unique $h\in\R$
  such that one 
  (and only one)
  of the following hold.
  \begin{enumerate}
   \item \label{item:insubwake}
      $h\in \Q\setminus\Z$ and $\s\in\W(\child)$,
   \item $h\in \Q\setminus\Z$ and $\s$ is a characteristic address of
         $\child$, or
   \item $h\in 
    (\R\setminus\Q)\cup
    (\Z\setminus\{0\})$ 
   and $\s=\bifaddr(W,h)$.
  \end{enumerate}
  In particular, if $\s$ is unbounded, then Property (\ref{item:insubwake})
   holds. \qed
 \end{cor}

\section{Internal Addresses}
 \label{sec:internaladdresses}

\begin{parameterarguments}
 In this section, we describe the global bifurcation structure
  of hyperbolic components. The basic question we are now interested in
  is as follows: given two hyperbolic components $V$ and
  $W$, one of
  which is contained in the wake of the other, how can we determine 
  which bifurcations occur ``between'' $V$ and $W$? 

 We will begin by dividing
  up the wake
  of a hyperbolic component $W$ into \emph{sector wakes} $\W(A)$ (one
  for every sector $A$ of $W$) such that every child component
  of $W$ belongs to the wake of the sector from which it bifurcates. 

 The wakes of two
  adjacent sectors are separated by a \emph{sector boundary}, 
  i.e.\ a periodic address which one should think of as the address
  of the parameter ray 
  landing at the 
  corresponding co-root of 
  $W$
  (except that, for now, 
   we do not know whether 
   this co-root parameter
   actually exists). 
  Every sector $A$ has a natural associated
   \emph{sector kneading sequence} $\K(A)$ (Definition
   \ref{deflem:kneadingentries}). We give a simple 
   description of how kneading sequences
   depend on the bifurcation structure of
   parameter space (Theorem \ref{thm:determiningbifurcationorder}).
  Finally, we will introduce \emph{internal addresses}, which
   organize the information encoded in kneading sequences in a
   ``human-readable'' way. 
\end{parameterarguments}

\begin{suppress}
 In this section, we will finally describe the global bifurcation structure
  of hyperbolic components. In order to give this decription, we first
  transfer the partition of a hyperbolic component into its sectors to
  the combinatorial level via the introduction of \emph{sector wakes}. 
  To every sector, we then associate a natural \emph{sector kneading sequence}.
  The main result of this section 
  (Theorem \ref{thm:determiningbifurcationorder}) then explains how the
  bifurcation
  structure of hyperbolic components and sectors
  can be read off from these kneading sequences.

 We also introduce \emph{internal addresses}, which contain the same 
  information as kneading sequences, but directly specify
  the position of
  a hyperbolic component within exponential parameter space in a 
  human-readable form. As a corollary, we give a simple necessary and
  sufficient condition for an attracting exponential map to have
  infinitely many essential periodic orbits. 
\end{suppress} 

\begin{parameterarguments}

\subsection*{Wakes of sectors and combinatoril arcs}
 Let $W$ be a hyperbolic component, let $h_0\in\Z$ and consider the
  sector
  \[ A = 
     Sec(W,h_0+frac{1}{2})
   =
    \Psi_W\bigl(\{z\in\H:
      \im z \in (2\pi h_0,
          2\pi (h_0+1))\}
                  \bigr) \]
  of $W$. The 
  \emph{sector boundaries}
  of $A$ are defined to be
   \[ \r^- := \lim_{h\searrow h_0} \bifaddr(W,h)\quad\text{and}\quad
      \r^+ := \lim_{h\nearrow h_0+1} \bifaddr(W,h). \]
 (Note that 
  $r^-=\bifaddr(W,h_0)$
  unless $h_0=0$, in which
  case $r^-$ is the
  lower characteristic
  address of $W$, and
  similarly for $r^+$.)
The \emph{wake} of $A$ is denoted by $\W(A) := (\r^-,\r^+)$.
    \index{wake!of a sector}%
\nomenclature[WA]{$\W(A)$}{(wake of a sector $A$)}%
    \index{sector boundary}%
  If $\r$ is a sector boundary of a sector of $A$, we also say that
   $\r$ is a sector boundary of $W$.
 
   Armed with this concept, we can introduce a natural
  (partial) order on sectors and hyperbolic components.

 \begin{defn}[%
 Combinatorial Arcs]
  If $A,B$ are hyperbolic components or sectors, we write 
  $A\prec B$ if $\W(A)\supset \W(B)$. The \emph{combinatorial arc}
  $[A,B]$ is the set of all
 hyperbolic components
 or sectors $C$ such that
  $A\prec C \prec B$. 
\nomenclature{$\prec$}{(combinatorial order)}%
\nomenclature{$[A,B]$}{(combinatorial arc)}%
\index{combinatorial arc}%

  Similarly, if $\s\in \cl{\W(A)}$, then the combinatorial arc
  $[A,\s)$ is the set of all $C$ with $A\prec C$ and
  $\s\in \cl{\W(C)}$. Note that $[A,\s)$ is linearly ordered by
  $\prec$. 
\nomenclature{$[A,\s)$}{(combinatorial arc)}%
 \end{defn}
 \begin{remark}
  We will often also consider
  open or half-open combinatorial arcs $(A,B)$, $[A,B)$ or $(A,B]$, in which
  one or both of the endpoints are excluded.
 \end{remark}

\subsection*{Sector Boundaries and Kneading Sequences}

\begin{lem}[Sector Boundaries and Itineraries]
 \label{lem:sectorboundaries}
  Let $W$ be a hyperbolic component of period $n$ and 
   kneading sequence $\K(W)=\u_1\dots\u_{n-1} *$. 
   Set $\s := \extaddr(W)$ and let
   $\r\in\Sequ$. Then the following
   are equivalent.
  \begin{enumerate}
   \item $\r$ is a sector boundary of $W$; \label{item:sectorboundary}
   \item $\itin_{\r}(\s)=\K(W)$ and        \label{item:itincriterion}
     $\itin_{\s}(\r)=\per{\u_1\dots \u_{n-1} \m}$
   for some $\m\in\Z$.
 \end{enumerate}
  Furthermore, every sector boundary $\r$ satisfies
  \begin{equation}
   \sigma^j(\r)\notin \cl{\W(W)} \label{eqn:iteratesoutsidewake}
  \end{equation}
 for $j=1,\dots,n-1$; in particular, $\r$ has (exact)
  period $n$.
\end{lem}
\proof
 It follows easily from Lemma
  \ref{lem:itinerarybifurcation} that for every
  $\m\in\Z$ there is a sector boundary $\r$ with
  $\itin_{\s}(\r)=\per{\u_1\dots \u_{n-1} \m}$, and every sector boundary
  $\r$
  has an itinerary of this form. Furthermore, every 
 sector boundary $\r$
  belongs to $\cl{\W(W)}$, and thus satisfies
  $\itin_{\r}(\s)=\K(W)$ 
 by Observation
 \ref{obs:itinchange}
 (recall that the iterates of $\s$ do not
  enter $\W(W)$). By
  Theorem \ref{thmdef:characteristicrays},
  the characteristic addresses of $W$ have period $n$.
 In particular,
  (\ref{item:sectorboundary})
  implies (\ref{item:itincriterion}).

 Now let $\kappa\in W$ and let $\r$ be a sector boundary of $W$
  which is not a characteristic address. A simple hyperbolic
  expansion argument (compare \cite[Proof of Theorem 6.2]{expattracting}
  or \cite[Theorem 4.2.4]{thesis}) shows that the dynamic ray
  $g_{\r}$ lands on the boundary of the characteristic Fatou component
  $U_1$ of $\Ek$. By the definition of characteristic addresses,
  this landing point $w$ is separated from the rest of its orbit
  by the characteristic rays of $\Ek$, and its orbit portrait
  is not essential. In particular, (\ref{eqn:iteratesoutsidewake})
  holds (and $\r$ and $w$ both have exact period $n$).

 To prove that 
  (\ref{item:itincriterion})
 implies (\ref{item:sectorboundary}), 
 suppose that $\r'$ is an 
  address which is not a sector boundary
  and has itinerary
  $\itin_{\s}(\r')=
   \per{\u_1 \dots \u_{n-1} \m}$.
  Note that
  $\r'$ is necessarily periodic by Lemma
  \ref{lem:itineraries}.
  There is a sector boundary $\r$ which has
  the same itinerary as 
  $\r'$; by Theorem \ref{thm:expper}, the dynamic rays
  $g_{\r}$ and $g_{\r'}$ have a common landing point.
  By the above, this implies that
  $\r$ is a characteristic address of
  $W$, and $\r'\notin\cl{\W(W)}$.  
  By replacing
  $\r$ with the other
  characteristic address
  of $W$,
  if necessary, we may
  suppose that $\r$ and
  $\r'$ do not enclose
  $\s$. It then follows
  from Observation
  \ref{obs:agreeingaddresses}
  that they 
  must enclose
  a forward iterate of
  $\s$. Thus
  $\itin_{\r'}(\s)\neq \itin_{\r}(\s) = \K(\s)$ by Observation
  \ref{obs:itinchange}, as required.
 \qed

\begin{deflem}[Kneading Entries] \label{deflem:kneadingentries}
 Let $W$ be a hyperbolic component of period $n$ and kneading sequence
  $\K(W)=\u_1\dots\u_{n-1}*$. 
 \begin{enumerate}
  \item Let $A$ be a sector of $W$ with sector boundaries
    $\r^-<\r^+$. Then there exists a number $\u(A)\in\Z$ (the
    \emph{kneading entry} of $A$) 
\nomenclature[uA]{$\u(A)$}{(kneading entry of $A$)}%
     \index{kneading entry}%
   such that
    \[ \K^+(\r^-) = \K^-(\r^+) = \per{\u_1\dots\u_{n-1} \u(A)}
                =: \K(A). \]
   The sequence $\K(A)$ is called the \emph{kneading sequence}
   of the sector $A$. \label{item:sectorkneading}
    \index{kneading sequence!of a sector}%
  \item If $n\geq 2$, then there exists a number 
        $\u(W)\in\{\u_1,\dots,\u_{n-1}\}$
    (the \emph{forbidden kneading entry of $W$})
\nomenclature[uW]{$\u(W)$}{(forbidden kneading entry)}%
\index{forbidden kneading sequence/entry}%
\index{kneading entry!forbidden}%
      such that
     \[ \K^-(\s^-) = \K^+(\s^+) =\per{\u_1\dots\u_{n-1}\u(W)}=:\KS(W) \]
     (where $\s^-$ and $\s^+$ are the characteristic 
  addresses of $W$). 
     The sequence $\KS(W)$ is called the \emph{forbidden kneading sequence
     of $W$}. 
\nomenclature[KW]{$\KS(W)$}{(forbidden kneading sequence)}%
\index{kneading sequence!forbidden}%

    The kneading entries of the sectors directly above and below 
     the central internal ray of $W$ are $\u(W)-1$ and
     $\u(W)+1$, respectively. If $A$ and $B$ are any other two
     adjacent sectors, with $A$ above $B$, the kneading entries satisfy
     $\u(A)=\u(B)+1$ (compare Figure \ref{fig:sectorkneadingsequences}).

    In particular, no two sectors have the same kneading entry and
     $\u(W)$ is the unique integer which is not assumed as the kneading
     entry of some sector of $W$.
  \item In the period one case, every integer $\u$ is realized as 
     the kneading entry of a sector of $W$, namely the sector
     at imaginary parts between $2\pi\u$ and $2\pi(\u+1))$. The period
     one component thus has no forbidden kneading sequence.
  \end{enumerate}
 For $\m\neq \u(W)$, we denote the unique sector $A$ satisfying
   $\u(A)=\m$ by $\Sec(W,\m)$.
\nomenclature[SecWm]{$\Sec(W,\m)$}{(sector with kneading entry $\m$)}%
\end{deflem}
\proof In the period one case, the sector boundaries are the addresses
  $\per{m}$ with $m\in\Z$, and the claims are trivial. So suppose that
  $n\geq 2$; to prove
  (\ref{item:sectorkneading}),
  let us fix our ideas by 
  supposing that
  $\r^-<\r^+<\s := 
   \extaddr(W)$. It
  follows from Lemma
  \ref{lem:itinerarybifurcation} that
   \[\itin_{\s}(\r^-)=\per{\u_1\dots\u_{n-1}(s_*-1)}
   \quad\text{and}\quad
     \itin_{\s}(\r^+)=\per{\u_1\dots\u_{n-1}s_*}, 
   \]
 where $s_* = s_*(A)$ is 
 the
 sector label of $A$.
  By Observation \ref{obs:itinchange} and
   Lemma \ref{lem:sectorboundaries}, we see that
   $\K^+(\r^-)=\itin_{\s}(\r^-)$ and
   $\K^+(\r^+)=\itin_{\s}(\r^+)$. Since $\r^+$ is periodic of
   period $n$, we thus have
   \[ \K^+(\r^-)=\per{\u_1\dots\u_{n-1}
     (s_*-1)}=\K^-(\r^+), \]
   so $\u(A)=s_*-1$ is the
   desired kneading entry.

 We already know from Lemma \ref{lem:charraypairs} that a number
  $\u(W)\in\Z$ with the required property exists. By Observation
  \ref{obs:agreeingaddresses}, the $(n-1)$-th iterates of 
  $\s^-$ and $\s^+$ must inclose an iterate of $\s$. Thus the
  entry $\u(W)$ must occur in $\K(\s)$; i.e.\
  $\u(W)\in\{\u_1,\dots,\u_{n-1}\}$. 
  The remainder of the
  statement follows from 
  the fact that, for any 
  address $\r$ of
  period $n$, the $n$-th 
  entries of $\K^-(\r)$ 
  and $\K^+(\r)$ differ
  exactly by one. \qed

\begin{figure}
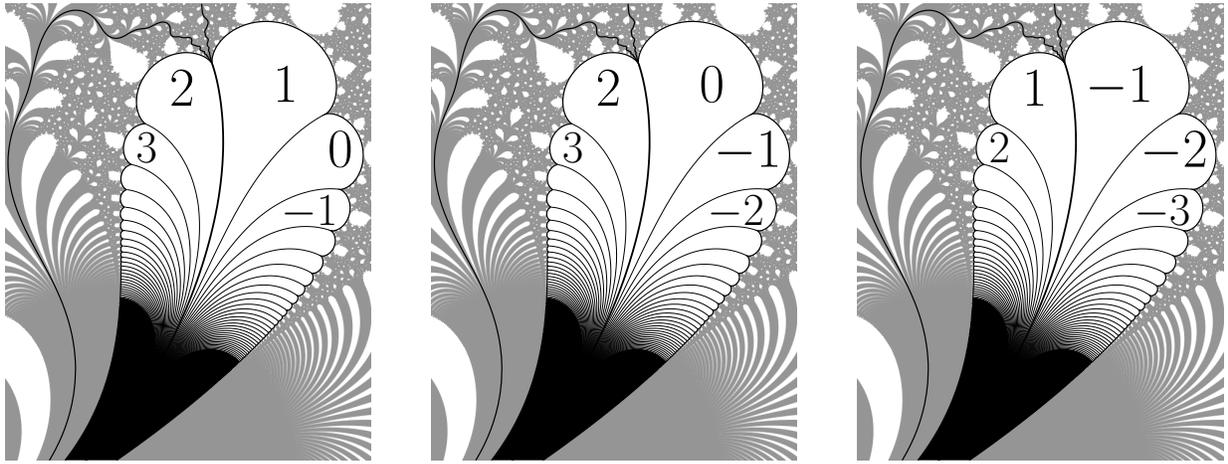

\subfigure[Sector Labels]{\input{sector_labels.epstex}%
             \label{fig:seclabels}}\hfill
\subfigure[Kneading Entries]{\input{sector_kneading.epstex}%
             \label{fig:sectorkneadingsequences}}\hfill
\subfigure[Sector Numbers]{\input{sector_numbers.epstex}%
             \label{fig:sectornumbers}}
 \caption{Three ways to label sectors, illustrated for the component
  at address $\s=0110\frac{1}{2}\infty$, with
  $\s^-=\per{011002}$ and $\s^+=\per{020101}$. The central internal ray of
  $\Hyp{\s}$ is
  emphasized.}
\end{figure}

\subsection*{Periodic addresses and hyperbolic components} 
 The last ingredient we 
  require for our analysis of bifurcation
  structure is the fact that every periodic address is the sector
  boundary of some hyperbolic component. 

\begin{prop}[Periodic Rays and Intermediate Addresses]
   \label{prop:periodicintermediate}
  Let $\r\in\Sequ$ be periodic of period $n$. 
   Then there exists an intermediate
   external address $\s$ of length $n$ such that $\r$ is a sector boundary of 
   $W:= \Hyp{\s}$. If there exists $\wt{\r}$ such that
   $\langle\r,\wt{\r}\rangle$ is a characteristic ray pair, then
   $\r$ and $\wt{\r}$ are the characteristic addresses of $W$.
\end{prop}
\proof Let $\per{\u_1\dots\u_{n-1}\bdyit{\u_n+1}{\u_n}}$ be the
  kneading sequence of $r$ and
   choose a parameter $\kappa$ for which the singular value lies on
   the dynamic ray $g_{\r}$.   
   (It is well-known --- and easy to see ---
    that such parameters exist for every periodic
    address; see 
    e.g.\ \cite{dghnew1}.\footnote{%
   It is possible to 
    formulate the following
    argument completely 
    combinatorially, 
    without appeal to the 
    existence of such 
    parameters.}
   In fact, this is true
    for every exponentially bounded external address, and these
    parameters form a corresponding
    \emph{parameter ray}; 
    see
   \cite{markus,markusdierk,markuslassedierk}.)
  Consider the piece of $g_{\r}$ which 
   connects the singular value with
   $+\infty$. The dynamic ray    
   $g_{\sigma^{n-1}(\r)}:(0,\infty)\to\C$ is a preimage component of
   this piece, and thus
   tends
   to $+\infty$ for
   $t\to\infty$,
   and to $-\infty$ for 
   $t\to 0$. 
   The ray 
   $g_{\sigma^{n-1}}(\r)$ and its translates thus 
  cut the 
   dynamic plane into 
  countably many domains
  which we call
  ``strips''
  (similarly as in 
  the definition
   of itineraries for 
  attracting parameters).

  Since 
  $g_{\sigma^{n-1}(\r)}
      =\Ek^{n-1}(g_{\r})$
  tends to $-\infty$ as
  $t\to 0$,
   the curve 
  $g_{\r}$ has an intermediate
   external
   address $\s$ of length 
  $n$ as $t\to 0$ (as well as
  the usual address $\r$ as
  $t\to\infty$). 
  Dynamic rays do not
   intersect, so both ends of the ray 
   $g_{\r}$ tend to $\infty$ within the same strip.
 In other words, the itineraries of
   $\r$ and $\s$ (with respect to $\r$) coincide in the first entry. 
   We can apply the
   same argument to $g_{\sigma(\r)},\dots, g_{\sigma^{n-2}(\r)}$, and 
   conclude that  $\itin_{\r}(\r)$
   and $\itin_{\r}(\s)$ coincide in the first $n-1$ entries; i.e.,
   $\itin_{\r}(\s)=\u_1\dots\u_{n-1}*$. To fix our ideas,
   let us suppose that $\r<\s$. 

 Note that for
  $j=0,\dots,n-1$, we have
  $\sigma^j(\s),
   \sigma^j(\r)\notin
  (\r,\s)$. Otherwise,
  we could choose a
  minimal such $j$, and
  Observation
  \ref{obs:agreeingaddresses}
  would imply that
  $\sigma^j$ maps 
  $[\r,\s]$ to
  $[\sigma^j(\r),
    \sigma^j(\s)]\subset
    [\r,\s]$, which is
  clearly impossible since
  no interval is invariant
  under the shift. 
 By Observation 
  \ref{obs:itinchange}, 
  it follows that
   \[ \K(\s) = 
      \itin_{\r}(\s) = \u_1\dots\u_{n-1}* \quad\text{and}\quad
      \itin_{\s}(\r)= 
         \K^-(\r) = \per{\u_1\dots\u_n}. \]
   By Lemma 
  \ref{lem:sectorboundaries},
   $\r$ is a sector boundary of $W:=\Hyp{\s}$. 

 \smallskip

 Finally, suppose that
  $\wt{\r}\in\Sequ$ is another address of period
  $n$ such that $\langle\r,\wt{\r}\rangle$ is a
  characteristic ray pair. Then
   $\itin^{-}_{\r}(\wt{\r}) = \per{\u_1\dots\u_n}$ by
  Lemma 
  \ref{lem:charraypairs}
   (\ref{eqn:itineraries}),
  and it follows that
  the cyclic order of $\r$, $\wt{\r}$ and
  $\s$ is preserved by $\sigma^{n-1}$. Since
  $\sigma^{n-1}(\wt{\r}) < \sigma^{n-1}(\wt(\r))$ and
  $\sigma^{n-1}(\s)=\infty$,
  this means that
  $\r<\s<\wt{\r}$. By Lemma
 \ref{lem:uniqueaddress},
  $\r$ and $\wt{\r}$ are the characteristic addresses of $\s$. 
 \qed

 \subsection*{Evolution of kneading sequences}
\begin{lem}[Nested 
  Components Have Different
    Kneading Sequences]
  \label{lem:noequalkneadingsequences}
 Let $V$ and $W$ be two hyperbolic components with
  $W\prec V$. Then $\K(V)\neq \K(W)$. 
\end{lem}
\proof This is a direct corollary of
 Lemma 
 \ref{lem:uniqueaddress}, 
 ``(\ref{item:equalitin})
    $\Rightarrow$ 
    \ref{item:charaddress})''.
 \qed
 
\begin{lem}[Kneading Sequences in Sector Wakes]
   \label{lem:kneadingandwakes}
  Let $W$ be a hyperbolic component, and let $A$ be a 
   sector of $W$.
  \begin{enumerate}
   \item Let $\s\in\W(A)$ 
  and suppose that 
  $m\geq 1$ is such that
  $\K(A)$ and $\K(\s)$ 
  have different
  $m$-th entries. Then there exists a hyperbolic component 
    $V\in [A,\s)$ of 
  period at most $m$. 
    \label{item:changeofkneadingentries}
   \item Suppose that $V\succ A$ is a hyperbolic component
    such that there are no
  components of periods 
  up to $m$ in $[A,V)$. 
    Then $\KS(V)$ agrees with $\K(A)$ in the first $m$ entries. 
    \label{item:componentkneading}
   \item Suppose that $h\in\Q\setminus\Z$. Then
    $\KS(\child)=\K(\Sec(W,h))$.
       \label{item:kneadingchild}
   \item Suppose that $h\in\R\setminus\Q$. Then
     $\K(\bifaddr(W,h))=\K(\Sec(W,h))$. 
       \label{item:kneadingbif}
  \end{enumerate}
 \end{lem}
 \begin{remark}
  Note that the addresses in (\ref{item:kneadingbif}) are non-periodic
   with periodic kneading sequences; compare
   the remark after Lemma \ref{lem:itineraries}. (It follows from
   Theorem \ref{thm:determiningbifurcationorder} below that, conversely,
   these are the only external addresses with this property.)
 \end{remark}
 \proof Let $K$ consist of all intermediate
  external addresses of length $\leq m$ in $\W(A)$. Then
  $K$ is a closed,
  and hence compact,
  subset of $\Sequb$,
  and $\U:=\{\W(\r):\r\in K\}$ is an open cover 
  of $K$. Since wakes of hyperbolic components are either
  nested or disjoint, it follows that $U := \bigcup \U$ has finitely
  many connected components, each of which is an element of $\U$. 

  By Proposition \ref{prop:periodicintermediate}, the set
   $M := \W(A)\setminus\cl{U}$ contains no periodic addresses of period
   $\leq m$; thus each of
   the first $m$ entries of $\K(\s)$ is locally constant when considered
   as a function
   of $\s\in M$. On the other hand, if $(\r^-,\r^+)\in\U$, then
   $\K^-(\r^-)=\K^+(\r+)$.
Thus the first $m$ kneading sequence entries
   remain constant throughout $M$ and agree with those of $\K(A)$. 

  This proves (\ref{item:changeofkneadingentries}). Items
   (\ref{item:componentkneading}) 
    through (\ref{item:kneadingbif}) are direct corollaries.
    \qed

 We are now ready to prove the main result of this section.

\begin{thm}[Determining 
  Components on a
  Combinatorial Arc]
 \label{thm:determiningbifurcationorder}
 Suppose that $A$ is a sector of a hyperbolic component $W$, and let 
  $\s\in\W(A)$. Let $j$ be the index of
  the first entry at which $\addu := \K(A)$ and $\addut := \K(\s)$ differ
  (or $j=\infty$ if no such entry exists).
 
 \begin{enumerate}
  \item 
   \label{item:mainbifstatement}
 Then there are no hyperbolic components of period less than $j$ on
  the combinatorial arc $[A,\s)$. If $j<\infty$, then there exists 
  a unique period $j$ component $V\in[A,\s)$. This component
  has forbidden kneading sequence $\KS(V)=\per{\u_1\dots\u_j}$;
  if furthermore
  $\ut_j\in\Z$, then $\s\in \W(\Sec(V,\ut_j))$.
 \item \label{item:forbiddenstatement}
 These statements remain true if $\s$ and $\addut$ are replaced
  by a hyperbolic component $W'$ and its forbidden kneading sequence $\KS(W')$.
 \end{enumerate}
\end{thm}
\proof 
 To prove
  (\ref{item:mainbifstatement}),
 let $m$ be the minimal period of a hyperbolic component on
  $[A,\s)$, and let $V\in[A,\s)$ be a component of period $m$.
  (If there are no hyperbolic components in
  $[A,\s)$, then $j=\infty$ by Lemma
  \ref{lem:kneadingandwakes}
  (\ref{item:changeofkneadingentries}), and there is nothing
  further to prove.) By Lemma \ref{lem:kneadingandwakes}
  (\ref{item:changeofkneadingentries}) and
  (\ref{item:componentkneading}), we have
  $m\leq j$ and $\K^*(V)=\per{\u_1\dots\u_m}$. In particular, $V$ is
  unique by Lemma \ref{lem:noequalkneadingsequences}. 

 If $\ut_m\notin\Z$, then
  $\ut_m\neq \u_m$, so
  $m=j$ and we are done.
  Otherwise,
  $\s\in\W(B)$ for some 
  sector $B$ of $V$. 
  By choice of $m$ and uniqueness of $V$, there are no components of period
  $\leq m$ on $[B,\s)$, and it follows from Lemma
  \ref{lem:kneadingandwakes} (\ref{item:changeofkneadingentries}) that
  $\ut_m = \u(B) \neq \u(V)=\u_m$; in particular, $m=j$.

Part (\ref{item:forbiddenstatement}) can be reduced
  to
  (\ref{item:mainbifstatement})
 by choosing $\s$ to be 
 an address
  just outside $\W(W)$.
  \qed

 \end{parameterarguments}

\begin{suppress}
 \subsection*{Wakes of Hyperbolic Components and Sectors}
 It is a well-known fact (see e.g.\ \cite{jackrays})
  that, for Mandel- and Multibrot sets,
  the wake of a hyperbolic component $W$ is exactly the subset of the plane
  in which the dynamic rays at the characteristic addresses of $W$ have
  a common repelling landing point.
  This result is also correct for exponential maps
  \cite{landing2new}, but the proof requires
  several deep results on parameter space, including those proved
  in this article and articles building on it, such as
  \cite{boundary}. The following observation gives a 
  combinatorial analog of this statement.

 \begin{lem}[Properties of $\W(W)$]
  Let $W$ be a hyperbolic component of period at least $2$ and external address
  $\s := \extaddr(W)$, with characteristic addresses $\s^-<\s^+$. 
  Let $\r\in\Sequ$.
   Then the following are equivalent:
  \begin{enumerate}
   \item $\r\in \W(W)$.
   \item $\displaystyle{\itin_{\r}(\s^+) = \itin_{\r}(\s^-)}$.
  \end{enumerate}
 \end{lem}
\proof (a) implies (b) by Observation \ref{obs:itinchange}. On the other hand,
   if $\r\notin \W(W)$, then there is an element of $\sigma^{-1}(\r)$ between
   $\sigma^{n-1}(\s^+)$ and $\sigma^{n-1}(\s^-)$, so the $n$-th itinerary
   entries of the two addresses differ. \qed

We wish to understand the bifurcation structure in the space of exponential
 maps; i.e., to understand which wakes a given hyperbolic component
 lies in. In order to carry this out, we subdivide the 
 (combinatorial) wake $\W(W)$ into pieces corresponding to the different
 sectors of $W$ as follows. 

\begin{defn}[Wakes of Sectors]
 If $A=\Sec(W,s_*)$ 
  is a sector of $W$, then the
  \emph{wake} of $A$ is defined to be the minimal interval
    \index{wake!of a sector}
 $\W(A)=(\r^1,\r^2)\subset\Sequ$ with
\nomenclature[WA]{$\W(A)$}{(wake of a sector $A$)}
  $\extaddr(W,s_*,\alpha)\in \W(A)$ for all rational $\alpha\in(0,1)$.
 The addresses $\r^1$ and $\r^2$ are called the
  \emph{sector boundaries} of $A$. Any sector boundary of a sector of $W$
  is also called a sector boundary of $W$.
    \index{sector boundary}
\end{defn}
 By Theorem \ref{thm:bifurcationstructure} (a) and (d), 
  the wakes of any two sectors
  of $W$ are disjoint, and the wakes of all sectors of $W$ together with
  their sector boundaries exhaust $\cl{\W(W)}\setminus\{\s\}$. 
  The sector boundaries
  of $W$ are exactly the characteristic addresses of $W$ together with
  all addresses of the form $\bifaddr(W,h)$ with $h\in\Z\setminus\{0\}$.

\begin{lem}[Preimages of $\W(W)$]  \label{lem:Im}
 Let $W$ be a hyperbolic component of period $n$ with
  kneading sequence $\K(W)=\u_1 \dots \u_{n-1} *$.
  Then for every $m\in\Z$ there exists a unique closed interval
  $I_m\subset \cl{\W(W)}$ with the following properties:
  \begin{enumerate}
   \item $\sigma^n$ maps $I_m$ bijectively onto $\cl{\W(W)}$;
     \label{item:bijectivity}
   \item every $\t\in I_m$ has $t_n=m$; 
     \label{item:nthentry}
   \item if $\t\in I_m$, then $\itin_{\r}(\t)$ begins with
    $\u_1\dots \u_{n-1}$ for every $\r\in\W(W)$.
      \label{item:itinerary}
  \end{enumerate}

  Furthermore, each $I_m$ contains a unique fixed point of $\sigma^n$ and satisfies
   \begin{equation}
    \sigma^j(I_m)\cap \W(W)=\emptyset \label{eqn:iteratesoutsidewake}
   \end{equation}
  for $j=1,\dots,n-1$.
  If $\t\in \Sequ$ is such that
   $\sigma^n(\t)\in \W(W)$ and 
   $\itin_{\r}(\t)$ begins with $\u_1\dots\u_{n-1}$, then
   $\t\in I_m$. 
 \end{lem}
 \proof Let $\s := \extaddr(W)$. Then
  $(m\s^-,m\s^+)$ contains no iterate of
  $\s$.
  By Observation \ref{obs:agreeingaddresses}, there
  exists an interval $I_m$ such that
  $\sigma^{n-1}$ maps $I_m$ bijectively to $[m\s^-,m\s^+]$ and such that
  $\itin_{\s}(\t)$ begins with $\u_1\dots \u_{n-1}$ for every
  $\t\in I_m$. So $I_m$ satisfies (\ref{item:bijectivity}) and
  (\ref{item:nthentry}). We should 
  also note that $\s^-,\s^+\notin \interior\left(\sigma^j(I_m)\right)$ 
  for every $j=1,\dots,n-1$. Thus 
  \begin{equation}
    \text{either $\sigma^j(I_m)\subset \cl{\W(W)}$ or
            $\sigma^j(I_m)\cap \W(W)=\emptyset$.} 
    \label{eqn:eitherinorout}
  \end{equation}

  The fact that $I_m\subset [\s^-,\s^+]$ follows from 
   a now familiar
   application of Observation \ref{obs:itinchange}.
   Since $\sigma^n(I_m)= [\s^-,\s^+]\supset I_m$, $I_m$ must contain
   a fixed point of $\sigma^n$, which is clearly unique
   since no interval is invariant under $\sigma^n$.

  Let $\t$ denote this
   fixed point, and let $\kappa\in W$.
   Then
   a hyperbolic expansion argument shows that
   $g_{\r}$ lands on the boundary of the characteristic 
   Fatou component $U_1$ of $\Ek$\footnote{%
   In fact, every ray with appropriate itinerary lands on $\partial U_1$.}
   (compare \cite[Proof of Theorem 6.2]{expattracting}, 
    \cite[Theorem 4.2.4]{thesis} or \cite[Theorem 6.3]{accessible}). 
   Since $\Ek^j(U_1)$ is separated from
   $U_1$ by the characteristic rays of $\kappa$ for $j=1,\dots,n-1$, we
   thus have
   $\sigma^j(t)\notin \W(W)$ for $j=1,\dots,n-1$. In view of
   (\ref{eqn:eitherinorout}), this proves (\ref{eqn:iteratesoutsidewake}).
   The latter in turn implies item (\ref{item:itinerary}) by
   Observation \ref{obs:itinchange}.

  Finally, if $\r$ and $\t$ are 
   as in the final claim,
   then $\t$ and
   $\wt{\t} := (\sigma^n_{I_{t_n}})^{-1}(\sigma^n(\t))$ are 
   $n$-th preimages of $\sigma^n(\t)$ which have the
   same itinerary under $\r$.  
   Thus $\t=\wt{\t}$ by Lemma
   \ref{lem:itineraries}. \qed

\begin{lem}[Sector Boundaries and Itineraries]
 \label{lem:sectorboundaries}
  Let $W$ be a hyperbolic component of period $n$ and 
   kneading sequence $\K(W)=\u_1\dots\u_{n-1} *$. 
   Set $\s := \extaddr(W)$ and let
   $\r\in\Sequ$. Then the following
   are equivalent.
  \begin{enumerate}
   \item $\r$ is a sector boundary of $W$; \label{item:sectorboundary}
   \item $\r$ is the unique fixed point of
     $\sigma^n$ in $I_{r_n}$;              \label{item:uniquefixedpoint}
   \item $\itin_{\r}(\s)=\K(W)$ and        \label{item:itincriterion}
     $\itin_{\s}(\r)=\per{\u_1\dots \u_{n-1} \m}$
   for some $\m\in\Z$.
 \end{enumerate}
  Furthermore, every sector boundary $\r$ has period $n$.
\end{lem}
\proof
 It follows easily from Lemma
  \ref{lem:itinerarybifurcation} that for every
  $\m\in\Z$ there is a sector boundary $\r$ with
  $\itin_{\s}(\r)=\per{\u_1\dots \u_{n-1} \m}$, and every sector boundary
  has an itinerary of this form. Recall that
  the characteristic addresses of $W$ have period $n$.
  Since there are no
  essential orbit portraits in the characteristic sector  
  of $\Ek$ (where $\kappa\in W$), it follows that any
  sector boundary which is not a characteristic address has
  the same period as its itinerary, and there is no other
  external address with the same itinerary. 
  In particular, 
  $\sigma^n(\r)=\r$ for every sector boundary of $W$.

 Therefore, (\ref{item:sectorboundary}) implies
  (\ref{item:uniquefixedpoint}) by the last statement in
  Lemma \ref{lem:Im}. 

 If $\t\in I_m$ with $\sigma^n(\t)=\t$, then 
  clearly its itinerary is as in
  (\ref{item:itincriterion}). Furthermore, 
  $\itin_{\r}(\s)=\K(W)$ for every
  $\r\in\cl{\W(W)}$ by Observation
  \ref{obs:itinchange}. So
  (\ref{item:uniquefixedpoint}) implies
  (\ref{item:itincriterion}). 

 Now let $\kappa\in W$ and
  suppose that $\r'$ is an 
  address which is not a sector boundary
  and has itinerary
  $\per{\u_1 \dots \u_{n-1} \m}$.
  Note that
  $\r'$ is necessarily periodic by Lemma
  \ref{lem:itineraries}. 
  There is a sector boundary $\r$ which has
  the same itinerary as 
  $\r'$, and by our initial remarks
  $\r$ must be a characteristic address of
  $W$, and $\r'\notin\cl{\W(W)}$.  
  By replacing
  $\r$ with $W$'s other
  characteristic address,
  if necessary, we may
  suppose that $\r$ and
  $\r'$ do not enclose
  $\s$. It then follows
  from Observation
  \ref{obs:agreeingaddresses}
  that they 
  must enclose
  a forward iterate of
  $\s$. Thus
  $\itin_{\r'}(\s)\neq \itin_{\r}(\s) = \K(\s)$ by Observation
  \ref{obs:itinchange}, as required.

 By (\ref{item:uniquefixedpoint}) and
  (\ref{eqn:iteratesoutsidewake}), every sector boundary has period $n$. 
 \qed

\begin{prop}[Periodic Rays and Intermediate Addresses]
   \label{prop:periodicintermediate}
  Let $\r\in\Sequ$ be periodic of period $n$. 
   Then there exists an intermediate
   external address $\s$ of length $n$ such that $\r$ is a sector boundary of 
   $W:= \Hyp{\s}$. If there exists $\wt{\r}$ such that
   $\langle\r,\wt{\r}\rangle$ is a characteristic ray pair, then
   $\r$ and $\wt{\r}$ are the characteristic addresses of $W$.
\end{prop}
\proof 
  Choose a parameter $\kappa$ for which the singular value lies on
   the dynamic ray $g_{\r}$.   
   (It is well-known --- and easy to see ---
    that such parameters exist for every periodic
    address; see 
    e.g.\ \cite{dghnew1}.\footnote{%
   It is possible to 
    formulate the following
    argument completely 
    combinatorially, 
    without appeal to the 
    existence of such 
    parameters.}
   In fact, this is true
    for every exponentially bounded external address, and these
    parameters form a corresponding
    \emph{parameter ray}; 
    see
   \cite{markus,markusdierk,markuslassedierk}.)
  Consider the piece of $g_{\r}$ which 
   connects the singular value with
   $+\infty$. The dynamic ray    
   $g_{\sigma^{n-1}(\r)}:(0,\infty)\to\C$ is a component of
   $\Ek^{-1}(g_{\r})$, 
  and tends
   to $+\infty$ for
   $t\to\infty$,
   and to $-\infty$ for $t\to 0$. The ray 
   $g_{\sigma^{n-1}}(\r)$ and its translates thus partition the 
   dynamic plane into ``strips'', similarly as we did in
   in Section \ref{sec:orbitportraits} for attracting parameters.

  Since 
  $g_{\sigma^{n-1}(\r)}
      =\Ek^{n-1}(g_{\r})$
  tends to $-\infty$ as
  $t\to 0$,
   the curve 
  $g_{\r}$
   address $\s$ of length 
  $n$ as $t\to 0$. 
 Since dynamic rays do not
   intersect, both ends of the ray 
   $g_{\r}$ tend to $\infty$ within the same strip of the abovementioned
   partition. In other words, the first itinerary entries of
   $\r$ and $\s$ (with respect to $\r$) coincide. We can apply the
   same argument to $g_{\sigma(\r)},\dots, g_{\sigma^{n-2}(\r)}$, and 
   conclude that  $\itin_{\r}(\r)$
   and $\itin_{\r}(\s)$ coincide in the first $n-1$ entries. To fix our ideas,
   let us suppose that $\r<\s$. 

 Note that for
  $j=1,\dots,n-2$,
  $\sigma^j(\s),
   \sigma^j(\r)\notin
  (\r,\s)$. Otherwise,
  we could choose a
  minimal such $j$, and
  Observation
  \ref{obs:agreeingaddresses}
  would imply that
  $\sigma^j$ maps 
  $[\r,\s]$ to
  $[\sigma^j(\r),
    \sigma^j(\s)]\subset
    [\r,\s]$, which is
  clearly impossible since
  no interval is invariant
  under the shift. 
 By Observation 
  \ref{obs:itinchange}, 
  it follows that
   \[ \itin_{\s}(\r)= 
     \K^-(\r) \text{ and }
      \K(\s) = 
      \itin_{\r}(\s). \]
   By Lemma 
  \ref{lem:sectorboundaries},
   $\r$ is a sector boundary of $\s$. 

 \smallskip

 If $\wt{\r}\in\Sequ$ is such that $\langle\r,\wt{\r}\rangle$ is a
  characteristic ray pair, then the cyclic order of $\r$, $\wt{\r}$ and
  $\s$ is preserved by $\sigma^{n-1}$, and thus 
  $\r<\s<\wt{\r}$. By Observation \ref{obs:itinchange}, 
  \[ \itin_{\s}(\r)=\itin_{\s}(\wt{\r}). \]
  Thus $\r$ and $\wt{\r}$ are the characteristic addresses of $\s$ by
  Lemma \ref{lem:uniqueaddress}. 
\qed

\subsection*{Kneading 
      Entries}

  If $\m\in\Z$, then the unique sector boundary  
   $\r\in I_{\m}$ has kneading sequence 
   \[ \per{\u_1\dots \u_{n-1} \bdyit{\m}{\m-1}}, \]  
  and we therefore denote this address $\r$ by
    $\Bdy\bigl(W,\bdyit{\m}{\m-1}\bigr)$. 
   \nomenclature[bdy]{$\Bdy\bigl(W,\bdyit{\m}{\m-1}\bigr)$}{(sector boundary)}

  If $A$ is a sector of $W$, then there exists a number
   $\u(A)$, called the \emph{kneading entry} of $A$, such that
\nomenclature[uA]{$\u(A)$}{(kneading entry of $A$)}
     \index{kneading entry}
   \[ \W(A) = \Bigl(\Bdy\bigl(W,\bdyit{\u(A)}{\u(A)-1}\bigr),
               \Bdy\bigl(W,\bdyit{\u(A)+1}{\u(A)}\bigr)\Bigr). \]
  (Compare Figure \ref{fig:sectorkneadingsequences}.)
  It makes sense to define the \emph{kneading sequence} of
   $A$ to be
   \[ \K(A) := \per{\u_1\dots \u_{n-1} \u(A)}. \]
    \index{kneading sequence!of a sector}

  Note that, if $\s\neq\infty$,
   no sector ever has the kneading sequence
   \[ \KS(W) := \itin_{\s}(\s^+)=\itin_{\s}(\s^-), \]
  which we therefore call the \emph{forbidden kneading sequence}
\nomenclature[KW]{$\KS(W)$}{(forbidden kneading sequence)}
  of $W$. Its $n$-th entry is called the \emph{forbidden kneading entry}
  and denoted by $\u(W)$. The period one component $\Hyp{\infty}$ does not
\nomenclature[uW]{$\u(W)$}{(forbidden kneading entry)}
  have a forbidden kneading sequence.

 For every $\m\neq \u(W)$, we denote by $\Sec(W,\m)$ the (unique) sector
  of $W$ whose kneading entry is $\m$.
  
 The preceding
   definitions are further justified by the following proposition.

 \begin{prop}[Kneading Sequences and Sectors]
  \label{prop:kneadingandsectors}
  Let $W$ be a hyperbolic component with kneading sequence
   $\K(W)=\u_1 \dots \u_{n-1} *$. 
  \begin{enumerate}
   \item Suppose that $\r\in \W(W)$ such that
    $\K(\r)=\u_1 \dots \u_{n-1} \m \dots$, where $\m\in\Z$. Then
    $\m\neq \u(W)$ and $\r\in \W(\Sec(W,\m))$. 
     \label{item:kneadingentry}
   \item Suppose that $h\in\Q\setminus\Z$. Then
    $\KS(\child)=\K(\Sec(W,h))$.
       \label{item:kneadingchild}
   \item Suppose that $h\in\R\setminus\Q$. Then
     $\K(\bifaddr(W,h))=\K(\Sec(W,h))$. 
       \label{item:kneadingbif}
  \end{enumerate}
 \end{prop}   
 \proof
  To prove  (\ref{item:kneadingentry}). 
  let
  $\m\in\Z\setminus\{\u(W)\}$, let $\t^-<\t^+$ be 
  the
  sector boundaries of $A:= \Sec(W,\m)$ and let 
  $\r\in \W(A)$ be an 
  address whose kneading 
  sequence begins with
  $\u_1,\dots,\u_{n-1}$. 
  We must show that the next
  entry of
  $\K(\r)$ is $\m$.

  Since $\sigma^{n-1}$ 
   preserves the cyclic 
   order of
   $\t^-,\r$ and $\t^+$, 
  it follows that 
  \[ \sigma^{n-1}(\t^-) <
     \sigma^{n-1}(\r) < 
     \sigma^{n-1}(\t^+). \]

  If $\sigma^n(\r)\notin \W(W)$, then
   \begin{align*}
     \sigma^{n-1}(\t^-) \leq \m\s^+ &<
      \sigma^{n-1}(\r) < (\m+1) \s^- 
      \leq \sigma^{n-1}(\t^+), \quad \text{and thus}\\
  \m\r \leq
     \m \s^+ &< \sigma^{n-1}(\r) < (\m+1) \s^-
     (\m+1)\r. \end{align*}

  Now suppose, on the other hand, that
   $\sigma^n(\r)\in \W(W)$. Note that, since
   $\sigma^n(\r)$ lies between $\t^-$ and
   $\t^+$, the entry
   $r_n$ is either $\m$ or $\m+1$; to
   fix our ideas, let us suppose that
   $r_n = \m$. Then 
   $\r$ belongs to the interval
   $I_{\m}$ of Lemma \ref{lem:Im},
   and $\sigma^n$ maps the interval
   $[\t^-,\r]\subset I_{\m}$ bijectively
   to $[\t^-,\sigma^n(\r)]$. Thus
   we must have $\sigma^n(\r)>\r$, so that
   the $n$-th itinerary entry of $\r$ is
   $r_n = \m$, as required.

\smallskip

 To prove (\ref{item:kneadingchild}), let
  $h\in\Q\setminus\Z$ and  $V:=\child$. 
  Recall that, if $\kappa\in\IR$, then 
  $\bifaddr(W,h)$ is the address of the broken attracting ray $\gamma$
  of $\Ek$
  which contains the singular value. Thus $\K(V)$ begins with
  $\u_1,\dots,\u_{n-1}$. By Proposition \ref{prop:childcomponents}
  (\ref{item:childcharacteristic}), the forbidden kneading sequence of 
  $V$ is periodic of period $n$; i.e., 
  $\KS(V)=\per{\u_1\dots \u_{n-1} \m}$ for some $\m\in\Z$.
  The claim now follows by (\ref{item:kneadingentry}). 

\smallskip

 Finally, (\ref{item:kneadingbif}) is easily deduced from
  (\ref{item:kneadingchild}) by stability of kneading sequence entries. \qed

\begin{cor}[Nested 
  Components Have Different
    Kneading Sequences]
  \label{cor:noequalkneadingsequences}
 Let $W$ be a hyperbolic component. If $V\neq W$ 
  is a hyperbolic component with
  $\extaddr(V)\in\W(W)$, then $\K(V)\neq \K(W)$. 
\end{cor}
\proof Suppose by contradiction that $V$ was a hyperbolic component with
  $\extaddr(V)\in\W(W)$ and $\K(V) = \K(W) = \u_1\dots \u_{n-1} *$.
  Set $\m := \u(W)$.
  If we choose $\s$ slightly below 
  $\Bdy\bigl(V,\bdyit{\m+1}{\m}\bigr)$, then
  $\K(\s)$ begins with $\u_1\dots\u_{n-1}\m$, contradicting
  Lemma \ref{lem:kneadingandwakes} (\ref{item:kneadingentry}).\qed

\subsection*{Bifurcation Order and Internal Addresses}
 \begin{defn}[Bifurcation Order and Combinatorial Arcs]
  If $A,B$ are hyperbolic components or sectors, we write 
  $A\prec B$ if $\W(A)\supset \W(B)$. The \emph{combinatorial arc}
  $[A,B]$ is the set of all $C$ such that
  $A\prec C \prec B$. 
\nomenclature{$\prec$}{(combinatorial order)}
\nomenclature{$[A,B]$}{(combinatorial arc)}

  Similarly, if $\s\in \cl{\W(A)}$, then the combinatorial arc
  $[A,\s)$ is the set of all $C$ with $A\prec C$ and
  $\s\in \cl{\W(C)}$.
   \index{combinatorial arc}
\nomenclature{$[A,\s)$}{(combinatorial arc)}
 \end{defn}
 \begin{remark}
  We will often also consider
  open or half-open combinatorial arcs $(A,B)$, $[A,B)$ or $(A,B]$, in which
  one or both of the endpoints are excluded.
 \end{remark}

\begin{figure}
\subfigure[Sector Labels]{\input{sector_labels.epstex}%
             \label{fig:seclabels}}\hfill
\subfigure[Kneading Entries]{\input{sector_kneading.epstex}%
             \label{fig:sectorkneadingsequences}}\hfill
\subfigure[Sector Numbers]{\input{sector_numbers.epstex}%
             \label{fig:sectornumbers}}
 \caption{Three ways to label sectors, illustrated for the component
  at address $\s=0110\frac{1}{2}\infty$, with
  $\s^-=\per{011002}$ and $\s^+=\per{020101}$. The central internal ray of
  $\Hyp{\s}$ is
  emphasized.}
\end{figure}

\end{suppress}

 \begin{parameterarguments}

\subsection*{Internal addresses}
 Repeated applications of 
  the preceding theorem enable us to determine the periods and combinatorial
  order of all hyperbolic components on the combinatorial arc between
  $\Hyp{\infty}$ and $\s$ solely from $\K(\s)$. 
  \emph{Internal addresses}, introduced for Mandel- and Multibrot sets
  in \cite{intaddr}, organize this information in a
  more convenient way.
 \end{parameterarguments}

\begin{suppress}%
 We are now ready to define internal addresses, which are equivalent to
  kneading sequences, as we show below, and give a (human-readable)
  description of where a given hyperbolic component --- or, more generally,
  external address --- is situated in parameter space.
\end{suppress}%
 Since neither 
  sector labels nor kneading entries can be easily identified in a 
  picture of parameter space, let us introduce a third labelling method
  for the sectors of a hyperbolic component.
  The \emph{sector number} of a given sector is the nonzero
  integer obtained by counting sectors in counterclockwise orientation,
  starting at the central internal ray.
\index{sector number}%
  In other words, the sector number of a sector
   $A$ is the integer
   $\u(A) - \u(W)\in\Z\setminus\{0\}$ 
   and can be found in parameter space by counting
   from the root of $W$ (compare Figure \ref{fig:sectornumbers}). 

\begin{suppress}%
 \begin{defn}[Internal Addresses] \label{defn:internaladdress}
  Let $\s\in\Sequb$. Consider the sequence $W_1=\Hyp{\infty} \prec
    W_2 \prec \dots$ of all hyperbolic components
   $W\in [\Hyp{\infty},\s)$ which have the property
   that all components on $(W,\s)$ have higher period than $W$.
   Denote by $n_j$ the period of $W_j$, and by
    $m_j$ the sector number of the sector of $W_j$ containing
    $\s$. (We adopt the convention that
    $m_j=\infty$ if $\s=\extaddr(W_j)$ and
    $m_j=k+\onehalf$ if 
    $\s=\Bdy\bigl(W_j,\bdyit{\u(W_j)+k+1}{\u(W_j)+k}\bigr)$.)
     \index{internal address}

   The \emph{internal address of $\s$} is the (finite or infinite) sequence
    \[ (n_1,m_1) \intaddnext (n_2,m_2) \intaddnext (n_3,m_3) \intaddnext \dots \]
   The internal address of an intermediate external address $\s$ is also
    called the internal address of the associated hyperbolic component
    $\Hyp{\s}$.
 \end{defn}
\end{suppress}%
\begin{parameterarguments}%
 \begin{defn}[Internal Addresses] \label{defn:internaladdress}
  Let $\s\in\Sequb$. Consider the sequence $W_1=\Hyp{\infty} \prec
    W_2 \prec \dots$ of all hyperbolic components
   $W\in [\Hyp{\infty},\s)$ which have the property
   that all components on $(W,\s)$ have higher period than $W$.
   Denote by $n_j$ the period of $W_j$, and by
    $m_j$ the sector number of the sector of $W_j$ containing
    $\s$. (We adopt the convention that
    $m_j=\infty$ if $\s=\extaddr(W_j)$ and
    $m_j=k+\onehalf$ if $\s$ is the 
    sector boundary with kneading sequence
    $\K(\s)=\per{\u_1\dots\u_{n_j-1}\bdyit{\u(W_j)+k+1}{\u(W_j)+k}}$.)
     \index{internal address}%

   The \emph{internal address of $\s$} is defined as
    the (finite or infinite) sequence
    \[ (n_1,m_1) \intaddnext (n_2,m_2) \intaddnext (n_3,m_3) \intaddnext \dots \]
   The internal address of an intermediate external address $\s$ is also
    called the internal address of the associated hyperbolic component
    $\Hyp{\s}$.
 \end{defn}
\end{parameterarguments}%
\begin{remark}[Remark 1]
 By Theorem \ref{thm:determiningbifurcationorder},
  we could alternatively define $W_{j+1}$ as the 
  unique period of lowest period on the combinatorial
  arc $(W_j,\s)$. In particular, the components are indeed
  ordered as stated and 
  the sequence $n_j$ is strictly increasing. (These facts
  also follow
  easily directly from the definition.) 
\end{remark}
\begin{remark}[Remark 2]
 Internal addresses do not label hyperbolic components
  uniquely, reflecting certain symmetries of parameter space:
  two child components of a
  given sector with the same denominator but differing numerators of
  the bifurcation angle have the same internal address. This is the
  only ambiguity and thus uniqueness can be achieved 
  by specifying bifurcation angles in the internal address, 
  see Theorem \ref{thm:angleduniqueness}. 
\end{remark}

\begin{suppress}

The key to computing internal addresses --- and bifurcation order 
  in general --- lies in the following statement.

\begin{thm}[Determining 
        Bifurcation Order]
 \label{thm:determiningbifurcationorder}
 Suppose that $A$ is a sector of a hyperbolic component $W$, and let 
  $\s\in\W(A)$. Let $j$ be the index of
  the first entry at which $\addu := \K(A)$ and $\addut := \K(\s)$ differ
  (or $j=\infty$ if no such entry exists).
 \begin{enumerate}
  \item If $V$ is a hyperbolic component on
    the combinatorial arc $[A,\s)$, then $\KS(V)$ begins with
    $\u_1\dots \u_{j-1}$. In particular, the period of $V$ is at
    least $j$. \label{item:periodatleastj}
  \item If $j<\infty$,
   then there exists a unique hyperbolic component $V$ of period
   $j$ on $[A,\s)$. If $\ut_j\in\Z$, then $\s\in\W(\Sec(V,\ut_j))$.
   \label{item:existenceuniqueness}
 \end{enumerate} 

 These statements remain true if $\s$ and $\addut$ are replaced
  by a hyperbolic component $W'$ and its forbidden kneading sequence $\KS(W')$.
\end{thm}
\proof
 Let $n$ be the period of $W$ and let
  $\r^-<\r^+$ be the characteristic addresses of $W$. Recall that
  \[\itin_{\s}(\r^{\pm}) = 
    \per{\u_1\dots \u_{n-1}
            (\u(W))}.\]

 \medskip
 
 \noindent\textsc{Claim 1. }
   For any
   $V\in [A,\s)$, $\KS(V)$ begins with
   $\u_1\dots \u_{j'-1}$, where $j'=\min(j,n)$. 
   In particular, the period of $V$ is at least $j'$.

\smallskip

  \emph{Proof of Claim 1. }
  Denote the characteristic ray pair of $V$ by 
  $\langle\wt{\r}^-,\wt{\r}^+\rangle$;
  recall that $\itin_{\s}(\wt{\r}^{\pm})=\KS(V)$. 
  Suppose by contradiction that there exists $m<j'$ such that
  the $m$-th entries of
  $\addu$ and $\KS(V)$ do not agree; we may suppose that $m$ is chosen
  minimal with this property. Then $\sigma^{m-1}$ preserves
  the circular order of the five
  addresses
  \[ \r^-<\wt{\r}^-<\s<\wt{\r}^+<\r^+. \]
  Thus
  \begin{align*}
   &\sigma^{m-1}(\r^-)< \sigma^{m-1}(\wt{\r}^-) <
     \sigma^{m-1}(\s)
  \hspace{.4cm}\text{or} \\
   & \sigma^{m-1}(\s) <  \sigma^{m-1}(\wt{\r}^+) <
     \sigma^{m-1}(\r^+). 
  \end{align*}
  In either case, the outer two addresses lie in the same interval of
  $\Sequ\setminus\sigma^{-1}(\s)$. So all three addresses have the
  same first itinerary entry, which contradicts the definition of
  $m$. 

 Thus $\KS(V)$ begins with $\u_1\dots \u_{j'-1}$, as required.
  It follows from Proposition \ref{prop:kneadingandsectors} 
  (\ref{item:kneadingentry}),
  applied to $V$ and $\s$, that $V$ cannot be of period less
  than $j'$. This proves Claim 1.

\medskip

\noindent
\textsc{Claim 2. } (\ref{item:periodatleastj}) and
   (\ref{item:existenceuniqueness}) hold when $j<n$. 

\smallskip

\emph{Proof of Claim 2. }
  Suppose that $j<n$. Then (\ref{item:periodatleastj}) follows directly from
  Claim 1. 
  Note also that Claim 1 implies that any hyperbolic component
  $V\in [A,\s)$ of period $j$ has kneading sequence
  $\K(V)=\u_1\dots \u_{j-1}*$. In particular, there is at most one such 
  component by Corollary \ref{cor:noequalkneadingsequences}. 

 Let us show that, for any such component $V$, 
  $\KS(V)=\per{\u_1\dots\u_j}$. That is, we need to show that
  $\u(V)=\u_j$. 
  Let
  $\wt{\r}:=\extaddr(V)$. Once again, the cyclic
  order
  \[ \r^-<\wt{\r}^-<\wt{\r}<\wt{\r}^+<\r^+ \]
  is preserved under $\sigma^{j-1}$, or in other words
  \[ \sigma^{j-1}(\wt{\r}^+) < \sigma^{j-1}(\r^+) <
     \sigma^{j-1}(\r^-) < \sigma^{j-1}(\wt{\r}^-), \]
  and it follows that these four addresses are all contained
  in the interval $(\u_n \wt{\r} , (\u_n+1)\wt{\r})$, as required. 

 To complete the proof of (\ref{item:existenceuniqueness}), we thus only need
  to show the existence of a component
  $V\in[A,\s)$ of period $\leq j$. To find this component, let us
  define
  two sequences $\s_k^-$ and $\s_k^+$ inductively by setting
  $\s_0^-:=\s_0^+:=\s$ and letting
  $\s_{k+1}^+$ (resp.\ $\s_{k+1}^-$) be the unique pullback of
  $\u_{j} \s_k^+$ (resp.\ $(\u_{j}+1)s_k^-$) under
  $\sigma^{j-1}$ whose itinerary begins with
  $\u_1\dots \u_{j-1}$. 

 It easily follows that
  \[ \r^- < \dots < \s_2^- < \s_1^- < \s < \s_1^+ < \s_2^+ < \dots <
     \r^+. \] 
 Thus $\s_k^-$ and $\s_k^+$ converge to
  addresses $\wt{\r}^-$ and $\wt{\r}^+$ with   
   $\itin_{\s}(\wt{\r}^-)=\itin_{\s}(\wt{\r}^+)=
    \per{\u_1\dots\u_j}$. Thus, these addresses belong to an
   essential orbit portrait, whose characteristic sector contains
   $\s$ and is contained in $\W(A)$. By 
   Proposition \ref{prop:periodicintermediate}, 
   this characteristic sector is the wake
   of some hyperbolic component of period $\leq j$, as required.
   This completes the proof of Claim 2; i.e., the case where $j<n$. 

\smallskip

  Since $j\neq n$ by Proposition \ref{prop:kneadingandsectors} 
  (\ref{item:kneadingentry}), it remains to consider
  the case $j>n$.  
  If 
  $\s=\bifaddr(W,h)$ for some $h\in\R\setminus\Z$, or if 
  $\s$ is a sector boundary of some child component of $W$, then 
  the claim follows easily from Proposition \ref{prop:kneadingandsectors}
  (\ref{item:kneadingchild}) and
  (\ref{item:kneadingbif}).

  Otherwise, $\s$ is contained in a sector $A'$ of some child component 
  $W':=\child$ by Corollary
  \ref{cor:subwakesfillwake}. Let $qn$ be the period of $W'$. By Proposition
  \ref{prop:kneadingandsectors} (\ref{item:kneadingchild}), $\KS(W')=\addu$.
  By Proposition \ref{prop:kneadingandsectors}
  (\ref{item:kneadingentry}), $\addut$ cannot agree with
  $\KS(W')$ for more than $qn-1$ entries; i.e.\ 
  $j\leq qn$. If $j< qn$, we can apply Claim 2, 
  replacing $W$ and $A$ by $W'$ and $A'$. 

  So consider the case where $j=qn$. It follows
   from Claim 1, applied to 
   $A'$ and $\s$, that there are no hyperbolic components of period
   $< qn$ on the combinatorial arc $[A',\s)$. This proves
   (\ref{item:periodatleastj}). Also note that
   $V:= W'$ is a period $j$ component with the required kneading sequence.
   Finally,
   there is no component of period $qn$ on the arc
   $[A',\s)$, as the kneading sequence of such a component
   would agree with $\K(W')$ by Claim 1. This is impossible by
   Corollary \ref{cor:noequalkneadingsequences}.

\smallskip

 To prove the last claim of the Proposition, 
  note that the statement with $\s$ and $\addu$ replaced by
  a hyperbolic component $W'$ and its forbidden kneading sequence can be
  reduced to the original formulation by taking $\s$ to be an address
  which is chosen just outside $\W(W')$. 
 \qed

\end{suppress}

\begin{cor}[Computing Internal Addresses]
  \label{cor:internaladdress}
  Two external addresses have the same internal address if and only
   if they have the same kneading sequence.  
 
  Furthermore,
   the internal address 
   $(1,m_1)\intaddnext (n_2,m_2)\intaddnext (n_3,m_3)\intaddnext\dots$
   of any $\s\in\Sequ$   
   can be determined inductively 
   from $\addu:=\K(\s)$ by the following procedure:

  Set $m_1 := \u_1$. 
   To compute $(n_{i+1},m_{i+1})$ 
   from $n_i$, continue the first $n_i$ entries of
   $\addu$ periodically to a periodic sequence $\addu^i$. 
   Then $n_{i+1}$ is the position of the first difference between
   $\u$ and $\addu^i$. Furthermore,
   \[ m_{i+1} = \begin{cases}
               \u_{n_{i+1}} - \u^i_{n_{i+1}} & \text{if }  \u_{n_{i+1}}\in\Z\\
               \infty & \text{if } \u_{n_{i+1}} = * \\
               k-\u^i_{n_{i+1}}+\onehalf 
                       &\text{if } \u_{n_{i+1}} = \bdyit{k+1}{k}.
            \end{cases}
   \]  
   (If $\K(\s)$ is periodic of period $n_{i+1}$, or if 
   $\s$ is intermediate of length $n_{i+1}$, then the
   algorithm terminates, and the internal address is finite.)
 \end{cor}
 \begin{remark}[Remark 1]
  $\addu^i$ is the kneading sequence of the sector of $W_i$ containing
  $\s$. The forbidden kneading sequence $\KS(W_i)$ can be obtained by
  repeating the first $n_i$ entries of $\addu^{i-1}$ periodically. 
  In particular, if $\s$ is an intermediate external address, then the
  forbidden kneading sequence of $\s$ consists of the first $n_{k-1}$ 
  entries of $\K(\s)$ repeated periodically, where $k$ is the length of
  the internal address of $\s$. 
 \end{remark}
 \begin{remark}[Remark 3]
  As an example, let us consider $\s=030\frac{1}{2}\infty$. Then we have
   $\K(\s)={\tt 0200*}$. Applying the above procedure, we obtain
   that $\addu^1=\per{0}$, $\addu^2=\per{02}$ and 
   $\addu^3=\per{0200}$, resulting in the internal address
    \[ (1,0) \intaddnext (2,2) \intaddnext (4,-2) \intaddnext (5,\infty). \]
 \end{remark}
 \begin{remark}[Remark 3]
  There is an obvious converse algorithm: given the internal address
   of $\s$, we can determine the kneading sequence $\K(\s)$ by inductively
   defining $\addu^i$. 
 \end{remark}
 \proof The correctness of the algorithm is an immediate corollary of
   Theorem \ref{thm:determiningbifurcationorder}. 
   In particular, the internal address of $\s$ depends only on
   $\K(\s)$. Conversely, applying this procedure to two different
   kneading sequences will produce different internal addresses. \qed

 We note the following consequences of 
  Theorem \ref{thm:determiningbifurcationorder} for further reference.
 
 \begin{cor}[Combinatorics
      of Nested Wakes] \label{cor:nestedwakes}
 Let $W$ and $V$ be hyperbolic components with $W\prec V$. Then all entries
 of $\K(W)$ also occur in 
 $\K(V)$.
\end{cor}
\proof Let $p$ be the 
  period of $W$, and let $\u_1\dots \u_{p-1}*$ be the kneading
  sequence of $W$. 
  The proof proceeds by induction on the number $n$ of hyperbolic
 components on the 
 combinatorial arc $(W,V)$
  which have period 
  less than $p$. 
 If $n=0$, then it follows from Theorem
 \ref{thm:determiningbifurcationorder} 
 that $\KS(V)$ begins with $\u_1\dots \u_{p-1}$, and we are done.
 (Recall from  Lemma
  \ref{deflem:kneadingentries} that all entries of $\KS(V)$ occur in
  $\K(V)$.)

If $n>0$,
 then let $V'\in (W,V)$ be a component of period $\leq p$. We can now
 apply the induction 
 hypothesis first to $W$ 
 and $V'$, and then to $V'$
 and $V$. \qed

\begin{cor}[Components on the Combinatorial Arc] \label{cor:unboundeditinerary}
 Let $\s\in\Sequ$ and $\addu := \K(\s)$. Suppose that $n\geq 1$ such
 that $\u_n\in\Z$ and
 $\u_j\in \Z\setminus\{\u_n\}$ for all $j<n$. Then there exists a hyperbolic
 component $W$ with $\KS(W)=\per{\u_1\dots \u_n}$ and $\s\in\W(W)$.
\end{cor}
\proof By the internal address algorithm, $n$ appears in the internal
 address of $\s$; let $V$ be the associated period $n$ component. The
 child component of $V$ containing $\s$ has the required property. \qed

 \subsection*{Infinitely 
  many essential periodic 
  orbits}

 To conclude this section, we will give a simple 
  necessary and sufficient criterion
  for an attracting exponential map to have infinitely many
  essential periodic orbits. 
  (A non-necessary sufficient condition
   under which this occurs was the main result of
   \cite{tying}.)

  \begin{prop}[Infinitely 
     Many Essential Orbits]
   Let $W$ be a hyperbolic 
   component and
   $\kappa\in W$.
  Then the
   characteristic ray pairs of essential periodic orbits of 
   $\Ek$ are exactly the characteristic ray pairs of hyperbolic components
   $V$ with $V\prec W$.    
   
  In particular, the number of
   essential periodic 
   orbits of
   $\Ek$ is finite
  if and only if 
  the internal address
    of $W$ is of the form
    \begin{equation}
  (1,m_1) \intaddnext (n_2,m_2) \intaddnext (n_3,m_3) \intaddnext \dots 
     (n_k,\infty), \label{eqn:finitelymanyorbits} \end{equation}
    with $n_j|n_{j+1}$ for all $j<k$. In this case, the number of
    essential periodic orbits is exactly $k-1$. 
  \end{prop}
 \begin{remark}
   Using the internal address algorithm, it is simple to convert 
    (\ref{eqn:finitelymanyorbits})
    to a (somewhat more complicated) statement about the
    kneading sequence of $W$.
 \end{remark}
  \proof The first statement follows immediately from Lemma
   \ref{lem:charraypairs} and Proposition \ref{prop:periodicintermediate}.
  In particular, $\Ek$ has only
   finitely many essential
   periodic orbits if and only if $W$ is contained in only
   finitely many wakes.
   This is the case if and only if
   $W$ can be reached by finitely many bifurcations from the period one
   component $\Hyp{\infty}$, which is exactly what the statement about
   internal addresses means. \qed

\appendix

\section{Further Topics} \label{app:furthertopics}
 In this appendix, we will treat some further developments which
  are naturally related to the discussion in this article. 
  In the previous sections, we have been careful to give self-contained
  combinatorial proofs for all presented theorems.
  With the exception of Theorem \ref{thm:ADDR1}, which is again
  given an independent proof, the subsequent results
   will not be required in the proofs of the further results referred
   to in the introduction. Thus, we will often simply sketch how to
   obtain them from well-known
   facts in the polynomial setting.

 \subsection*{Addresses of Connected Sets}

 An important application of the results in this article is to obtain
  control over the combinatorial position of curves and, more generally,
  connected sets within exponential parameter space. 

 To make this precise, suppose that $A\subset\C$ is connected and contains
  at most one attracting or indifferent parameter. Let $\s\in\Sequ$ and
  suppose that there exist two hyperbolic components
  $W_1,W_2$ with $\extaddr(W_1)<\s<\extaddr(W_2)$ and the following property:
  there is $R>0$ such that every component
   $U$ of
   \[ \{\re z > R\}\setminus\Bigl(\IRWH{W_1}{0}\bigl((-\infty,-1]\bigr) \cup
                              \IRWH{W_2}{0}\bigl((-\infty,-1]\bigr)\Bigr) \]
  which is unbounded but has bounded imaginary parts satisfies
  $U\cap A = \emptyset$. In this case, we say that
  \emph{$A$ is separated from $\s$}. We define
    \[ \ADDR(A) := \{\s\in\Sequ: \text{$A$ is not separated from $\s$}\}. \]
\nomenclature[AddrA]{$\ADDR(A)$}{(addresses associated to a connected set $A$)}
  Note that $\ADDR(A)$ is a closed subset of $\Sequ$, and that
   $\ADDR(A)$ is empty if and only if $A$ is bounded. 

  \begin{remark}
   If $\gamma:[0,\infty)\to\C$ is a curve to infinity which contains
    at most one attracting or indifferent parameter, then 
    $\ADDR(\gamma)$ consists of a single external address;
    compare also \cite[Section 2]{boundary}. 
    In particular, if $G_{\s}$ is a \emph{parameter ray tail} as defined in
    \cite{markus}, then 
    $\ADDR(\gamma)=\{\s\}$.
  \end{remark}

 \begin{thm}[Addresses of Connected Sets]
   \label{thm:ADDR1}
  Let $A\subset\C$ be connected and contain at most one
   attracting or indifferent parameter. Then either
   \begin{enumerate}
    \item All addresses in $\ADDR(A)$ have the same kneading sequence, or
    \item $\ADDR(A)$ consists of the characteristic addresses of some
     hyperbolic component.
   \end{enumerate}
 \end{thm}
 \proof We claim that, for any hyperbolic component $W$, 
   either $\ADDR(A)\subset \W(W)$ or $\ADDR(A)\cap \W(W)=\emptyset$.
   In fact, we prove the following stronger fact.

 \smallskip\noindent
 \emph{Claim.} If $\ADDR(A)\cap\W(W)\neq\emptyset$, then there exists
  $h_0\in \R\cup\{\infty\}$ such that 
   $\ADDR(A)=\{\bifaddr(W,h_0)\}$ if $h\notin \Q\setminus\Z$ and
   $\ADDR(A)\subset \{\W(\childH{h_0})\}$ otherwise. 
    (We adopt the convention
   that $\bifaddr(W,\infty)=\extaddr(W)$.)
 
 \smallskip

 \noindent
 \emph{Proof.}
  Let us suppose that 
   $\extaddr(W)\neq \infty$ 
   (the proof in the period one case is completely analogous),
   and suppose that $\s \in \ADDR(A) \setminus \extaddr(W)$ 
   (if no such address $\s$ exists, then
   there is nothing to prove). To fix our ideas, let us assume that
   $\s<\extaddr(W)$. 

  By Corollary \ref{cor:subwakesfillwake}, there exists
   $h_0\in (\R\setminus\{0\})$ such that 
   $\s=\bifaddr(W,h_0)$ if $h_0\notin \Q\setminus\Z$ and
   $\s\in \W(\childH{h_0})$ otherwise. 

   By Lemma \ref{lem:internallanding}, we can choose rational
   $h^-$ and $h^+$ with 
    $h^- < h_0 < h_1^+$ 
   arbitrarily close to $h_0$
   such that $\Psi_W(ih^{\pm})\in\C$.
   If $A$ contains a (necessarily unique) indifferent parameter 
   $\kappa_0$, then we may suppose that these values are chosen
   such that 
   $\kappa_0\notin \Psi_W\bigl(i[h^-,h^+]\bigr)\setminus 
                   \Psi_W(ih_0)$. 

  Consider the Jordan arc
    \[ \gamma := \IRWH{\childH{h^-}}{0} \cup
        \{\Psi_W(ih^-)\}\cup \wt{\gamma} \cup 
        \{\Psi_W(ih^+)\}\cup \IRWH{\childH{h^+}}{0}, \]
  where $\wt{\gamma}\subset W$ is some curve connecting
   $\Psi_W(ih^-)$ and $\Psi_W(ih^+)$. Let
   $U$ denote the component of $\C\setminus\gamma$ which does
  not contain a left half plane. Since 
    \[ \r^+ := \bifaddr(W,h^+) > \s > 
       \r^- := \bifaddr(W,h), \]
  it follows
  by the definition of $\ADDR(A)$ that $U\cap A\neq \emptyset$. Also,
   $A\cap \gamma=\emptyset$ and $A$ is connected, so
   $A\subset U$. Therefore $A$ is separated from every address in 
   $\Sequ\setminus [\r^-,\r^+]$, and so
    $\ADDR(A)\subset [\r^-,\r^+]$. Letting $h_1^+$ and $h_1^-$ 
    tend to
    $h_0$, we have 
   \[ \ADDR(A)
      \subset \bigl[\lim_{h\nearrow h_0} \bifaddr(W,h)\, , \, 
                    \lim_{h\searrow h_0} \bifaddr(W,h)\bigr], \]
    as required. 

\smallskip

 Let us distinguish two cases.

 \noindent
  \emph{Case 1: Some address in $\ADDR(A)$ has an infinite internal address. }
  It then follows from the claim that all addresses have the same
   internal address, and by Corollary \ref{cor:internaladdress}, they
   also all share
   the same kneading
   sequence.

 \smallskip
 \noindent
  \emph{Case 2: All addresses in $\ADDR(A)$ have a bounded internal address. }
  It follows from the claim and the definition of internal addresses
  that there exists some hyperbolic component
  $W$ such that $\ADDR(A)\subset \cl{\W(W)}$ while
  $\ADDR(A)$ is not contained in the wake of any child component of $W$. 
  Suppose that $\ADDR(A)$ contains more than one external address.

  If $\ADDR(A)\not\subset \W(W)$, then 
  $\ADDR(A)$ must consist of the two
  characteristic addresses of $W$.
  If $\ADDR(A)\subset \W(W)$, then by the claim there exists
  $h_0\in \Q\setminus\Z$ with
    \[ \ADDR(A)\subset \cl{\W(\childH{h_0})}\setminus \W(\childH{h_0}), \]
  and again $\ADDR(A)$ consists of the two characteristic addresses of
  $\childH{h_0}$. \qed

 We record the following special case for use in \cite{markuslassedierk}
  (see there or in \cite{markus} for definitions). 
  \begin{cor}[Parameter Rays Accumulating at a Common Point]
   Suppose that $G_{\s^1}$ and $G_{\s^2}$ are parameter rays which
    have a common accumulation point $\kappa_0$. Then 
    $|s^1_j - s^2_j| \leq 1$
    for all $j\geq 1$.
  \end{cor}
  \proof Let $A := G_{\s^1}\cup G_{\s^2} \cup \{\kappa_0\}$. 
   Then $\ADDR(A)$ contains $\s^1$ and $\s^2$, and by the previous
   theorem, either $\K(\s^1)=\K(\s^2)$ or $\s^1$ and $\s^2$ are the
   characteristic addresses of some hyperbolic component. In either case,
   the claim follows. \qed

 In fact, we can sharpen Theorem
  \ref{thm:ADDR1} to the following statement.

 \begin{thm}[Addresses of Connected Sets II] \label{thm:ADDR2}
  Let $A\subset\C$ be connected and contain at most one
   attracting or indifferent parameter. Then exactly one of the following
   holds.
   \begin{enumerate}
    \item $\ADDR(A)$ is empty or consists of a single external address.
    \item $\ADDR(A)$ consists of two bounded external
      addresses, both of which
      have the same kneading sequence.
    \item $\ADDR(A)$ consists of the characteristic addresses of some
     hyperbolic component.
    \item $\ADDR(A)$ consists of at least three but finitely many
      preperiodic addresses. Moreover, 
      there exists a postsingularly finite
      parameter $\kappa_0\in\C$ such that, for every $\s\in\ADDR(A)$,
      the parameter ray $G_{\s}$ lands at $\kappa_0$. 
      \label{item:misiurewicz}
   \end{enumerate}
 \end{thm}
 \begin{remark}
   The \emph{Squeezing Lemma}, proved in 
    \cite{boundary}, also shows that
    $\ADDR(A)$ cannot contain intermediate or 
    exponentially unbounded addresses.
 \end{remark}
 \sketch Recall from the proof of Theorem
   \ref{thm:ADDR1} that all addresses in 
   $\ADDR(A)$ are contained in exactly the same wakes, and that the claim is
   true when their common internal address $a$ is finite. It thus
   suffices to consider the case where $a$ is infinite. 

  First suppose that $\ADDR(A)$ contains some unbounded infinite
   external address $\s$. We need to show that $\ADDR(A)=\{\s\}$. 
   Let $(W_i)_{i\geq 1}$ be the hyperbolic components appearing in the 
   internal address of $\s$, and set
   \[ I := \bigcap_i \cl{\W(W_i)}. \]
   Then $I$ is a closed connected subset of $\Sequ$, and 
   $\ADDR(A)\subset\{\s\}$. We claim that $I$ contains no intermediate
   external addresses (and thus consists of a single point). 

  Indeed, if $\r\in I$ was an intermediate external address, then by 
   Corollary \ref{cor:unboundeditinerary}, every entry of $\K(W_i)$ is one of the finitely
   many symbols of $\K(\r)$. However, $\K(W_i)\to \K(\s)$ by
   Corollary \ref{cor:internaladdress}, and $\K(\s)$ is unbounded. 
   This is a contradiction.

\smallskip

  Now suppose that $\ADDR(A)$ consists of more than two external
   addresses. By the previous step, all addresses in $\ADDR(A)$ are
   bounded; let $M$ be an upper bound on the size of the entries in
   their (common) kneading sequence, and let $d := 2M+2$. Then
   the map
    \[ \s \mapsto \sum_{j\geq 1} \frac{s_k}{d^j} \ (\mod 1)\]
   takes $\ADDR(A)$ injectively
   to a set $\wt{A}\subset \R/\Z$ which has the property
   that, for any wake of a hyperbolic component $W$
   \emph{in the Multibrot set $\M_d$ of degree $d$}, either
   $\wt{A}\subset \W(W)$ or $\wt{A}\cap \W(W)=\emptyset$.
   It follows from the Branch Theorem for Multibrot sets
   \cite[Theorem 9.1]{intaddr} that $\wt{A}$ (and $\ADDR(A)$) consists of 
   finitely many preperiodic addresses.

  It is easy to see that, in this case, for any $\s,\r\in\ADDR(A)$, 
   $\itin_{\s}(\r)=\K(\s)$. By the main result of \cite{expmisiurewicz}, 
   there
   exists a parameter $\kappa_0$ for which the dynamic ray $\gs$ lands at
   the singular value (and thus the singular orbit is finite
   for this parameter). It follows from \cite[Proposition 4.4]{expper} 
   that all
   rays $g_{\r}$ with $\r\in\ADDR(A)$ also land at the singular value.
   It follows easily from Hurwitz's theorem and the stability of
   orbit portraits (compare Proposition \ref{prop:orbitstability}) 
   that all parameter rays $G_{\r}$ with 
   $\r\in\ADDR(A)$ land at $\kappa_0$ (compare \cite[Theorem 5.14.5]{thesis} 
   or
   \cite[Theorem IV.6.1]{habil} for details).
   \qed

 We believe that the ideas of \cite{boundary} can be extended
  to show that
  that all hypothetical ``queer'', i.e.\ nonhyperbolic, 
  components (which conjecturally do not
  exist) must be bounded; compare the discussion 
  in \cite[Section 8]{boundary}. 
  The following corollary, which states that
  such a component could be unbounded in 
  at most two directions, is a first step in this direction. 
 
 \begin{cor}[Nonhyperbolic Components]
   \index{nonhyperbolic components}
   \index{queer components}
  Suppose that $U$ is a nonhyperbolic component in exponential parameter space
     (or more generally, any connected subset of parameter space which contains
     no attracting, indifferent or escaping parameters). Then
     $\ADDR(U)$ consists of at most two external addresses. 
 \end{cor}
 \proof This follows directly from the previous theorem except in the
  case of item (\ref{item:misiurewicz}). In the latter case, it follows
  since $U$ cannot intersect any of the parameter rays landing at the
  associated Misiurewicz points, and thus is separated from
  all but at most two of these addresses. \qed

 \subsection*{Angled Internal Addresses}

 Internal addresses do not label hyperbolic components uniquely.
 For completeness, we will now discuss a way of decorating internal addresses
 to restore uniqueness.

  \begin{defn}[Angled Internal Address]
   \label{defn:angledinternaladdresses}
  Let $\s\in\Sequ$, and 
   let $W_j$ be the components in the internal address of $W_j$. Then the 
   \emph{angled internal address} of $\adds$ is
 \[ (1,h_1)\intaddnext (n_2,h_2)\intaddnext (n_3,h_3) \intaddnext \dots, \] 
  where $W_{j+1}\subset \W(\childWH{W_j}{h_j})$.
   \index{internal address!angled}
\end{defn}

\begin{thm}[Uniqueness of Angled Internal Addresses]
  \label{thm:angleduniqueness}
 No two hyperbolic components share the same angled internal address.
\end{thm}
\sketch Suppose that $W_1\neq W_2$ have the same
  angled internal address $a$. 
  Let $M$ be an upper bound for the entries of 
  $\K(W_1)=\K(W_2)$, and set $d:= 2M+2$. Similarly as in the proof
  of Theorem \ref{thm:ADDR2}, it then follows that the Multibrot Set
  $\M_d$ contains two different hyperbolic components
  $\wt{W}_1$ and $\wt{W}_2$ which both have the same angled internal address.
  This contradicts \cite[Theorem 9.2]{intaddr}. \qedd

 \subsection*{A Combinatorial Tuning Formula}
 Let $\s$ be an intermediate external address of length
  $\geq 2$. We will give an analog of
 the concept of tuning for polynomials, on a combinatorial level. 
 For every $i$, let us denote by
  $r^i_{i-1}$ the first $n$ entries of the sector boundary 
  $\Bdy(\s,\bdyit{i}{i-1})$.

 A map $\tau:\Sequb\to \W(\s)$ is called a tuning map for $\s$, if
 $\tau(-\infty)=\s$ and
\nomenclature[tau]{$\tau$}{(tuning map)}
 \[ \tau(k\r)=\begin{cases}
             r^{\u_n+k}_{\u_n+k-1}\tau(\r) & \tau(\r) > \s \vspace{2mm} \\
             r^{\u_n+k+1}_{\u_n+k}\tau(\r) & \tau(\r) < \s \vspace{2mm} \\
             r^{\u_n+k+\onehalf}_{\u_n+k-\onehalf}\tau(\r) & \tau(\r) = \s.
           \end{cases}. \]
   \index{tuning map}
 There are exactly two such maps, which are uniquely defined by
 choosing $\tau(\per{0})$ to be either $\per{r_{\u_n}^{u_n+1}}$ or
 $\per{r_{\u_n-1}^{\u_n}}$. (Note that under a tuning map, some addresses
 which are not exponentially bounded will be mapped to addresses which
 are exponentially bounded. This is related to the fact that
 topological renormalization fails for exponential maps; see
 \cite[Section 4.3]{thesis} or \cite{nonlanding}.)

\begin{thm}[Tuning Theorem] \label{thm:tuning}
 If the internal address of $\r$ is 
  $(1,m_1)\mapsto (n_2,m_2) \mapsto (n_3,m_3) \mapsto \dots$, 
 and
 the internal address of $\s$ is 
  $(1,\wt{m}_1) \mapsto (\wt{n}_2,\wt{m}_2) \mapsto \dots\mapsto 
   (n,\infty)$,
 then the
 internal address of $\tau(\r)$ is
 \[(1,\wt{m}_1) \mapsto (\wt{n}_2),\wt{m}_2) \mapsto \dots\mapsto 
    (n,m_1') \mapsto
  (n*\wt{n}_2,m_2) \mapsto 
  (n*\wt{n}_3,m_3) \mapsto \dots,\]
 where $m_1'$ is $m_1+1$ or $m_1$, depending on whether
   $m_1\geq 0$ or $m_1 < 0$. 
\end{thm}
\sketch
 This can be easily inferred  from the 
 well-known tuning formula for Multibrot sets 
 (see e.g.\ \cite[Theorem 8.2]{jackrays} or
  \cite[Proposition 6.7]{intaddr}).

 Alternatively, it is not difficult to give a direct combinatorial
  proof of this fact; see \cite[5.11.2]{thesis}. \qedd

\section{Combinatorial Algorithms} 
 \label{app:algorithms}
 In this article, we have seen several ways to describe hyperbolic
  components, and our results allow us to compute any of these from any
  other. Since these algorithms have not always been made explicit
  in the previous treatment, we collect them here. 
%
%

First note the kneading sequence of a hyperbolic component 
 $W$ can be easily computed from its intermediate external address 
 according to the definition (Definition \ref{defn:combinatorialitinerary}). 
 How to compute
 the internal address of a hyperbolic component from its kneading sequence
 was shown in Corollary \ref{cor:internaladdress}, which also
 describes how to obtain the forbidden kneading sequence $\KS(W)$ from
 $\K(W)$. 
 Furthermore, $\KS(W)$ can easily be computed from its
 characteristic ray pair $\langle \s^-, \s^+ \rangle$. 

\begin{algorithm}[Computing $\extaddr(W)$ given $\langle \s^-,\s^+\rangle$]
 \label{alg:intermediatefromcharacteristic}
 \emph{Given:} Kneading sequence 
   $\K(W)=\u_1\dots \u_{n-1} *$ and some $\r\in \cl{\W(W)}$. 

 \smallskip
 \noindent
 \emph{Aim:} Compute $\extaddr(W)$. 

 \smallskip
 \noindent
 \emph{Algorithm:} Let $\s^n := \infty$, and compute
   $\s^{n-1},\dots, \s^{1}$ inductively by choosing
   $\s^{j-1}$ to be the unique preimage of $\s^j$ in
   $\bigl(\u_{j-1}\r,(\u_{j-1}+1)\r\bigr)$. Then
   $\extaddr(W)=\s^1$.
\end{algorithm}
\proof By Observation \ref{obs:itinchange} and 
   Lemma \ref{lem:itineraries}, $\extaddr(W)$ is the
   unique address $\s$ with $\itin_{\r}(\s)=\K(W)$. This is
   exactly the address $\s^1$ computed by this algorithm. \qed

\begin{algorithm}[Computing $\langle \s^-,\s^+\rangle$ given
                  $\extaddr(W)$] 
   \label{alg:characteristicfromintermediate}
 \emph{Given:} Kneading sequence
   $\K(W)=\u_1\dots \u_{n-1}*$, some $\r\in \cl{\W(W)}$ and some
    $s_*\in\Z$. 

 \smallskip
 \noindent
 \emph{Aim:} Compute $\Bdy\bigl(W,\bdyit{s_*}{s_*-1}\bigr)$.

 \smallskip
 \noindent
 \emph{Algorithm: } Compute the unique preimage $\wt{\r}$ of $s_*\r$
  whose itinerary (with respect to $\r$) begins with
  $\u_1\dots \u_{n-1}$, as in Algorithm
  \ref{alg:intermediatefromcharacteristic}. 
  The sought address is obtained by continuing
  the first $n$ entries of $\wt{\r}$ periodically.
\end{algorithm}
\begin{remark}
 To actually compute the characteristic addresses of $W$,
  first compute the forbidden kneading entry $\u_n = \u(W)$, and then
  apply the algorithm to $s_*=\u_n$ and $s_*=\u_n+1$.
\end{remark}
\proof Using Observation \ref{obs:agreeingaddresses}
  (
  similarly to Proposition \ref{prop:itinerarymonotonicity}
 (\ref{item:itinerarymonotonicity}), one shows that
  the interval between $\wt{\r}$ and the required address is mapped
  bijectively by $\sigma^n$. \qed

\begin{algorithm}[Compute the Angled Internal Address of $W$]
 \emph{Given:} An intermediate external address $\s$.

 \smallskip
 \noindent
 \emph{Aim:} Compute the angled internal address of $\s$. 

 \smallskip
 \noindent
 \emph{Algorithm:} Let $n$ be the length
  of $\s$, and calculate the kneading sequence 
   $\addu=\u_1\dots \u_{n-1} *:= \K(\s)$. Set $n_1 := 1$
   and $m_1 := \u_1$. 

  Given, for some $j\geq 1$, two numbers $n_j<n$ and 
   $p^j\in\Z$, we calculate three numbers $n_{j+1}>n_j$, 
   $m_{j+1}\in\Z$ and
   $h_j\in \Q$. The algorithm terminates when $n_{j+1}$ is equal to $n$;
   at this point the angled internal address of $\s$ will be
   given by
   \[ (n_1,h_1) \mapsto (n_2,h_2) 
      \mapsto \dots \mapsto (n_{j},h_j)\mapsto (n,\infty). \]

 \emph{Step 1: Calculation of $n_{j+1}$ and the corresponding
  sector number.}
  Define $n_{j+1}$ to be the first index at which
   $\addu$ and $\addu^j := \per{\u_1\dots \u_{n_j}}$ differ.
   If $n_{j+1}\neq n$, then set 
   $m_{j+1} := \u_{n_{j+1}} - \u^j_{n_{j+1}}$.

 \emph{Step 2: Determining the denominator.} 
  Let us inductively
   define a finite sequence $\ell_1 \leq \ell_2 \leq \dots \leq 
   \ell_r$ as follows. Set $\ell_1 := n_{j+1}$; 
   if $n_j$ does not divide $\ell_k$, let
   $\ell_{k+1}$ be the first index at which
   $\addu^j$ and $\per{\u_1\dots \u_{\ell_k}}$ differ. Otherwise,
   we terminate, setting $r:= k$. 
   Set $q_j := \ell_r / n_j$. 

 \emph{Step 3: Determining the numerator.} 
  If $q_j=2$, set $p_j := 1$. Otherwise, calculate
    \[ x := \#\{k\in\{2,\dots,q_j-2\}:
                \sigma^{kn_j}(\s) \text{ lies between }
                \s \text{ and } \sigma(\s)\}. \]
  Set $p_n := x+1$ if $\sigma(\s)>\s$, and
   $p_n := n_j - x - 1$ otherwise.

 \emph{Step 4: Determining $h_j$.}
  $h_j$ is defined to be $m_j+p_n/q_n$ if
  $j=1$ or $m_j < 0$, and $m_j-1 + p_n/q_n$ otherwise. 
\end{algorithm}
\proof Let $W_j$ denote the $j$-th component in the internal 
  address of $\s$. By Corollary \ref{cor:internaladdress},
  the value of $n_j$ computed by our algorithm will be the
  period of $W_j$, and if $n_j\neq n$, then 
  $m_j$ is the number of the sector containing $\s$.

  If $n_j\neq n$, let $V_j$ be the child component of $W$ containing
   $\s$. By Theorem \ref{thm:determiningbifurcationorder},
   this component has period $\ell_r= q_j n_j$. 
   By Proposition \ref{prop:itinerarymonotonicity}
   (\ref{item:itinerarymonotonicity}), the interval
   $\W(V_j)$ is mapped bijectively by $\sigma^{(q_j-2)n_j}$.
   (Note that this shows, in particular, that
    $\ell_r < 2n_j+n$, and thus
    limits the iterations necessary in Step 2.)
   Therefore, the order of the iterates of $\s$ and those of
    $\extaddr(V_j)$ is the same. It follows that
    $V_j=\childWH{W_j}{h_j}$, where $h_j$ is defined
    as indicated.  \qed
\begin{remark}
 In Step 3, we could have instead first calculated the 
  intermediate external address of the bifurcating component
  in question (using Algorithm \ref{alg:bifurcation} below), and then 
  calculated its rotation number. This introduces an
  extra step in the algorithm, but makes the proof somewhat
  simpler. 
\end{remark}

It remains to indicate how an intermediate external address can be
 recovered from its angled internal address. Let us first note how to
 handle the special case of computing a combinatorial bifurcation.

\begin{algorithm}[Computing Combinatorial Bifurcations]
  \label{alg:bifurcation}
 \emph{Given:} $\s:= \extaddr(W)$, $s_*\in\Z$ and
   $\alpha=\frac{p}{q}\in (0,1)\cap\Q$.

 \smallskip
 \noindent
 \emph{Aim:} Compute $\bifaddr(W,s_*,\alpha)$.

 \smallskip
 \noindent
 \emph{Algorithm: } Compute the unique address whose itinerary under $\s$
   is as given by Lemma \ref{lem:itinerarybifurcation}. 
   (Alternatively, compute the
   sector boundaries of the given sector and apply the Combinatorial
   Tuning Formula.) \qed
\end{algorithm} 

\begin{algorithm}[Computing $\extaddr(W)$ from an angled internal address]
 \emph{Given:} The angled internal address
    $a=(n_1,h_1)\mapsto \dots \mapsto (n,\infty)$ of some hyperbolic
    component $W$.

 \smallskip
 \noindent
 \emph{Aim:} Compute the unique intermediate external address $\s$ whose
    angled internal address is $a$.

 \smallskip
 \noindent
 \emph{Algorithm:} Let $W_i$ denote the component represented by the
   $i$-th entry of $a$. We will compute $\s^i := \extaddr(W_i)$ inductively
   as follows;
   note that $\s^1 = \infty$. 

  Compute the upper characteristic address $\r^+$ of the component 
   $V_i := \childWH{W_i,h_i}$
   (by applying Algorithms \ref{alg:bifurcation} 
    and \ref{alg:characteristicfromintermediate}, or by using the
    combinatorial tuning formula). 

  Let $\t^1\leq \r^+$ be maximal such that $\t^1$ is periodic of
   period at most $n_{i+1}$. If the period of $\t^1$ is strictly less than
   $n_{i+1}$, then $\t^1$ is the upper characteristic address of some
   hyperbolic component. Compute the lower characteristic address $\t^{1-}$
   of this component and let 
   $\t^2\leq \t^{1-}$ be maximal such that $\t^2$ is periodic of period
   at most $n_{i+1}$.

  Continue until an address $\t^k$ is computed which is
   periodic of period $n_{i+1}$. This is the upper characteristic 
   addresses of $\s^{i+1}$, and we can compute $\s^i$ using Algorithm 
    \ref{alg:intermediatefromcharacteristic}.
\end{algorithm}
\proof Recall that $\s^{i+1}$ is not contained in the wake
  of any component $U\prec V_i$
   which has smaller period than $n_{i+1}$, and by the
   uniqueness of angled internal addresses
   (Theorem \ref{thm:angleduniqueness}), there are only finitely many
   periodic addresses $\t^j$ which are encountered in each step.
   Thus, the
   algorithm will indeed terminate and
   compute an address which is periodic of period
   $n_{i+1}$. The associated hyperbolic component
   then has angled internal address
    \[ (n_1,h_1) \mapsto \dots \mapsto (n_j,h_j) \mapsto 
        (n_{j+1},\infty) \]
   and the claim follows by Theorem \ref{thm:angleduniqueness}. \qed

 Finally, let us mention that it is very simple to decide whether
  a given hyperbolic component is a child component, and to calculate
  the address of the parent component in this case. In particular,
  it is easy to decide whether two given hyperbolic components 
  have a common parabolic boundary point (or, in fact, \emph{any}
  common boundary point, compare \cite[Proposition 8.1]{boundary}).

  \begin{algorithm}[Primitive and Satellite Components]
   Let $\s\in\Sequ$ be an intermediate external address with
   kneading sequence $\u_1\dots \u_{n-1} *$. If there exists some
   $\u_n\in \Z$ such that $\per{\u_1\dots \u_{n-1} \u_n}$ is periodic
   with period $j<n$, then $\Hyp{\s}$ is a child component
   of $\Hyp{\sigma^{n-j}(\s)}$.

   Otherwise, $\Hyp{\s}$ is a primitive component.
  \end{algorithm}
  \proof This is an immediate consequence of
    Corollary \ref{cor:parabolicbifurcation} 
    and Theorem \ref{thm:bifurcationstructure}. \qed

\newpage

{\def\nompreamble{\footnotesize}
\printnomenclature[3.5cm]}

\printindex


\small
\bibliographystyle{hamsplain}
\small{\bibliography{c:/Latex/Biblio/biblio}}

\end{document}